\documentclass[12pt]{article}
\usepackage{cite,epic,amsmath,amsthm,amssymb,texdraw}

\newtheorem{thm}{Theorem}[section]
\newtheorem{prop}[thm]{Proposition}
\newtheorem{lem}[thm]{Lemma}
\newtheorem{dfn}[thm]{Definition}
\newtheorem{coro}[thm]{Corollary}
\newcommand{\nc}{\newcommand}

\font\germ=eufm10
\nc{\goth}[1]{{\germ #1}}
\nc{\qbinom}[2]{{\genfrac{[}{]}{0pt}{}{{#1}}{{#2}}}_{q}}
\nc{\bra}[1]{\langle #1 |}
\nc{\ket}[1]{| #1 \rangle}
\nc{\qpf}[1]{({#1}\, ; \, q^4)_{\infty}}
\nc{\br}[1]{\begin{array}{#1}}
\nc{\Ps}[2]{\Phi^{(#1,#2)}}
\nc{\Psz}[2]{\Phi^{(#1,#2)}(\zeta)}
\nc{\Vz}[1]{V^{(#1)}_{\zeta}}
\nc{\wo}[6]{W^I_{#6}\left(\left.
            \begin{array}{ll}{#1}&{#2}\\{#3}&{#4}\end{array}
            \right|{#5}\right)}
\nc{\wt}[6]{W^{II}_{#6}\left(\left.
            \begin{array}{ll}{#1}&{#2}\\{#3}&{#4}\end{array}
            \right|{#5}\right)}
\nc{\wth}[6]{W^*_{#6}\left(\left.
            \begin{array}{ll}{#1}&{#2}\\{#3}&{#4}\end{array}
            \right|{#5}\right)}
\nc{\wec}[5]{\overline{W}^1_{\ell}\left(\left.
            \begin{array}{ll}{#1}&{#2}\\{#3}&{#4}\end{array}
            \right|{#5}\right)}
\nc{\isomo}{\buildrel {\sim} \over \longrightarrow}
\nc{\Czp}{|_0^+}
\nc{\Czm}{|_0^-}
\nc{\Cop}{|_1^+}
\nc{\Com}{|_1^-}
\nc{\Bzp}{\lfloor_0^+}
\nc{\Bzm}{\lfloor_0^-}
\nc{\Bop}{\lfloor_1^+}
\nc{\Bom}{\lfloor_1^-}
\nc{\Tzp}{\lceil_0^+}
\nc{\Tzm}{\lceil_0^-}
\nc{\Top}{\lceil_1^+}
\nc{\Tom}{\lceil_1^-}
\nc{\bpp}{\bullet^{++}}
\nc{\bpm}{\bullet^{+-}}
\nc{\bmp}{\bullet^{-+}}
\nc{\bmm}{\bullet^{--}}
\nc{\bpmpm}{\bullet^{\pm\pm}}
\nc{\bpmmp}{\bullet^{\pm\mp}}
\nc{\Tb}{\lceil}
\nc{\Bb}{\lfloor}
\nc{\ttn}{\kappa}
\nc{\ttm}{\kappa'}
\nc{\tts}{k}
\nc{\ssa}{s}
\nc{\ssb}{t}
\nc{\ssj}{\varepsilon}
\nc{\Aff}{\operatorname{Aff}}
\nc{\Affz}{\operatorname{Aff}(B^{(0)})}
\nc{\Affo}{\operatorname{Aff}(B^{(1)})}
\nc{\calA}{\mathcal{A}}
\nc{\calP}{\mathcal{P}}
\nc{\calR}{\mathcal{R}}
\nc{\calN}{\mathcal{S}}
\nc{\calNo}{\mathcal{N}}
\nc{\calD}{\mathcal{D}}
\nc{\ot}{\otimes}
\nc{\ve}{\varepsilon}
\nc{\be}{\begin{eqnarray}}
\nc{\ee}{\end{eqnarray}}
\nc{\bea}{\begin{eqnarray}}  
\nc{\ena}{\end{eqnarray}}
\nc{\C}{\mathbf{C}}
\nc{\Q}{\mathbf{Q}}
\nc{\R}{\textit{R}}
\nc{\half}{\ensuremath{\frac{1}{2}}}
\nc{\z}{\zeta}
\nc{\La}{\Lambda}
\nc{\s}{\sigma}
\nc{\Hom}{\operatorname{Hom}}
\nc{\End}{\operatorname{End}}
\nc{\vac}{|\textrm{vac}\rangle}
\nc{\dvac}{\langle\textrm{vac}|}
\nc{\id}{\operatorname{id}}
\nc{\la}{\lambda}  
\nc{\bqa}{\begin{eqnarray}}  
\nc{\eqa}{\end{eqnarray}}  
\nc{\ra}{\rightarrow}  
\nc{\lra}{\longrightarrow}
\nc{\uqp}{U^{\prime}_q (\widehat{sl}_2)}
\nc{\uq}{U_q(\slth)}
\nc{\vsl}{V(\sigma(\lambda))}
\nc{\Rbar}{\bar{R}}
\nc{\vl}{V(\lambda)}  
\nc{\bu}{\bullet}
\nc{\an}{{\ell}}
\nc{\ak}{{\ell}}
\nc{\am}{{\ell'}}
\nc{\al}{{\ell'}}
\nc{\ep}{\varepsilon} 
\nc{\de}{\delta}  
\nc{\slth}{\widehat{\goth{sl}}_2\hskip 1pt}
\nc{\ws}{\;\;}
\nc{\qu}{{1\ov 4}}
\nc{\tf}{\tilde{f}}
\nc{\te}{\tilde{e}}
\nc{\ez}{|_{0}}
\nc{\eo}{|_{1}}
\nc{\cz}{[_{0}}
\nc{\co}{[_{1}}
\nc{\hif}{\hb{ if }}
\nc{\hev}{\hb{ is even }}
\nc{\hod}{\hb{ is odd }}
\nc{\HR}{{_h}\bar{R}}
\nc{\PR}{\bar{R}}
\nc{\cP}{\mathcal{P}}
\nc{\er}{\end{array}}
\nc{\Tr}{{\rm Tr}}
\nc{\hb}{\hbox}
\nc{\cT}{\mathcal{T}}
\nc{\zi}{\zeta^{-1}}
\nc{\nn}{\nonumber} 
\nc{\Z}{\mathbf{Z}}
\nc{\bR}{\bar{R}}
\nc{\curlra}{\buildrel{\sim}\over\longrightarrow}
\nc{\bs}{\tilde{s}}
\nc{\cD}{\mathcal{D}}  
\nc{\cH}{\mathcal{H}}
\nc{\cF}{\mathcal{F}}
\nc{\cO}{\mathcal{O}}
\nc{\epp}{\varepsilon^{\prime}} 
\nc{\bi}{\bar{i}}
\nc{\bj}{\bar{j}}
\nc{\ol}{\overline}
\nc{\Om}{\Omega}
\nc{\pl}{\prod\limits} 
\nc{\sli}{\sum\limits} 

\begin{document}
\begin{flushright}
 q-alg/9811175 \\[10mm]
\end{flushright}
\begin{center}
{\Large \bf Vertex Models with Alternating Spins\\[8mm] }
{Jin Hong\,$^1$, Seok-Jin Kang\,$^1$,
Tetsuji Miwa\,$^2$ and Robert Weston\,$^3$\\[8mm]
{\it Dedicated to Professor Mikio Sato on his seventieth birthday.}}
\end{center}
\footnotetext[1]{Department of Mathematics, Seoul National
University, Seoul 151-742, Korea.}
\footnotetext[2]{Research Institute for Mathematical Sciences,
Kyoto University, Kyoto 606, Japan.}
\footnotetext[3]{Department of Mathematics, Heriot-Watt University,
Edinburgh EH14 4AS, UK. }

\begin{abstract}
\noindent 
The diagonalisation of the transfer matrices of solvable
vertex models with alternating spins is given. The crystal
structure of (semi-)infinite tensor products of finite dimensional
$U_q(\widehat{sl}_2)$ crystals with alternating dimensions
is determined. Upon this basis the vertex models are formulated and
then solved by means of $U_q(\widehat{sl}_2)$ intertwiners.
\end{abstract}
\nopagebreak

\setcounter{equation}{0}
\section{Introduction}

In~\cite{DFJMN}, the diagonalisation of the XXZ Hamiltonian,
\bea
H_{XXZ}&=&-\half\sum_{k=-\infty}^{\infty}
\Bigl(\sigma^x_{k+1}\sigma^x_k+\sigma^y_{k+1}\sigma^y_k+
\Delta\sigma^z_{k+1}\sigma^z_k\Bigr),
\ena
in the anti-ferromagnetic regime $(\Delta=\frac{q+q^{-1}}{2}<-1)$ was carried
out by making use of the representation theory of the quantum affine algebra
$\uq$. The key observation in this method was the identification
of the semi-infinite tensor product of the two-dimensional representation
$V^{(1)}\simeq \C^2$ of $\uq$ with the level one irreducible highest
weight representation $V(\La_i)$ $(i=0,1)$ of the same algebra~\cite{FM},
\bea\label{IDEN}
\cdots \ot \C^2\ot\C^2\ot\C^2\ot\C^2\simeq V(\La_i).
\ena
Using~\eqref{IDEN}, the corner transfer matrix $A(\z)$ of the corresponding
six-vertex model was identified with the grading operator
\bea\label{CTM}
A(\z)&\sim&\z^{-D},
\ena
and the half transfer matrix $\Phi(\z)$ was identified with the vertex operator
\bea\label{HTM}
\Phi(\z):V(\La_i)\ra V(\La_{1-i})\ot V^{(1)}_\z,
\ena
where $V^{(1)}_\z$ is the evaluation representation.
The choice of $i=0,1$ corresponds to the choice of the boundary condition
at infinity.

Under these identifications, the transfer matrix $T(\z)$
was identified with the composition of the vertex operators
acting on the tensor product of the highest and lowest weight representations,
\bea
T(\z)&:&V(\La_i)\ot V(\La_j)^*
\ra V(\La_{1-i})\ot V^{(1)}_\z\ot V(\La_j)^*
\ra V(\La_{1-i})\ot V(\La_{1-j})^*,
\ena
and then diagonalised by making use of another vertex
operator~\cite{DFJMN, JM}
\bea
\Psi^*(\xi)&:&V^{(1)}_\xi\ot V(\La_j)\ra V(\La_{1-j}).
\ena

A similar method was also applied to other models such as
the higher spin generalisation of the XXZ model~\cite{idzal93} and the ABF
models~\cite{JMO93}.
In the former, for which the local spaces are $V^{(n)}\simeq\C^{n+1}$,
the spaces of physical states in the semi-infinite volume
with the chosen boundary conditions
were identified with the level~$n$ irreducible highest weight representations.
On the other hand, in the latter, they were identified
with the coset spaces of GKO type (see also~\cite{LP96}).

In this paper, we study yet another example of this sort. We consider
the vertex models with alternating spins. This requires new insights,
both physical and mathematical, and leads to new results
in the connection between solvable lattice models and  representation
theory.

The origin of our study is~\cite{Nak}, in which a spin-$\half$ chain
with a few higher spin components (or \emph{impurities} in physical terms)
was studied by using the vertex operator
\bea
\Phi^{(n-1,n)}(\z):V^{(n-1)}_\z\ot V(\La_i) \ra V(\La_{1-i}) \ot V^{(n)}_\z.
\ena
This operator explains the $n$-fold degeneracy of the vacuum states
with a chosen boundary condition
when a spin-$\frac{n}{2}$ impurity is inserted in the spin-$\half$ chain.
In~\cite{MW97}, the above vertex operator was identified with the half
transfer matrix of the vertex model that has semi-infinite spin-$\half$
horizontal lines and a spin-$\frac{n}{2}$ vertical line.
In this paper we consider a vertex model with alternating spins $\frac{m}{2}$
and $\frac{n}{2}$ $(m>n)$, and diagonalise the corresponding transfer matrices.
Such models were constructed and analysed using the Bethe
Ansatz in \cite{VW92,AM93,VMN93,VMN94}.
The first step in our
solution is the identification of the semi-infinite tensor product
\bea
\cdots \ot \C^{m+1}\ot\C^{n+1}\ot\C^{m+1}\ot\C^{n+1}\label{SEMIPATH}
\ena
having an appropriate boundary condition, with the tensor product of
level~${m-n}$ and level~$n$ highest weight representations
\bea
V(\la^{(m-n)}_a)\ot V(\la^{(n)}_b).\label{TENSORPRO}
\ena
Here we set
\bea
\la^{(\ell)}_a=(\ell-a)\La_0+a\La_1.
\ena
Formula~\eqref{CTM} is again valid in this situation.

In the second step, we identify the half transfer matrices
(see Figure 3) having alternating spins for the horizontal lines
and spin-$\frac{n}{2}$ (Case A) or $\frac{m}{2}$ (Case B) for the vertical
line, with the following vertex operators.

\noindent
Case A:
\bea
\phi^A(\z)&:&
      V(\la^{(m-n)}_a)\ot V(\la^{(n)}_b)
\xrightarrow{\id \ot \Phi^{(0,n)}(\z)}
         V(\la^{(m-n)}_a)
         \ot V(\la^{(n)}_{n-b})
         \ot V^{(n)}_\z.
\ena

\noindent
Case B:
\bea
\phi^B(\z)&:&
      V(\la^{(m-n)}_a) \ot V(\la^{(n)}_b)
\xrightarrow{\Phi^{(0,m-n)}(\z)\ot \id}
           V(\la^{(m-n)}_{m-n-a})
           \ot V^{(m-n)}_\z
           \ot V(\la^{(n)}_b)\nonumber\\
& &\xrightarrow{\id\ot\Phi^{(m-n,m)}(\z)}
           V(\la^{(m-n)}_{m-n-a})
           \ot V(\la^{(n)}_{n-b})
           \ot V^{(m)}_\z,
\ena

Finally, we have two (full) transfer matrices $T^A(\z)$ and $T^B(\z)$
for cases A and B.
We can think of these operators as acting on the direct sum of the
vectors spaces,
\bea
&&\Hom_\C\Bigl(V(\la^{(m-n)}_a)\ot V(\la^{(n)}_b),
V(\la^{(m-n)}_c)\ot V(\la^{(n)}_d)\Bigr)\nonumber\\
&&\simeq
V(\la^{(m-n)}_c)\ot V(\la^{(n)}_d)
\ot
\Bigl(V(\la^{(m-n)}_a)\ot V(\la^{(n)}_b)\Bigr)^*.\label{PHYS}
\ena
$T^A(\z)$ and $T^B(\z)$ are mutually commuting and expressed in terms
of the operators $\phi^A(\z)$ and $\phi^B(\z)$, respectively.

The operator $T^A(\z)$ can be viewed as the limit where the number of
insertions of the higher spin components becomes infinite.
Therefore, in this case,
one can expect that the vacuum states are infinitely degenerate,
and it is indeed so. The same is true for $T^B(\z)$. However, if we consider
the product $T(\z)=T^B(\z)T^A(\z)$ the infinite degeneracy resolves,
and we have a unique vacuum for a fixed boundary condition.
The vacuum states are given by
\bea
(-q)^D\in\End_\C\Bigl(V(\la^{(m-n)}_a)\ot V(\la^{(n)}_b)\Bigr).
\ena

The excited states are constructed upon these vacua.
Consider two kinds of vertex operators with spin $0$ and $\half$,
respectively.

\noindent
Spin-$0$ case:
\bea
\psi^{(0)}(\xi)&:&
V(\la^{(m-n)}_a)\ot V(\la^{(n)}_b)
\ra
V(\la^{(m-n)}_{a'})\ot V^{(1)}_\xi\ot V(\la^{(n)}_b)
\ra
V(\la^{(m-n)}_{a'})\ot V(\la^{(n)}_{b'}).\nonumber\\
\ena
Spin-$\half$ case:
\bea
\psi^{(\half)}(\xi)&:&
V^{(1)}_\xi\ot V(\la^{(m-n)}_a)\ot V(\la^{(n)}_b)
\ra
V(\la^{(m-n)}_{a'})\ot V(\la^{(n)}_b).
\ena
Acting on the vacua,  the operators $\psi^{(0)}(\xi)$
and $\psi^{(\half)}(\xi)$ create particles with spin $0$ and $\half$,
respectively.
We give the exchange relations for these operators. The vacuum states
$(-q)^D$, the operators $\psi^{(0)}(\xi)$ and $\psi^{(\half)}(\xi)$
and their exchange relations are the diagonalisation data of the transfer
matrix $T(\z)$ in the sense of the vertex operator approach~\cite{DFJMN}.
 From the view point of the representation theory this data gives the
irreducible decomposition of the space of physical states~\eqref{PHYS}
with respect to the action of $U_q(\widehat{sl}_2)$.
We call this description of the physical space the \emph{particle} picture
in comparison with the \emph{local} picture consisting of the alternating
infinite tensor product
of $\C^{m+1}$ and $\C^{n+1}$.
We should say that the equivalence of the local and
particle pictures is a conjecture because we have no argument for the
completeness of the particle decomposition except in the crystal limit
$q=0$ (see (ii) below).

Many of the results in this paper have been announced in~\cite{HKMW98a}.
In this paper we give proofs for them. 
(On the other hand, we will not discuss the mixing of ground states,
one of the main results in~\cite{HKMW98a}.
We have
nothing to add to the result and a complete proof is already given
there.)
To be precise, we prove the following:

(i) A crystal isomorphism between the space of semi-infinite paths $P_{a,b}$
and the crystal $B(\la^{(m-n)}_a)\ot B(\la^{(n)}_b)$.

The crystal structure of $P_{a,b}$ represents  by definition
the semi-infinite tensor product of the alternating finite crystals
$B^{(m)}$ and $B^{(n)}$. Therefore, the crystal isomorphism
mentioned above gives supporting evidence for the conjecture that
there is an isomorphism between~\eqref{SEMIPATH} and~\eqref{TENSORPRO}.
We give two proofs.
The first one uses the RSOS paths which describe the highest weight vectors
in $B(\la^{(m-n)}_a)\ot B(\la^{(n)}_b)$.
The second proof is more direct; however, the identification of the corner
transfer matrix is made only in the first proof.

(ii) Crystal decomposition of the full-infinite path spaces.

We decompose each path uniquely to a union of ground state paths patched
together at the `walls' between the ground states. Under the crystal action
these walls behave like elements in the affinizations~\cite{KMN92} of
$B^{(0)}$ or $B^{(1)}$.
This observation gives supporting evidence for our conjecture that
the particle structure of the alternating vertex model consists of
the spin-$0$ and spin-$\half$ particles.

(iii) Commutativity of the vertex operator
\begin{equation}
\begin{aligned}
V^{(l)}_\z\ot V(\la^{(k)}_a)\ot V(\La_i)&\ra
V(\la^{(k)}_{k-a})\ot V^{(l+k)}_\z\ot V(\La_i)\\
&\ra V(\la^{(k)}_{k-a})\ot V(\La_{1-i})\ot V^{(l+k+1)}_\z
\end{aligned}
\end{equation}
with the DVA (deformed Virasoro algebra) actions~\cite{JS97} on
$V(\la^{(k)}_a)\ot V(\la_i)$ and $V(\la^{(k)}_{k-a})\ot V(\la_{1-i})$.
This fact is used to derive the properties of the vertex operators
of higher level from those of level $1$.

The plan of the paper is as follows.
In Section~\ref{sec2}, the vertex models with alternating spins are
formulated.
The ground states and the eigenvalues of the corner transfer matrices
are determined.
In Section~\ref{sec3}, the path space, i.e., the $q\ra0$ limit of the
model, is studied.
In Section~\ref{sec4}, we prepare some properties of the level-$1$
vertex operators.
In Section~\ref{sec5}, the commutativity with the DVA is proved.
The diagonalisation of the transfer matrices is discussed in
Section~\ref{sec6}.
In Section~\ref{sec7} we give the crystal isomorphism between the local and
particle pictures. Finally, we present a brief summary of our results in
Sectiion~\ref{sec8}.

\emph{Acknowledgements.}
We thank Boris Feigin, Michio Jimbo and Masaki Kashiwara for discussions.
Two of us (TM and RW) thank Feodor Smirnov for his hospitality in Paris.
TM thanks everyone in GARC at Seoul National University for their
hospitality.
JH, SJK, and RW thank the people of RIMS at Kyoto University for their
hospitality.
This research was partially supported by Basic Science Research
Institute Program, Ministry of Education of Korea, BSRI-98-1414 and
GARC-KOSEF at Seoul National University, and also by
the Grant-in-Aid for Scientific Research A-08304001
from the Ministry of Education, Science and Culture of Japan.
RW acknowledges the support of the EPSRC through Advanced
Fellowship B/96/AF/2235; he also acknowledges funding from the RIMS/Isaac
Newton Inst./JSPS/Roy. Soc. Exchange Program.


\setcounter{equation}{0}
\section{The Vertex Model}\label{sec2}

In this section, we recall the definition of the alternating spin vertex model
of reference~\cite{HKMW98a}.
We define the path space, the corner transfer matrices, and the corner 
transfer matrix Hamiltonian.

\subsection{The \R-matrices}\label{sec21}
The Boltzmann weights of our vertex model are given in terms of $\uqp$
\R-matrices (as usual $\uqp$ refers to the subalgebra of $\uq$
generated by $e_i,f_i,t_i$ $(i=0,1)$; our comultiplication
is that of~\cite{MW97}).
We use the spin-$\frac{n}{2}$ principal $\uqp$ evaluation module
$V^{(n)}_{\z}$
defined, in terms of weight vectors $u^{(n)}_{i}$ $(i=0,1,\cdots,n)$, in
Section 3.1 of~\cite{MW97}.

In this paper, we consider the spectral parameter $\z$ (or $z=\z^2$)
mainly as a generic complex number, and $U_q(\widehat{sl}_2)$
as a $\C$-algebra.
However, in Sections~\ref{sec4} and~\ref{sec5}, when we develop the theory
of vertex operators, we treat the spectral parameter as an auxiliary variable.
Namely, when we consider the evaluation module $V^{(n)}_\z$,
we always extend the field of coefficients to a ring by adding $\z$ and
$\z^{-1}$. Therefore, we consider $V^{(n)}_\z$ as the rank $n+1$
$U'_q(\widehat{sl}_2)$ module over the extended ring.

The necessary \R-matrix 
is given by the $\uqp$ intertwiner
$R^{(\ak,\al)}(\z_1/\z_2):V^{(\ak)}_{\z_1}\ot V^{(\al)}_{\z_2}\ra
 V^{(\al)}_{\z_2}\ot V^{(\ak)}_{\z_1}$
(note that $R^{(\ak,\al)}(\z)$ here is $P R^{(\ak,\al)}(\z)$
in the notation of~\cite{MW97}).
We fix the normalisation by the requirement 
$R^{(\ak,\al)}(\z) = \bR^{(\ak,\al)}(\z)/\kappa^{(\ak,\al)}(\z)$, where
$\bar{R}^{(\ak,\al)}(\z)(u^{(\ak)}_0 \ot u^{(\al)}_0)
=(u^{(\al)}_0 \ot u^{(\ak)}_0)$,
and
\be
\kappa^{(\ak,\al)}(\z)=\z^{\min(\ak,\al)}
 \frac{\qpf{q^{2+\ak+\al}\z^2}\qpf{q^{2+|\ak-\al|}\z^{-2}}}%
      {\qpf{q^{2+\ak+\al}\z^{-2}} \qpf{q^{2+|\ak-\al|}\z^2} }.\nn
\ee

This choice of  normalisation has two nice consequences. The first
is that the partition function per site of our vertex model is equal to 1.
The second is that the \R-matrix has the properties of crossing
symmetry and unitarity:
\bea
R^{(\ak,\al)}(\z)^{i,j}_{i',j'} &=&   
R^{(\al,\ak)}(-q^{-1}\zi)^{l'-j',i}_{l'-j,i'},\label{crossing}\\
\sli_{i',j'} R^{(\ak,\al)}(\z)^{i',j'}_{i_1,j_1} 
R^{(\al,\ak)}(\zi)_{j',i'}^{j_2,i_2} &=& \de_{i_1,i_2}
 \de_{j_1,j_2}.\label{unit}\ena
Here we use the components defined by 
\be R^{(\ak,\al)}(\z)(u^{(\ak)}_i \ot u^{(\al)}_j)=\sli_{i',j'}
R^{(\ak,\al)}(\z)^{i,j}_{i',j'}(u^{(\al)}_{j'}\ot u^{(\ak)}_{i'}).\nn
\ee

We wish to give an expansion of $\PR^{(\ak,\al)}(\z)$ in terms of certain
projectors.
In order to do this it is useful to introduce a homogeneous evaluation
module $(V_n)_z$ with weight vectors $v^{(n)}_{i}$ $(i=0,1,\cdots,n)$.
The action of $\uqp$ on $(V_n)_z$ is given by
\be
&& f_1 v_j^{(n)} = [n-j]  v_{j+1}^{(n)}, \quad
e_1 v_j^{(n)} = [j]  v_{j-1}^{(n)},\quad t_1 v_j^{(n)}=q^{n-2j}
v_j^{(n)},\\
&& f_0= z^{-1} e_1,\quad e_0= z f_1, \quad t_0 =t_1^{-1}.
\ee
We shall refer to the associated $U_1=\langle e_1,f_1,t_1\rangle$
module as $V_n$.
The $\uqp$-modules
$(V_n)_z$ and $V^{(n)}_{\z}$ are isomorphic. The isomorphism is given by
\be C_n(\z): V^{(n)}_{\z} &\curlra& (V_n)_z,\\
u_j^{(n)} &\longmapsto& c_j^{(n)} \z^{j} v_j^{(n)},
\ee
where $c_j^{(n)}={{\qbinom{n}{j}}}^{\!\!\!\!\half} q^{\frac{j}{2}(n-j)}$ and
we identify $\z^2=z$\; (in this paper, we shall use the notation
$[a]_q=(q^a-q^{-a})/(q-q^{-1})$, and $[a]_q!$ and $\qbinom{a}{b}$ for 
the standard $q$-factorial and $q$-binomial
coefficients).
Consider the $\uqp$ intertwiner
$\HR^{(\ak,\al)}(z_1/z_2): (V_{\ak})_{z_1} \ot (V_{\al})_{z_2}
\lra (V_{\al})_{z_2}\ot (V_{\ak})_{z_1}$ defined uniquely by
\be
v_0^{(\ak)} \ot v_0^{(\al)} \longmapsto v_0^{(\al)} \ot v_0^{(\ak)}.\ee
The \R-matrix $\PR^{(\ak,\al)}(\z_1/\z_2)$ is given in terms of this
intertwiner by
\be
\PR^{(\ak,\al)} (\z_1/\z_2)
 = \left( C_{\al}(\z_2)^{-1}\ot C_{\ak}(\z_1)^{-1}\right) 
\HR^{(\ak,\al)}((\z_1/\z_2)^2) \left( C_{\ak}(\z_1)\ot C_{\al}(\z_2)\right).%
\label{prmatrix}
\ee

To proceed, we note that there is a $U_1$ highest weight vector
$\Omega_p\in V_{\ak}\ot V_{\al}$:
\be \Omega_p =  \sli_{i=0}^p \frac{(-1)^i q^{(\ak+1-i)i} }{[i]_q! [p-i]_q!}
v_i^{(\ak)} \ot v_{p-i}^{(\al)},\quad 0\leq p \leq \hbox{min}(\ak,\al),\ee
that has the properties,
\be
&& e_1\Om_p =0,\quad t_1 \Om_p =q^{\ak+\al-2p} \Om_p,\\
&& (1\ot e_1) \Om_p=\Om_{p-1},\quad (e_1\ot t_1) \Omega_p =
-q^{\al+\ak-2(p-1)} \Om_{p-1}.\label{omprop}
\ee
Let $P_p^{(\ak,\al)}$ be
the unique $U_1$ linear map
$P_p^{(\ak,\al)}: V_{\ak}\ot V_{\al} \lra V_{\al} \ot V_{\ak}$
with the properties
\be
P_p^{(\ak,\al)}:\Om_p &\longmapsto& \Om'_p ,\\
P_r^{(\ak,\al)}:\Om_p &\longmapsto& 0,\quad r\ne p,\ee
where $\Om'_p$ is the corresponding highest weight vector in
$V_{\al} \ot V_{\ak}$. 
Then one can follow the argument of~\cite{Ji85} to expand
$R^{(\ak,\al)}(\z)$ in terms of the projectors $P_p^{(\ak,\al)}$.
We find
\be
\HR^{(\ak,\al)}(z_1/z_2) = \sli_{p=0}^{\min\{\ak,\al\}}
\left(\pl_{j=0}^{p-1} \frac{z-q^{\al+\ak-2j}}{ 
1-z q^{\al+\ak-2j}}  \right) P_p^{(\ak,\al)},
\label{rgen}
\ee
where $z=z_1/z_2$.

In the definition of our vertex model, we will use 
the \R-matrix $R^{(\ak,\al)}(\z)$, with $\z$ and $q$ restricted to lie in the
regions $-1<q<0$, $1<\z<-q^{-1}$. 
If we expand
\be
\Rbar^{(\ak,\al)}(\z)=\Rbar^{(\ak,\al)}_0 + (\z-1) \Rbar^{(\ak,\al)}_1
+ O((\z-1)^2),\label{Rexp}\ee
we find
\be 
\lim_{q\ra 0}
\Rbar^{(\ak,\al)}_0 (u_i^{(\ak)} \ot u_j^{(\al)}) =\left\{\begin{array}{lll}
u^{(\al)}_i\ot u^{(\ak)}_j
&\hb{ if }& i+j\leq \ak,\al,\\
u^{(\al)}_{2i+j-\ak}\ot u^{(\ak)}_{\ak-i}&
\hb{ if }& \ak\leq i+j\leq \al,\\
u^{(\al)}_{\al-j}\ot u^{(\ak)}_{i+2j-\al}&
\hb{ if }& \al\leq i+j\leq \ak,\\
u^{(\al)}_{\al-\ak+i}\ot u^{(\ak)}_{j-\al+\ak}&
\hb{ if }& \ak,\al\leq i+j.\end{array}\right.
\label{zerot}
\ee
and
\be 
\lim_{q\ra 0}
\Rbar^{(\ak,\al)}_1(u_i^{(\ak)} \ot u_j^{(\al)})
= \left\{\begin{array}{lll}
(i+j)\, u_i^{(\al)}\ot u_j^{(\ak)} &\hbox{ if }& i+j \leq \ak,\al\\
\ak\, u_{2i+j-\ak}^{(\al)}\ot u_{\ak-i}^{(\ak)}
&\hbox{ if }& \ak\leq i+j\leq \al.\\
\al\, u_{\al-j}^{(\al)}\ot u_{i+2j-\al}^{(\ak)}
&\hbox{ if }& \al\leq i+j \leq \ak.\\
(\al+\ak-i-j)\, u_{\al-\ak+i}^{(\al)}\ot u_{j-\al+\ak}^{(\ak)}
&\hbox{ if }& \ak,\al \leq i+j.
\end{array}
\right. \label{rlimit2}
\ee
These formulas come from equations~\eqref{prmatrix},
\eqref{rgen} and the explicit formula for the projectors $P^{(\ak,\al)}_p$ in the
$q\ra 0$ limit ($P^{(\ak,\al)}_p$ become  diagonal in the
basis $v^{(\ak)}_i \ot v^{(\al)}_j$ in this limit).

The matrix element $R^{(\ak,\al)}(\z_1/\z_2)^{i,j}_{i',j'}$
is the Boltzmann weight associated with the following configuration
of spin variables $i,i'\in\{0,\cdots,\ak\}$ and $j,j'\in\{0,\cdots,\al\}$,
and spectral parameters $\z_1$ and $\z_2$ around a vertex.

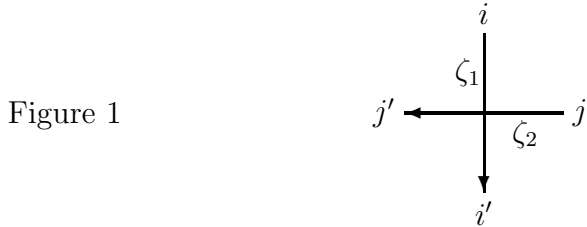
\setlength{\unitlength}{1.2pt}%
 \begin{picture}(100,80)(-100,13)
\thicklines
\put(50,75){\vector(0,-1){50}}
\put(75,50){\vector(-1,0){50}}
\put(50,78){\makebox(0,0)[b]{$i$}}
\put(78,50){\makebox(0,0)[l]{$j$}}
\put(50,22){\makebox(0,0)[t]{$i'$}}
\put(22,50){\makebox(0,0)[r]{$j'$}}
\put(49,63){\makebox(0,0)[r]{$\z_1$}}
\put(63,48){\makebox(0,0)[t]{$\z_2$}}
\put(-100,45){\makebox(0,0)[lb]{Figure 1}}
\end{picture}

\noindent From~\eqref{rlimit2}, we see that if we choose $\z$ and $q$ close
to $1$ and $0_-$ respectively, and consider the case when $\ak\leq \al$,
then the largest Boltzmann weights will be
$R^{(\ak,\al)}(\z_1/\z_2)^{i,j}_{\ak-i,2i+j-\ak}$
with $\ak\leq i+j\leq \al$.
Similarly for $\al\leq \ak$, the largest Boltzmann
weights will be $R^{(\ak,\al)}(\z_1/\z_2)^{i,j}_{i+2j-\al,\al-j}$ with  
$\al\leq i+j\leq \ak$.

\subsection{Definition of the vertex model}
In reference~\cite{HKMW98a}, 
we define the alternating spin vertex model as the vertex model
associated with the two-dimensional lattice consisting of
alternating spin-$\frac{n}{2}$ and spin-$\frac{m}{2}$ lines
(in both the horizontal and vertical directions), where $0<n<m$.
In fact, we choose two vertical and two horizontal spin-$\frac{n}{2}$ lines
next to each other at the centre of our lattice (see Figure 2, in which
the spin-$\frac{n}{2}$ and spin-$\frac{m}{2}$ lines are shown as solid and
dashed lines respectively).
This simplifies our discussion of the corner transfer matrix.

\setlength{\unitlength}{0.00067in}
\begingroup\makeatletter\ifx\SetFigFont\undefined%
\gdef\SetFigFont#1#2#3#4#5{%
  \reset@font\fontsize{#1}{#2pt}%
  \fontfamily{#3}\fontseries{#4}\fontshape{#5}%
  \selectfont}%
\fi\endgroup%
{\renewcommand{\dashlinestretch}{30}
\begin{picture}(7437,6500)(0,-350)
\drawline(4350,5412)(4350,1212)
\drawline(4950,5412)(4950,1212)
\drawline(4950,1212)(4950,12)
\drawline(4920.000,132.000)(4950.000,12.000)(4980.000,132.000)
\drawline(4350,1212)(4350,12)
\drawline(4320.000,132.000)(4350.000,12.000)(4380.000,132.000)
\drawline(6150,5412)(6150,12)
\drawline(6120.000,132.000)(6150.000,12.000)(6180.000,132.000)
\drawline(3150,5412)(3150,12)
\drawline(3120.000,132.000)(3150.000,12.000)(3180.000,132.000)
\drawline(1950,4212)(7350,4212)
\drawline(2070.000,4242.000)(1950.000,4212.000)(2070.000,4182.000)
\drawline(1950,3012)(7350,3012)
\drawline(2070.000,3042.000)(1950.000,3012.000)(2070.000,2982.000)
\drawline(1950,2412)(7350,2412)
\drawline(2070.000,2442.000)(1950.000,2412.000)(2070.000,2382.000)
\drawline(1950,1212)(7350,1212)
\drawline(2070.000,1242.000)(1950.000,1212.000)(2070.000,1182.000)
\dashline{60.000}(5550,5412)(5550,12)
\drawline(5520.000,132.000)(5550.000,12.000)(5580.000,132.000)
\dashline{60.000}(6750,5412)(6750,12)
\drawline(6720.000,132.000)(6750.000,12.000)(6780.000,132.000)
\dashline{60.000}(3750,5412)(3750,12)
\drawline(3720.000,132.000)(3750.000,12.000)(3780.000,132.000)
\dashline{60.000}(2550,5412)(2550,12)
\drawline(2520.000,132.000)(2550.000,12.000)(2580.000,132.000)
\dashline{60.000}(1950,612)(7350,612)
\drawline(2070.000,642.000)(1950.000,612.000)(2070.000,582.000)
\dashline{60.000}(1950,1812)(7350,1812)
\drawline(2070.000,1842.000)(1950.000,1812.000)(2070.000,1782.000)
\dashline{60.000}(1950,3612)(7425,3612)
\drawline(2070.000,3642.000)(1950.000,3612.000)(2070.000,3582.000)
\dashline{60.000}(1950,4812)(7350,4812)
\drawline(2070.000,4842.000)(1950.000,4812.000)(2070.000,4782.000)
\put(6975,387){\makebox(0,0)[lb]{\smash{{{\SetFigFont{8}{14.4}{\rmdefault}%
{\mddefault}{\updefault}$\bar{j_1}$}}}}}
\put(6975,987){\makebox(0,0)[lb]{\smash{{{\SetFigFont{8}{14.4}{\rmdefault}%
{\mddefault}{\updefault}$\bar{i}$}}}}}
\put(6975,1587){\makebox(0,0)[lb]{\smash{{{\SetFigFont{8}{14.4}{\rmdefault}%
{\mddefault}{\updefault}$\bar{j}$}}}}}
\put(6975,2187){\makebox(0,0)[lb]{\smash{{{\SetFigFont{8}{14.4}{\rmdefault}%
{\mddefault}{\updefault}$\bar{i}$}}}}}
\put(6975,2787){\makebox(0,0)[lb]{\smash{{{\SetFigFont{8}{14.4}{\rmdefault}%
{\mddefault}{\updefault}$i$}}}}}
\put(6975,3387){\makebox(0,0)[lb]{\smash{{{\SetFigFont{8}{14.4}{\rmdefault}%
{\mddefault}{\updefault}$j$}}}}}
\put(6975,3987){\makebox(0,0)[lb]{\smash{{{\SetFigFont{8}{14.4}{\rmdefault}%
{\mddefault}{\updefault}$i$}}}}}
\put(6975,4587){\makebox(0,0)[lb]{\smash{{{\SetFigFont{8}{14.4}{\rmdefault}%
{\mddefault}{\updefault}$j_1$}}}}}
\put(6375,387){\makebox(0,0)[lb]{\smash{{{\SetFigFont{8}{14.4}{\rmdefault}%
{\mddefault}{\updefault}$j_1$}}}}}
\put(6375,987){\makebox(0,0)[lb]{\smash{{{\SetFigFont{8}{14.4}{\rmdefault}%
{\mddefault}{\updefault}$i$}}}}}
\put(6375,1587){\makebox(0,0)[lb]{\smash{{{\SetFigFont{8}{14.4}{\rmdefault}%
{\mddefault}{\updefault}$j$}}}}}
\put(6375,2187){\makebox(0,0)[lb]{\smash{{{\SetFigFont{8}{14.4}{\rmdefault}%
{\mddefault}{\updefault}$i$}}}}}
\put(6375,2787){\makebox(0,0)[lb]{\smash{{{\SetFigFont{8}{14.4}{\rmdefault}%
{\mddefault}{\updefault}$\bi$}}}}}
\put(6375,3387){\makebox(0,0)[lb]{\smash{{{\SetFigFont{8}{14.4}{\rmdefault}%
{\mddefault}{\updefault}$\bj$}}}}}
\put(6375,3987){\makebox(0,0)[lb]{\smash{{{\SetFigFont{8}{14.4}{\rmdefault}%
{\mddefault}{\updefault}$\bi$}}}}}
\put(6375,4587){\makebox(0,0)[lb]{\smash{{{\SetFigFont{8}{14.4}{\rmdefault}%
{\mddefault}{\updefault}$\bj_1$}}}}}
\put(5775,387){\makebox(0,0)[lb]{\smash{{{\SetFigFont{8}{14.4}{\rmdefault}%
{\mddefault}{\updefault}$\bar{j}$}}}}}
\put(5775,987){\makebox(0,0)[lb]{\smash{{{\SetFigFont{8}{14.4}{\rmdefault}%
{\mddefault}{\updefault}$\bar{i}$}}}}}
\put(5775,1587){\makebox(0,0)[lb]{\smash{{{\SetFigFont{8}{14.4}{\rmdefault}%
{\mddefault}{\updefault}$\bar{j_1}$}}}}}
\put(5775,2187){\makebox(0,0)[lb]{\smash{{{\SetFigFont{8}{14.4}{\rmdefault}%
{\mddefault}{\updefault}$\bar{i}$}}}}}
\put(5775,2787){\makebox(0,0)[lb]{\smash{{{\SetFigFont{8}{14.4}{\rmdefault}%
{\mddefault}{\updefault}$i$}}}}}
\put(5775,3387){\makebox(0,0)[lb]{\smash{{{\SetFigFont{8}{14.4}{\rmdefault}%
{\mddefault}{\updefault}$j_1$}}}}}
\put(5775,3987){\makebox(0,0)[lb]{\smash{{{\SetFigFont{8}{14.4}{\rmdefault}%
{\mddefault}{\updefault}$i$}}}}}
\put(5775,4587){\makebox(0,0)[lb]{\smash{{{\SetFigFont{8}{14.4}{\rmdefault}%
{\mddefault}{\updefault}$j$}}}}}
\put(5175,387){\makebox(0,0)[lb]{\smash{{{\SetFigFont{8}{14.4}{\rmdefault}%
{\mddefault}{\updefault}$j$}}}}}
\put(5175,987){\makebox(0,0)[lb]{\smash{{{\SetFigFont{8}{14.4}{\rmdefault}%
{\mddefault}{\updefault}$i$}}}}}
\put(5175,1587){\makebox(0,0)[lb]{\smash{{{\SetFigFont{8}{14.4}{\rmdefault}%
{\mddefault}{\updefault}$j_1$}}}}}
\put(5175,2187){\makebox(0,0)[lb]{\smash{{{\SetFigFont{8}{14.4}{\rmdefault}%
{\mddefault}{\updefault}$i$}}}}}
\put(5175,2787){\makebox(0,0)[lb]{\smash{{{\SetFigFont{8}{14.4}{\rmdefault}%
{\mddefault}{\updefault}$\bi$}}}}}
\put(5175,3387){\makebox(0,0)[lb]{\smash{{{\SetFigFont{8}{14.4}{\rmdefault}%
{\mddefault}{\updefault}$\bj_1$}}}}}
\put(5175,3987){\makebox(0,0)[lb]{\smash{{{\SetFigFont{8}{14.4}{\rmdefault}%
{\mddefault}{\updefault}$\bi$}}}}}
\put(5175,4587){\makebox(0,0)[lb]{\smash{{{\SetFigFont{8}{14.4}{\rmdefault}%
{\mddefault}{\updefault}$\bj$}}}}}
\put(4575,387){\makebox(0,0)[lb]{\smash{{{\SetFigFont{8}{14.4}{\rmdefault}%
{\mddefault}{\updefault}$\bar{j_1}$}}}}}
\put(4575,987){\makebox(0,0)[lb]{\smash{{{\SetFigFont{8}{14.4}{\rmdefault}%
{\mddefault}{\updefault}$\bar{i}$}}}}}
\put(4575,1587){\makebox(0,0)[lb]{\smash{{{\SetFigFont{8}{14.4}{\rmdefault}%
{\mddefault}{\updefault}$\bar{j}$}}}}}
\put(4575,2187){\makebox(0,0)[lb]{\smash{{{\SetFigFont{8}{14.4}{\rmdefault}%
{\mddefault}{\updefault}$\bar{i}$}}}}}
\put(4575,2787){\makebox(0,0)[lb]{\smash{{{\SetFigFont{8}{14.4}{\rmdefault}%
{\mddefault}{\updefault}$i$}}}}}
\put(4575,3387){\makebox(0,0)[lb]{\smash{{{\SetFigFont{8}{14.4}{\rmdefault}%
{\mddefault}{\updefault}$j$}}}}}
\put(4575,3987){\makebox(0,0)[lb]{\smash{{{\SetFigFont{8}{14.4}{\rmdefault}%
{\mddefault}{\updefault}$i$}}}}}
\put(4575,4587){\makebox(0,0)[lb]{\smash{{{\SetFigFont{8}{14.4}{\rmdefault}%
{\mddefault}{\updefault}$j_1$}}}}}
\put(3975,387){\makebox(0,0)[lb]{\smash{{{\SetFigFont{8}{14.4}{\rmdefault}%
{\mddefault}{\updefault}$j$}}}}}
\put(3975,987){\makebox(0,0)[lb]{\smash{{{\SetFigFont{8}{14.4}{\rmdefault}%
{\mddefault}{\updefault}$i$}}}}}
\put(3975,1587){\makebox(0,0)[lb]{\smash{{{\SetFigFont{8}{14.4}{\rmdefault}%
{\mddefault}{\updefault}$j_1$}}}}}
\put(3975,2187){\makebox(0,0)[lb]{\smash{{{\SetFigFont{8}{14.4}{\rmdefault}%
{\mddefault}{\updefault}$i$}}}}}
\put(3975,2787){\makebox(0,0)[lb]{\smash{{{\SetFigFont{8}{14.4}{\rmdefault}%
{\mddefault}{\updefault}$\bi$}}}}}
\put(3975,3387){\makebox(0,0)[lb]{\smash{{{\SetFigFont{8}{14.4}{\rmdefault}%
{\mddefault}{\updefault}$\bj_1$}}}}}
\put(3975,3987){\makebox(0,0)[lb]{\smash{{{\SetFigFont{8}{14.4}{\rmdefault}%
{\mddefault}{\updefault}$\bi$}}}}}
\put(3975,4587){\makebox(0,0)[lb]{\smash{{{\SetFigFont{8}{14.4}{\rmdefault}%
{\mddefault}{\updefault}$\bj$}}}}}
\put(3375,387){\makebox(0,0)[lb]{\smash{{{\SetFigFont{8}{14.4}{\rmdefault}%
{\mddefault}{\updefault}$\bar{j}$}}}}}
\put(3375,987){\makebox(0,0)[lb]{\smash{{{\SetFigFont{8}{14.4}{\rmdefault}%
{\mddefault}{\updefault}$\bar{i}$}}}}}
\put(3375,1587){\makebox(0,0)[lb]{\smash{{{\SetFigFont{8}{14.4}{\rmdefault}%
{\mddefault}{\updefault}$\bar{j_1}$}}}}}
\put(3375,2187){\makebox(0,0)[lb]{\smash{{{\SetFigFont{8}{14.4}{\rmdefault}%
{\mddefault}{\updefault}$\bar{i}$}}}}}
\put(3375,2787){\makebox(0,0)[lb]{\smash{{{\SetFigFont{8}{14.4}{\rmdefault}%
{\mddefault}{\updefault}$i$}}}}}
\put(3375,3387){\makebox(0,0)[lb]{\smash{{{\SetFigFont{8}{14.4}{\rmdefault}%
{\mddefault}{\updefault}$j_1$}}}}}
\put(3375,3987){\makebox(0,0)[lb]{\smash{{{\SetFigFont{8}{14.4}{\rmdefault}%
{\mddefault}{\updefault}$i$}}}}}
\put(3375,4587){\makebox(0,0)[lb]{\smash{{{\SetFigFont{8}{14.4}{\rmdefault}%
{\mddefault}{\updefault}$j$}}}}}
\put(2775,387){\makebox(0,0)[lb]{\smash{{{\SetFigFont{8}{14.4}{\rmdefault}%
{\mddefault}{\updefault}$j_1$}}}}}
\put(2775,987){\makebox(0,0)[lb]{\smash{{{\SetFigFont{8}{14.4}{\rmdefault}%
{\mddefault}{\updefault}$i$}}}}}
\put(2775,1587){\makebox(0,0)[lb]{\smash{{{\SetFigFont{8}{14.4}{\rmdefault}%
{\mddefault}{\updefault}$j$}}}}}
\put(2775,2187){\makebox(0,0)[lb]{\smash{{{\SetFigFont{8}{14.4}{\rmdefault}%
{\mddefault}{\updefault}$i$}}}}}
\put(2775,2787){\makebox(0,0)[lb]{\smash{{{\SetFigFont{8}{14.4}{\rmdefault}%
{\mddefault}{\updefault}$\bi$}}}}}
\put(2775,3387){\makebox(0,0)[lb]{\smash{{{\SetFigFont{8}{14.4}{\rmdefault}%
{\mddefault}{\updefault}$\bj$}}}}}
\put(2775,3987){\makebox(0,0)[lb]{\smash{{{\SetFigFont{8}{14.4}{\rmdefault}%
{\mddefault}{\updefault}$\bi$}}}}}
\put(2775,4587){\makebox(0,0)[lb]{\smash{{{\SetFigFont{8}{14.4}{\rmdefault}%
{\mddefault}{\updefault}$\bj_1$}}}}}
\put(2175,387){\makebox(0,0)[lb]{\smash{{{\SetFigFont{8}{14.4}{\rmdefault}%
{\mddefault}{\updefault}$\bar{j_1}$}}}}}
\put(2175,987){\makebox(0,0)[lb]{\smash{{{\SetFigFont{8}{14.4}{\rmdefault}%
{\mddefault}{\updefault}$\bar{i}$}}}}}
\put(2175,1587){\makebox(0,0)[lb]{\smash{{{\SetFigFont{8}{14.4}{\rmdefault}%
{\mddefault}{\updefault}$\bar{j}$}}}}}
\put(2175,2187){\makebox(0,0)[lb]{\smash{{{\SetFigFont{8}{14.4}{\rmdefault}%
{\mddefault}{\updefault}$\bar{i}$}}}}}
\put(2175,2787){\makebox(0,0)[lb]{\smash{{{\SetFigFont{8}{14.4}{\rmdefault}%
{\mddefault}{\updefault}$i$}}}}}
\put(2175,3387){\makebox(0,0)[lb]{\smash{{{\SetFigFont{8}{14.4}{\rmdefault}%
{\mddefault}{\updefault}$j$}}}}}
\put(2175,3987){\makebox(0,0)[lb]{\smash{{{\SetFigFont{8}{14.4}{\rmdefault}%
{\mddefault}{\updefault}$i$}}}}}
\put(2175,4587){\makebox(0,0)[lb]{\smash{{{\SetFigFont{8}{14.4}{\rmdefault}%
{\mddefault}{\updefault}$j_1$}}}}}
\put(2400,5037){\makebox(0,0)[lb]{\smash{{{\SetFigFont{8}{14.4}{\rmdefault}%
{\mddefault}{\updefault}$j_1$}}}}}
\put(2400,4437){\makebox(0,0)[lb]{\smash{{{\SetFigFont{8}{14.4}{\rmdefault}%
{\mddefault}{\updefault}$\bj_1$}}}}}
\put(2400,3837){\makebox(0,0)[lb]{\smash{{{\SetFigFont{8}{14.4}{\rmdefault}%
{\mddefault}{\updefault}$j$}}}}}
\put(2400,3237){\makebox(0,0)[lb]{\smash{{{\SetFigFont{8}{14.4}{\rmdefault}%
{\mddefault}{\updefault}$\bj$}}}}}
\put(2400,2637){\makebox(0,0)[lb]{\smash{{{\SetFigFont{8}{14.4}{\rmdefault}%
{\mddefault}{\updefault}$j_1$}}}}}
\put(2400,2037){\makebox(0,0)[lb]{\smash{{{\SetFigFont{8}{14.4}{\rmdefault}%
{\mddefault}{\updefault}$\bj$}}}}}
\put(2400,1437){\makebox(0,0)[lb]{\smash{{{\SetFigFont{8}{14.4}{\rmdefault}%
{\mddefault}{\updefault}$j$}}}}}
\put(2400,837){\makebox(0,0)[lb]{\smash{{{\SetFigFont{8}{14.4}{\rmdefault}%
{\mddefault}{\updefault}$\bj_1$}}}}}
\put(2400,237){\makebox(0,0)[lb]{\smash{{{\SetFigFont{8}{14.4}{\rmdefault}%
{\mddefault}{\updefault}$j_1$}}}}}
\put(3000,5037){\makebox(0,0)[lb]{\smash{{{\SetFigFont{8}{14.4}{\rmdefault}%
{\mddefault}{\updefault}$i$}}}}}
\put(3000,4437){\makebox(0,0)[lb]{\smash{{{\SetFigFont{8}{14.4}{\rmdefault}%
{\mddefault}{\updefault}$\bi$}}}}}
\put(3000,3837){\makebox(0,0)[lb]{\smash{{{\SetFigFont{8}{14.4}{\rmdefault}%
{\mddefault}{\updefault}$i$}}}}}
\put(3000,3237){\makebox(0,0)[lb]{\smash{{{\SetFigFont{8}{14.4}{\rmdefault}%
{\mddefault}{\updefault}$\bi$}}}}}
\put(3000,2637){\makebox(0,0)[lb]{\smash{{{\SetFigFont{8}{14.4}{\rmdefault}%
{\mddefault}{\updefault}$i$}}}}}
\put(3000,2037){\makebox(0,0)[lb]{\smash{{{\SetFigFont{8}{14.4}{\rmdefault}%
{\mddefault}{\updefault}$\bi$}}}}}
\put(3000,1437){\makebox(0,0)[lb]{\smash{{{\SetFigFont{8}{14.4}{\rmdefault}%
{\mddefault}{\updefault}$i$}}}}}
\put(3000,837){\makebox(0,0)[lb]{\smash{{{\SetFigFont{8}{14.4}{\rmdefault}%
{\mddefault}{\updefault}$\bi$}}}}}
\put(3000,237){\makebox(0,0)[lb]{\smash{{{\SetFigFont{8}{14.4}{\rmdefault}%
{\mddefault}{\updefault}$i$}}}}}
\put(3600,5037){\makebox(0,0)[lb]{\smash{{{\SetFigFont{8}{14.4}{\rmdefault}%
{\mddefault}{\updefault}$j$}}}}}
\put(3600,4437){\makebox(0,0)[lb]{\smash{{{\SetFigFont{8}{14.4}{\rmdefault}%
{\mddefault}{\updefault}$\bj$}}}}}
\put(3600,3837){\makebox(0,0)[lb]{\smash{{{\SetFigFont{8}{14.4}{\rmdefault}%
{\mddefault}{\updefault}$j_1$}}}}}
\put(3600,3237){\makebox(0,0)[lb]{\smash{{{\SetFigFont{8}{14.4}{\rmdefault}%
{\mddefault}{\updefault}$\bj_1$}}}}}
\put(3600,2637){\makebox(0,0)[lb]{\smash{{{\SetFigFont{8}{14.4}{\rmdefault}%
{\mddefault}{\updefault}$j$}}}}}
\put(3600,2037){\makebox(0,0)[lb]{\smash{{{\SetFigFont{8}{14.4}{\rmdefault}%
{\mddefault}{\updefault}$\bj_1$}}}}}
\put(3600,1437){\makebox(0,0)[lb]{\smash{{{\SetFigFont{8}{14.4}{\rmdefault}%
{\mddefault}{\updefault}$j_1$}}}}}
\put(3600,837){\makebox(0,0)[lb]{\smash{{{\SetFigFont{8}{14.4}{\rmdefault}%
{\mddefault}{\updefault}$\bj$}}}}}
\put(3600,237){\makebox(0,0)[lb]{\smash{{{\SetFigFont{8}{14.4}{\rmdefault}%
{\mddefault}{\updefault}$j$}}}}}
\put(4200,5037){\makebox(0,0)[lb]{\smash{{{\SetFigFont{8}{14.4}{\rmdefault}%
{\mddefault}{\updefault}$i$}}}}}
\put(4200,4437){\makebox(0,0)[lb]{\smash{{{\SetFigFont{8}{14.4}{\rmdefault}%
{\mddefault}{\updefault}$\bi$}}}}}
\put(4200,3837){\makebox(0,0)[lb]{\smash{{{\SetFigFont{8}{14.4}{\rmdefault}%
{\mddefault}{\updefault}$i$}}}}}
\put(4200,3237){\makebox(0,0)[lb]{\smash{{{\SetFigFont{8}{14.4}{\rmdefault}%
{\mddefault}{\updefault}$\bi$}}}}}
\put(4200,2637){\makebox(0,0)[lb]{\smash{{{\SetFigFont{8}{14.4}{\rmdefault}%
{\mddefault}{\updefault}$i$}}}}}
\put(4200,2037){\makebox(0,0)[lb]{\smash{{{\SetFigFont{8}{14.4}{\rmdefault}%
{\mddefault}{\updefault}$\bi$}}}}}
\put(4200,1437){\makebox(0,0)[lb]{\smash{{{\SetFigFont{8}{14.4}{\rmdefault}%
{\mddefault}{\updefault}$i$}}}}}
\put(4200,837){\makebox(0,0)[lb]{\smash{{{\SetFigFont{8}{14.4}{\rmdefault}%
{\mddefault}{\updefault}$\bi$}}}}}
\put(4200,237){\makebox(0,0)[lb]{\smash{{{\SetFigFont{8}{14.4}{\rmdefault}%
{\mddefault}{\updefault}$i$}}}}}
\put(4800,5037){\makebox(0,0)[lb]{\smash{{{\SetFigFont{8}{14.4}{\rmdefault}%
{\mddefault}{\updefault}$\bi$}}}}}
\put(4800,4437){\makebox(0,0)[lb]{\smash{{{\SetFigFont{8}{14.4}{\rmdefault}%
{\mddefault}{\updefault}$i$}}}}}
\put(4800,3837){\makebox(0,0)[lb]{\smash{{{\SetFigFont{8}{14.4}{\rmdefault}%
{\mddefault}{\updefault}$\bi$}}}}}
\put(4800,3237){\makebox(0,0)[lb]{\smash{{{\SetFigFont{8}{14.4}{\rmdefault}%
{\mddefault}{\updefault}$i$}}}}}
\put(4800,2637){\makebox(0,0)[lb]{\smash{{{\SetFigFont{8}{14.4}{\rmdefault}%
{\mddefault}{\updefault}$\bi$}}}}}
\put(4800,2037){\makebox(0,0)[lb]{\smash{{{\SetFigFont{8}{14.4}{\rmdefault}%
{\mddefault}{\updefault}$i$}}}}}
\put(4800,1437){\makebox(0,0)[lb]{\smash{{{\SetFigFont{8}{14.4}{\rmdefault}%
{\mddefault}{\updefault}$\bi$}}}}}
\put(4800,837){\makebox(0,0)[lb]{\smash{{{\SetFigFont{8}{14.4}{\rmdefault}%
{\mddefault}{\updefault}$i$}}}}}
\put(4800,237){\makebox(0,0)[lb]{\smash{{{\SetFigFont{8}{14.4}{\rmdefault}%
{\mddefault}{\updefault}$\bi$}}}}}
\put(5400,5037){\makebox(0,0)[lb]{\smash{{{\SetFigFont{8}{14.4}{\rmdefault}%
{\mddefault}{\updefault}$\bj$}}}}}
\put(5400,4437){\makebox(0,0)[lb]{\smash{{{\SetFigFont{8}{14.4}{\rmdefault}%
{\mddefault}{\updefault}$j$}}}}}
\put(5400,3837){\makebox(0,0)[lb]{\smash{{{\SetFigFont{8}{14.4}{\rmdefault}%
{\mddefault}{\updefault}$\bj_1$}}}}}
\put(5400,3237){\makebox(0,0)[lb]{\smash{{{\SetFigFont{8}{14.4}{\rmdefault}%
{\mddefault}{\updefault}$j_1$}}}}}
\put(5400,2637){\makebox(0,0)[lb]{\smash{{{\SetFigFont{8}{14.4}{\rmdefault}%
{\mddefault}{\updefault}$\bj$}}}}}
\put(5400,2037){\makebox(0,0)[lb]{\smash{{{\SetFigFont{8}{14.4}{\rmdefault}%
{\mddefault}{\updefault}$j_1$}}}}}
\put(5400,1437){\makebox(0,0)[lb]{\smash{{{\SetFigFont{8}{14.4}{\rmdefault}%
{\mddefault}{\updefault}$\bj_1$}}}}}
\put(5400,837){\makebox(0,0)[lb]{\smash{{{\SetFigFont{8}{14.4}{\rmdefault}%
{\mddefault}{\updefault}$j$}}}}}
\put(5400,237){\makebox(0,0)[lb]{\smash{{{\SetFigFont{8}{14.4}{\rmdefault}%
{\mddefault}{\updefault}$\bj$}}}}}
\put(6000,5037){\makebox(0,0)[lb]{\smash{{{\SetFigFont{8}{14.4}{\rmdefault}%
{\mddefault}{\updefault}$\bi$}}}}}
\put(6000,4437){\makebox(0,0)[lb]{\smash{{{\SetFigFont{8}{14.4}{\rmdefault}%
{\mddefault}{\updefault}$i$}}}}}
\put(6000,3837){\makebox(0,0)[lb]{\smash{{{\SetFigFont{8}{14.4}{\rmdefault}%
{\mddefault}{\updefault}$\bi$}}}}}
\put(6000,3237){\makebox(0,0)[lb]{\smash{{{\SetFigFont{8}{14.4}{\rmdefault}%
{\mddefault}{\updefault}$i$}}}}}
\put(6000,2637){\makebox(0,0)[lb]{\smash{{{\SetFigFont{8}{14.4}{\rmdefault}%
{\mddefault}{\updefault}$\bi$}}}}}
\put(6000,2037){\makebox(0,0)[lb]{\smash{{{\SetFigFont{8}{14.4}{\rmdefault}%
{\mddefault}{\updefault}$i$}}}}}
\put(6000,1437){\makebox(0,0)[lb]{\smash{{{\SetFigFont{8}{14.4}{\rmdefault}%
{\mddefault}{\updefault}$\bi$}}}}}
\put(6000,837){\makebox(0,0)[lb]{\smash{{{\SetFigFont{8}{14.4}{\rmdefault}%
{\mddefault}{\updefault}$i$}}}}}
\put(6000,237){\makebox(0,0)[lb]{\smash{{{\SetFigFont{8}{14.4}{\rmdefault}%
{\mddefault}{\updefault}$\bi$}}}}}
\put(6600,5037){\makebox(0,0)[lb]{\smash{{{\SetFigFont{8}{14.4}{\rmdefault}%
{\mddefault}{\updefault}$\bj_1$}}}}}
\put(6600,4437){\makebox(0,0)[lb]{\smash{{{\SetFigFont{8}{14.4}{\rmdefault}%
{\mddefault}{\updefault}$j_1$}}}}}
\put(6600,3837){\makebox(0,0)[lb]{\smash{{{\SetFigFont{8}{14.4}{\rmdefault}%
{\mddefault}{\updefault}$\bj$}}}}}
\put(6600,3237){\makebox(0,0)[lb]{\smash{{{\SetFigFont{8}{14.4}{\rmdefault}%
{\mddefault}{\updefault}$j$}}}}}
\put(6600,2637){\makebox(0,0)[lb]{\smash{{{\SetFigFont{8}{14.4}{\rmdefault}%
{\mddefault}{\updefault}$\bj_1$}}}}}
\put(6600,2037){\makebox(0,0)[lb]{\smash{{{\SetFigFont{8}{14.4}{\rmdefault}%
{\mddefault}{\updefault}$j$}}}}}
\put(6600,1437){\makebox(0,0)[lb]{\smash{{{\SetFigFont{8}{14.4}{\rmdefault}%
{\mddefault}{\updefault}$\bj$}}}}}
\put(6600,837){\makebox(0,0)[lb]{\smash{{{\SetFigFont{8}{14.4}{\rmdefault}%
{\mddefault}{\updefault}$j_1$}}}}}
\put(6600,237){\makebox(0,0)[lb]{\smash{{{\SetFigFont{8}{14.4}{\rmdefault}%
{\mddefault}{\updefault}$\bj_1$}}}}}

\put(0,2712){\makebox(0,0)[lb]{\smash{{{\SetFigFont{12}{14.4}{\rmdefault}%
{\mddefault}{\updefault}Figure 2}}}}}

\end{picture}
}

\noindent Vertical lines will carry a spectral parameter equal to $\z$,
and horizontal
lines a spectral parameter equal to 1. We restrict our discussion to
the anti-ferromagnetic region $-1<q<0$, $1<\z<-q^{-1}$.
The 
different 
local Boltzmann weights associated with the intersection vertices of
this lattice are given by the \R-matrices
$R^{(n,m)}(\z)$,
$R^{(m,n)}(\z)$, $R^{(n,n)}(\z)$, and $R^{(m,m)}(\z)$. 

A ground state of such a vertex model is a configuration of the
spin variables for which all of the local vertex configurations
are associated with one of the largest Boltzmann weights discussed
above. There are $(m-n+1)(n+1)$ different anti-ferromagnetic 
ground states for our model,
each labelled by a pair of integers $(a,b)$, where 
$0\leq a\leq m-n$ and $0\leq b\leq n$. 
The spin configuration in the $(a,b)$ ground state is given
in Figure 2 in which we use the notation
$i=n-b$, $j=m-n-a+b$, $j_1=a+b$, $\bi=b$, $\bj=n+a-b$, $\bj_1=m-a-b$.
Define $Z_{a,b;N}$ to be the partition function (i.e., the weighted
configuration sum) of
such a lattice which consists of $N$ vertices, and whose boundary
spins are fixed to the
values
of the $(a,b)$ ground state. With our normalisation of the
\R-matrices, the partition function per unit site, $\lim_{N\ra
  \infty}Z_{a,b;N}^{1/N}$, is equal to 1. We are interested in the
infinite-volume lattice with partition function
$Z_{a,b}=\lim_{N\ra\infty}Z_{a,b;N}$. $Z_{a,b}$ is divergent. 
However, this divergence cancels for 
correlation functions since
they are given as ratios (see~\cite{HKMW98a}).

We can identify $Z_{a,b}$ with the trace of corner transfer matrices
\be
Z_{a,b}=\Tr_{\cH_{a,b}}(A_{NE}(\z) A_{SE}(\z) A_{SW}(\z) A_{NW}(\z)).\ee
Let us explain the various elements in this formula. $\cH_{a,b}$ is the
space of eigenstates of the corner transfer matrix $A_{NW}(\z)$
associated
with the North-West quadrant of the lattice. In the limit $q\ra 0$, we
can identify $\cH_{a,b}$ with the path space $P_{a,b}$. The latter is
defined to be the set of paths $\ket{p}=\cdots
p(3)\,p(2)\,p(1)$ with the following restrictions:
\be
p(\tts)&\in& \{0,1,\cdots,n\}\quad \hb{if }k\hb{ is odd},\label{gs1}\\
p(\tts)&\in&\{0,1,\cdots,m\}\quad \hb{if }k\hb{ is even},\\
p(\tts)&=&\bar{p}(\tts;a,b),\ws
\tts\gg 0,\ws\hb{where}\\[2mm]
\bar{p}(\tts;a,b)&=&\left\{
\begin{array}{lll}n-b &\hb{ if }k\hb{ is odd};\\
                  a+b&\hb{ if } k \equiv 0\pmod 4;\\
                  m-n-a+b&\hb{ if }  k\equiv 2\pmod 4.\end{array}\right.
 \label{gs2}\ena
A path $\ket{p}\in P_{a,b}$ corresponds to a particular choice of
the spin variables on the half-infinite column of horizontal edges
running North from the centre of our lattice. The boundary condition
$p(\tts)=\bar{p}(\tts;a,b),\ws \tts\gg 0$ corresponds to the choice of 
the $(a,b)$ ground state.
If $A_{NW}(\z)$ acts on some $\ket{p}\in P_{a,b}$, then it will produce an
infinite linear combination of paths $\ket{p'}\in P_{a,b}$. One term
will be of order $q^0$ (see~\eqref{zerot}), 
and all the others of higher order in
$q$. The infinite linear combination is not in $P_{a,b}$.
For $q\ne 0$, $A_{NW}(\z)$ should be renormalised as a map 
$\cH_{a,b}\ra \cH_{a,b}$, where the space
$\cH_{a,b}$ will be identified in terms of the representation
theory of $\uq$ in Section~\ref{sec6}.

The corner transfer
matrices corresponding to the other quadrants can be identified as the
maps $A_{SW}(\z):\cH_{a,b}\ra \cH_{m-n-a,n-b}$,
$A_{SE}(\z):\cH_{m-n-a,n-b}\ra \cH_{m-n-a,n-b}$, and
$A_{NE}(\z):\cH_{m-n-a,n-b}\ra \cH_{a,b}$.
One can construct heuristic arguments along the lines of those
in~\cite{JMN,Fodal93} (which rely upon the crossing and unitarity properties 
of our \R-matrix; given by~\eqref{crossing} and~\eqref{unit}
respectively),
to yield the following relations among the different
corner transfer matrices:
\be A_{SW}(\z)=C
A_{NW}(-q^{-1}\z^{-1}),\; A_{SE}(\z)=C A_{NW}(\z) C,\;
A_{NE}(\z)= A_{NW}(-q^{-1} \z^{-1})C.\nonumber\\
\label{z2}\ee
Here, $C$ is the `conjugation operator': In the limit $q\ra 0$, it is
the operator $P_{a,b}\ra P_{m-n-a,n-b}$  defined by
\be
p(\tts)\to \left\{
\begin{array}{ll}
n-p(\tts) & \hb{ if $\tts$ is odd};\\
m-p(\tts) &\hb{ if $\tts$ is even}.
\end{array}\right.
\nn\ee
When $q\ne 0$, it will be the operator $\cH_{a,b}\ra \cH_{m-n-a,n-b}$ 
which exchanges the fundamental weights $\La_0\leftrightarrow \La_1$ of
$\uq$.

The corner transfer matrix $A_{NW}(\z)$ has a 
remarkably simple form in the infinite-volume limit. Baxter's 
argument (see~\cite{Bax82}) applied here implies $A_{NW}(\z)=c(\z)
\z^{-H_{CTM}}$. Here, $c(\z)$ is a divergent scalar.
$H_{CTM}$ is the corner transfer matrix Hamiltonian, which 
is independent of $\z$ and has a non-negative integer spectrum.
Using~\eqref{z2}, we then 
have that up to a divergent scalar the infinite-volume 
partition function $Z_{a,b}$ is proportional to
$\Tr_{\cH_{a,b}}((-q)^{2 H_{CTM}})$.

\subsection{The corner transfer matrix Hamiltonian $H_{CTM}$}
The corner transfer matrix Hamiltonian $H_{CTM}$ is  
defined by $H_{CTM}=-\frac{d}{d\z}
A_{NW}(\z)|_{\z=1}: \cH_{a,b}\ra \cH_{a,b}$.
Its action on a path  $\ket{p}\in
P_{a,b}$ can be calculated from~\eqref{Rexp}:
\be H_{CTM} = -\sli_{s=1}^{\infty} s \left(  H_{1;2s+1,2s,2s-1}  +
H_{2;2s+2,2s+1,2s}  + 2  H_{3;2s+1,2s}  \right).\label{HCTM}\ee 
Here, $ H_{1;2s+1,2s,2s-1}$ acts as the identity on $\ket{p}\in
P_{a,b}$ except at the positions $2s+1$, $2s$, $2s-1$, where its action,
written in terms of $\Rbar^{(\ak,\al)}_0$
and $\Rbar^{(\ak,\al)}_1$ as defined in~\eqref{Rexp}, is given by
\be  H_{1} &=& (\Rbar^{(m,n)}_0 \ot 1) (1\ot
\Rbar_1^{(n,n)}) (\Rbar_0^{(n,m)}\ot 1).
\ee
Similarly, $ H_{2;2s+2,2s+1,2s} $ acts as the identity except at
the positions $2s+2$, $2s+1$, $2s$, where it acts as
\be  H_{2} &=& (1\ot \Rbar^{(m,n)}_0 ) (
\Rbar_1^{(m,m)}\ot 1) (1\ot \Rbar_0^{(n,m)}).\ee
Finally, $H_{3;2s+1,2s} $ acts as the identity except at the positions
$2s+1$, $2s$, where it acts as 
\be
H_{2s+1,2s}&=& \Rbar^{(m,n)}_0  \Rbar_1^{(n,m)} =
\Rbar^{(m,n)}_1  \Rbar_0^{(n,m)}.\ee
The equality of the last two expressions follows from the unitarity property
\eqref{unit}.

In the limit $q\ra 0$, $H_1,H_2,H_3:P_{a,b}\ra P_{a,b}$ act diagonally.
Let us use the notation
\bea
\lim_{q\ra 0}H_1(u^{(n)}_i \ot u^{(m)}_j \ot u^{(n)}_k)
&=&h_1(i,j,k) (u^{(n)}_i \ot u^{(m)}_j \ot u^{(n)}_k),\\
\lim_{q\ra 0}H_2(u^{(m)}_i \ot u^{(n)}_j \ot u^{(m)}_k)
&=&h_2(i,j,k) (u^{(m)}_i \ot u^{(n)}_j \ot u^{(m)}_k),\\
\lim_{q\ra 0}H_3(u^{(n)}_i \ot u^{(m)}_j)&=&
h_3(i,j)(u^{(n)}_i \ot u^{(m)}_j). \ena
Using~\eqref{zerot} and~\eqref{rlimit2}, we find 
\bea
h_1(i,j,k)&=&\left\{\begin{array}{lll} 
                       {\{}k+j\}_n & \hbox{ if }  & i+j\leq m,n;\\
                       {\{}k+i+2j-m\}_n &  \hbox{ if } &m\leq i+j \leq n;\\
                       {\{}k+n-i\}_n &  \hbox{ if } & n\leq i+j \leq m;\\
                       {\{}k+n-m+j\}_n & \hbox{ if } & n,m \leq i+j, %
\end{array}\right. 
\label{h1}\\[2mm]
h_2(i,j,k)&=&\left\{\begin{array}{lll} 
                       {\{}i+j\}_m & \hbox{ if }  & j+k\leq m,n;\\
                       {\{}i+m-k\}_m &  \hbox{ if } &m\leq j+k \leq n;\\
                       {\{}i+k+2j-n\}_m &  \hbox{ if } & n\leq j+k \leq m;\\
                       {\{}i+j+m-n\}_m & \hbox{ if } & n,m \leq j+k, %
\end{array}\right.
\label{h2}\\[2mm]
h_3(i,j)&=&\left\{\begin{array}{lll} 
                       i+j & \hbox{ if }  & i+j\leq m,n;\\
                       m &  \hbox{ if } &m\leq i+j \leq n;\\
                       n &  \hbox{ if } & n\leq i+j \leq m;\\
                       m+n-i-j & \hbox{ if } & n,m \leq i+j. %
\end{array}\right.
\label{h3}\ena
Here we have used the notation 
\be \{a\}_b=\left\{\begin{array}{lll} 
                        a & \hbox{ if }  & a\leq b;\\
                        2b-a & \hbox{ if }  & b\leq a. %
\er\right.\ee
In the next section we shall make use of the `crystal energy' of 
a path $\ket{p}\in P_{a,b}$, which we denote by $h(p)$ and define as
\be
h(p)=-\sli_{s=1}^{\infty} s\!\! &\Big(
 h_1\big(p(2s+1),p(2s),p(2s-1)\big)
-h_1\big(\bar{p}(2s+1),\bar{p}(2s),\bar{p}(2s-1)\big)\nn\\
&+ 
h_2\big(p(2s+2),p(2s+1),p(2s)\big)
-h_2\big(\bar{p}(2s+2),\bar{p}(2s+1),\bar{p}(2s)\big)\nn\\[3mm]
&\!\!+ 
2 h_3\big(p(2s+1),p(2s)\big)
-2 h_3\big(\bar{p}(2s+1),\bar{p}(2s)\big)\Big).\label{hp}\ee
Here, we have abbreviated $\bar{p}(\tts;a,b)$ to $\bar{p}(\tts)$.


\setcounter{equation}{0}
\section{The Path Space $P_{a,b}$}\label{sec3}

The path space $P_{a,b}$ was defined by \eqref{gs1}--\eqref{gs2} in the
previous section.
We shall now go on to consider this space in more detail. In particular, 
we wish to understand the action of $\uq$ on $P_{a,b}$  in the limit
$q\ra 0$. The theory which systematically describes the $q\ra 0$ limit
of $\uq$ was developed by Kashiwara and others, and is known as the theory
of crystal bases~\cite{K91,KMN92,DJO93}.
The main content of this section is a proof of the crystal isomorphism
$P_{a,b}\simeq B(\la_a^{(m-n)})\ot B(\la_b^{(n)})$.
Here, $\la_j^{(k)}=(k-j)\La_0 +j\La_1$,
$j\in\{0,1,\cdots,k\}$, is a level~$k$ dominant integral weight and
$B(\la_j^{(k)})$ is the crystal associated with the highest weight
module $V(\la_j^{(k)})$ (see~\cite{K91}).
We shall use  a principal grading operator $D$,
defined on $B(\lambda^{(m-n)}_a)\otimes B(\lambda^{(n)}_b)$ by
\begin{equation}
D=-\rho+(\rho,\lambda^{(m-n)}_a+\lambda^{(n)}_b),
\end{equation}
where $\rho=\Lambda_0+\Lambda_1$ and $(\cdot,\cdot)$ is the symmetric bilinear
form used in \cite{JM}.
We denote by $B^{(k)}$ the crystal of the $k+1$ dimensional
$U'_q(\widehat{sl}_2)$ module $(V_k)_z$ with $z=1$.
Set $\s\la^{(k)}_j=j\La_0+(k-j)\La_1$.

We give two proofs.
The first makes use of a relation between our models and the
fusion RSOS models.
The second proceeds by examining the crystal isomorphism
\begin{equation}
B(\la_a^{(m-n)})\ot B(\la_b^{(n)}) \simeq B(\s^N \la_a^{(m-n)})\ot
B(\la_b^{(n)}) \ot (B^{(m)}\ot B^{(n)})^{\ot N},
\end{equation}
where $N\in \Z_{>0}$.

\subsection{Identification of $P_{a,b}$ with the tensor product
            of crystals with highest weights}
Let us give the rules for the crystal action of
$\tf_i,\te_i$ ($i=0,1$) on a path $\ket{p}\in P_{a,b}$ 
(for the definition of $\tf_i,\te_i$ and for a detailed discussion of 
the theory of crystal bases, see~\cite{K91,KMN92,Can}):
First, for each $\tts>0$, 
replace each $p(\tts)$ by the sequence of $1$'s and $0$'s
\begin{equation}\label{eq3.1}
p(\tts)\ra \underbrace{1\cdots 1}_{\# 1}\;\underbrace{0\cdots 0}_{\# 0}\; ,
\end{equation}
where
\begin{equation}
(\# 1, \# 0) = 
\begin{cases}
        ({n-p(\tts)},{p(\tts)}) &\text{for $i=0$, $\tts$ odd},\\
        ({m-p(\tts)},{p(\tts)}) &\text{for $i=0$, $\tts$ even},\\
        ({p(\tts)},{n-p(\tts)}) &\text{for $i=1$, $\tts$ odd},\\
        ({p(\tts)},{m-p(\tts)}) &\text{for $i=1$, $\tts$ even}.
\end{cases}
\end{equation}
Then, remove repeatedly all occurrences of adjacent $01$
pairs until we have a sequence of the form $1\cdots1\,0\cdots0$.
On the remaining sequence, use the rule
\bea
        &&\tf_i (\underbrace{1\cdots 1}_{j}\underbrace{0\cdots 0}_{k}) = 
        \underbrace{1\cdots 1}_{j+1}\underbrace{0\cdots 0}_{k-1},\\
        &&\te_i (\underbrace{1\cdots 1}_{j}\underbrace{0\cdots 0}_{k}) = 
        \underbrace{1\cdots 1}_{j-1}\underbrace{0\cdots 0}_{k+1}.
\ena
Finally, put the $01$ pairs back into their original positions and
rebuild the modified path using the inverse of the replacement given
in~\eqref{eq3.1}.

If we remove $01$ from the sequence $\{\bar p(\tts;a,b)\}_{\tts\geq \tts_0}$,
then for $i=1$, we get the sequences,
$\underbrace{0\cdots0}_{a+b}$ (if $\tts_0\equiv1\pmod4$),
$\underbrace{0\cdots0}_{n+a-b}$ (if $\tts_0\equiv2\pmod4$),
$\underbrace{0\cdots0}_{m-n-a+b}$ (if $\tts_0\equiv3\pmod4$),
$\underbrace{0\cdots0}_{m-a-b}$ (if $\tts_0\equiv4\pmod4$).
In all cases, $\te_1$ annihilates the sequence. For $i=0$,
the same is true with the replacement of $a$ by $m-n-a$,
and of $b$ by $n-b$.

For $\tts\geq1$, we use the notation
\be
{\rm wt}_\tts(p)=\left\{ \begin{array}{ll}
(n-2p)(\La_1-\La_0)& \hb{ if }\tts \hb{ is odd;}\\
(m-2p)(\La_1-\La_0)& \hb{ if }\tts \hb{ is even}.
\er\right.
\ee

\begin{dfn}
A path $|p\rangle\in P_{a,b}$ is called admissible if
the sequence of weights $\{\la(\tts)\}_{\tts\geq1}$ defined by
\bea
\la(\tts+1)+{\rm wt}_\tts(p(\tts))&=&\la(\tts),\\
\la(1)&=&\la^{(m)}_{a+b}+\sum_{\tts\geq1}
\Bigl({\rm wt}_\tts(p(\tts))-{\rm wt}_\tts(\bar p(\tts;a,b))\Bigr),
\ena
satisfies
\bea\label{ADM}
\langle h_1,\la(\tts+1)\rangle\geq p(\tts),\quad
\langle h_0,\la(\tts+1)\rangle\geq 
\left\{\br{ll}
n-p(\tts)& \hb{ if } s \hb{ is odd};\\ m-p(\tts)& \hb{ if }s \hb{ is even}.
\er\right.
\ena
\end{dfn}
If $|p\rangle$ is admissible, we have
$\la(\tts)\in\{\la^{(m)}_j;0\leq j\leq m\}$.
Note that the path $\ket{\bar{p}}\in P_{a,b}$ is
admissible, and that the corresponding sequence of weights
is given by the period $4$ repetition
\be
\cdots \la^{(m)}_{m-a-b} \; \la^{(m)}_{m-n-a+b} \;
\la^{(m)}_{n+a-b} \;\la^{(m)}_{a+b}.\nn
\ee

With these definitions in hand, we can proceed to state and prove 
the following theorem:
\begin{thm}\label{isoprop}
There is a crystal isomorphism
$ P_{a,b} \simeq B(\la^{(m-n)}_a)\ot B(\la^{(n)}_b)$, under which the
principal grading is given by $D\ket{p}=h(p)\ket{p}$.
\end{thm}
\noindent The proof is given after preparing the following lemma.
\begin{lem}\label{hwlem}
A path $\ket{p}\in P_{a,b}$ is highest, i.e.,
$\te_i \ket{p}=0$, for  $i=0,1,$
if and only if it is admissible.
\end{lem}
\begin{proof}
 First note that the tensor product rule for crystals (see~\cite{K91})
implies that if a path $\ket{p}\in P_{a,b}$ is highest, and if we
split the tensor product expression for the path at any arbitrary
point $l$ to write 
\be \ket{p}=\big(\cdots \ot p(l+2) \ot p(l+1) \ot p(l)\big)\ot\big(p(l-1)\ot
p(l-2)\ot \cdots p(1)\big),\ee
then $\big(\cdots \ot p(l+2) \ot p(l+1) \ot p(l)\big)$ must also be highest.

Now suppose $l\gg 0$ such that $p(\tts)=\bar p(\tts;a,b)$ for
$\tts\geq l$.
Then, the reduction of the sequence $(p(\tts))_{\tts\geq l}$ for $i=0,1$
gives rise to $\underbrace{0\cdots0}_{\langle
 h_i,\bar\la(l;a,b)\rangle}$.
Since $\ket{p}$ is highest, the path 
 $(p(\tts))_{\tts\geq l-1}$
must also be highest.
For this to be true, it is necessary and sufficient that 
\bea
\langle h_1,\bar{\la}(l;a,b)\rangle\geq p(l-1)\quad
\langle h_0,\bar{\la}(l;a,b)\rangle\geq
\left\{\br{ll}
n-p(l-1)& \hif l \hev;\\
m-p(l-1)& \hif l \hod.\er\right.
\ena
Namely, we have \eqref{ADM} for $\tts=l-1$. Setting $\la(l-1)=
\bar\la(l;a,b)+{\rm wt}(p(l-1))$ we can repeat this argument.
Continuing in the same way to $\la(l-2),\la(l-3)$, etc.,
we can prove the lemma.
\end{proof}

\noindent\emph{Proof of Theorem~\ref{isoprop}\hskip2pt.}
First let us consider the conditions \eqref{ADM} in more detail.
If we write $\la(\tts)=\la^{(m)}_{a(\tts)}$
(where $a(\tts)\in\{0,1,\cdots m\}$), then
the conditions for $\tts$ odd become
\be a(\tts+1)+(n-2p(\tts))&=&a(\tts),\label{c1}\\
    a(\tts+1)&\ge & p(\tts),\label{c2}\\
    m-a(\tts+1) &\ge & n-p(\tts).\label{c3}\ee 
Eliminating $p(\tts)$, we find
\be & a(\tts+1)-a(\tts) & \in \{-n,-n+2,\cdots,n\},\label{c4}\\
    n\leq &a(\tts+1)+a(\tts)&\leq 2m-n.\label{c5}\ee
On the other hand, the admissibility conditions for $\tts$ even become
\be a(\tts+1)+(m-2p(\tts))&=&a(\tts),\label{c6}\\
    a(\tts+1)&\ge & p(\tts),\label{c7}\\
    m-a(\tts+1) &\ge & m-p(\tts).\label{c8}\ee 
Eliminating $p(\tts)$  gives just \be a(\tts+1)=m-a(\tts).\label{c9}\ee

 From these considerations, it follows that an admissible sequence
of weights can be written in the form
\be \cdots  \s(\la^{(m)}_{r(5)})\quad \s(\la^{(m)}_{r(4)})
\quad \la^{(m)}_{r(4)} \quad \la^{(m)}_{r(3)} \quad \s(\la^{(m)}_{r(3)})
\quad \s(\la^{(m)}_{r(2)}) \quad \la^{(m)}_{r(2)} \quad\la^{(m)}_{r(1)},\ee
where the path $\ket{r}=\cdots r(4) \, r(3)\, r(2) \, r(1)$ lies in the space
$R_{a,b}$, defined as the set of paths for which
\be
&& r(\tts) \in\{0,\cdots,m\},\nn\\
&&r(\tts+1)-r(\tts)\in \{-n, -n+2, \cdots, n\},\label{rsosrest}\\
&& n\leq r(\tts+1)+r(\tts) \leq 2m-n,\nn\\
&&r(\tts)=\bar{r}(\tts;a,b),\ws
\tts\gg 0,\ws\hb{where}\label{cab}
\\[3mm]
&&\bar{r}(\tts;a,b)=\left\{\br{ll} 
a+b & \tts \hb{ odd};\\
a+n-b & \tts \hb{ even}.\er \right.\ee
That is, we can identify an admissible path $\ket{p}\in P_{a,b}$ with a 
path $\ket{r}\in R_{a,b}$ by defining
\be r(\tts)= \left\{ \br{lll} 
                      &a(2\tts-1)&\hb{ if } \tts \hb{ is odd;}\\
                      &m-a(2\tts-1)&\hb{ if } \tts \hb{ is even.}
                  \er\right.
\label{rdef}\ee
The restrictions on $\ket{r}\in R_{a,b}$ are those on the space of 
states of the
$\uq$ fusion RSOS models. In~\cite{JMO93}, such a
model is labelled by two integers $(\ell,N)$ and by the level $k=\ell+N$. 
The connection with our notation is that $(\ell,N,k) \leftrightarrow
(m-n,n,m)$. 

Let us denote by $\Omega\big(B(\la_a^{(m-n)})\ot B(\la_b^{(n)})\big)$
the space of highest weight elements in
the crystal $B(\la_a^{(m-n)})\ot B(\la_b^{(n)})$.
Then, it is a well-known theorem that
\be
R_{a,b}&\simeq & \Omega\big(B(\la_a^{(m-n)})\ot B(\la_b^{(n)})\big),
\quad \hb{where}\\
D \ket{r} &=&
  \left( \half \sli_{\tts>0} \tts |r(\tts+2)-r(\tts)|\right) \ket{r}.%
\label{DDD}
\ee
This theorem appears at least implicitly in many of the original works
on RSOS models (see~\cite{ABF,DJKMO87} for example). A statement
and proof using the language of crystals is given in~\cite{JMMO91}.

We now require two
lemmas concerning the crystal energy $h(p)$ of a path in $P_{a,b}$. 
\begin{lem}\label{cren}
The crystal energy of an admissible path $\ket{p}\in P_{a,b}$ is given by
\be\label{RSOS}
h(p)=\half \sum_{\tts=1}^\infty \tts{|r(\tts+2)-r(\tts)|}.
\label{wt1}\ee
\end{lem}
\begin{proof}
The crystal energy $h(p)$ of a path is defined by \eqref{h1}--\eqref{hp}.
If $\ket{p}$ is admissible we have 
\be
a(2\tts+2)+n-2p(2\tts+1) &=& a(2\tts+1),\\
a(2\tts+1)+m-2p(2\tts) &=& a(2\tts).
\ee
Adding these equations and using \eqref{c5} and \eqref{c9} gives
$n\leq p(2\tts+1) +p(2\tts) \leq m$. Hence $h_1,h_2$ and $h_3$ are given
by 
\be
h_1(p(2\tts+1),p(2\tts),p(2\tts-1))&=&
   \{p(2\tts-1)+n-p(2\tts+1)\}_n,\\
h_2(p(2\tts+2),p(2\tts+1),p(2\tts))&=&
   \{p(2\tts+2)+p(2\tts)+2 p(2\tts+1)-n\}_m,\\
h_3(p(2\tts+1),p(2\tts))&=& n.
\ee
Using \eqref{c1}--\eqref{c9} it is simple to show that
\be p(2\tts-1)+n-p(2\tts+1) &=& n+ a(2\tts+3) -a(2\tts-1),\\
p(2\tts+2)+p(2\tts)+2 p(2\tts+1)-n&=& m.\ee
Writing $a(2\tts-1)$ in terms of $r(\tts)$ using \eqref{rdef} then gives
\be
h_1(p(2\tts+1),p(2\tts),p(2\tts-1))= n-\half|r(\tts+2)-r(\tts)|,\ee
which completes the proof.
\end{proof}

\begin{lem}\label{pmo}
The action of $\tf_i$ on a path increases $h(p)$ by $1$, 
and that of $\te_i$ decreases it by $1$.
\end{lem}
\begin{proof}
$\tf_1$ acts on a path $\ket{p}\in P_{a,b}$ by changing 
$p(\tts)\ra p(\tts)+1$ at a single value of $\tts$. Suppose this happens for
$\tts=2l+1$.
Then the following 
inequalities must hold:
\be
&&n-p(2l+1) > p(2l), \\
&&m-p(2l+2)\leq p(2l+1).
\ee

Using these inequalities, \eqref{h1}--\eqref{h3}, and the property
\be
\{a+1\}_n -\{a\}_n=\left\{\begin{array}{ll}
                1 \quad &\hbox{ if } a < n,\\
               -1 \quad &\hbox{ if } a \geq n,%
\end{array}\right.
\ee 
we arrive at 
\be 
h_1(p(2l+3),p(2l+2),p(2l+1)+1) -  
h_1(p(2l+3),p(2l+2),p(2l+1)) &=& -1, \nn\\
h_2(p(2l+2),p(2l+1)+1,p(2l)) - 
 h_2(p(2l+2),p(2l+1),p(2l)) &=& -1, \nn\\
h_1(p(2l+1)+1,p(2l),p(2l-1)) -  
h_1(p(2l+1),p(2l),p(2l-1)) &=& 0,\nn \\
h_3(p(2l+1)+1,p(2l)) -  h_3(p(2l+1),p(2l)) &=& 1.\nn
\ee
 From this we see that $h(p)\ra h(p)+1$
when we change $p(2l+1)\rightarrow p(2l+1)+1$.
The proofs for the case when $\tts=2l$, and for
$\tf_0,\te_0,\te_1$ are similar.
\end{proof}

We have shown that the space of highest paths
in $P_{a,b}$ is isomorphic to the
space $R_{a,b}$ which is in turn isomorphic
to $\Omega\big(B(\la_a^{(m-n)})\ot B(\la_b^{(n)})\big)$. Combining this
result with \eqref{DDD}, and Lemmas \ref{cren} and
\ref{pmo} completes the proof of Theorem~\ref{isoprop}. \hfill $\square$

\subsection{The crystal structure of the path space}

In this subsection, we will give another proof of the crystal isomorphism
$P_{a,b} \cong B(\la) \ot B(\mu)$ between the path
space $P_{a,b}$ and the tensor product of the crystals 
$B(\la) \ot B(\mu)$, where 
$\la=a\La_{1} + (m-n-a) \La_{0}$ and
$\mu=b\La_{1} + (n-b) \La_{0}$. 
We first recall some of the fundamental results on the crystals for
the quantum affine algebra $\uq$ (cf.~\cite{KMN92}).

Let $\ell>0$ be a positive integer and let 
$B^{(\ell)}=\{[j]^{(\ell)} \, | \, 0 \leq j \leq \ell \}$ be the perfect
crystal of level $\ell$ for the quantum affine algebra $\uq$.
The crystal structure of $B^{(\ell)}$ is given in the following:

\begin{center}
\begin{texdraw}
\fontsize{10}{13}\selectfont
\setunitscale 1
\drawdim mm
\arrowheadsize l:2.4 w:1.1 \arrowheadtype t:F
\textref h:C v:C
\htext(0 0){$[0]^{(\ell)}$}
\htext(18 0){$[1]^{(\ell)}$}
\htext(70 0){$[\ell-1]^{(\ell)}$}
\htext(90 0){$[\ell]^{(\ell)}$}
\move(4 1)\ravec(9 0)\rmove(0 -2)\ravec(-9 0)
\move(22 1)\ravec(9 0)\rmove(0 -2)\ravec(-9 0)
\move(77 1)\avec(85 1)\rmove(0 -2)\avec(77 -1)
\move(62 -1)\ravec(-9 0)\rmove(0 2)\ravec(9 0)
\htext(8 3){$1$}\htext(9 -3){$0$}
\htext(26 3){$1$}\htext(27 -3){$0$}
\htext(80 3){$1$}\htext(81 -3){$0$}
\htext(57 3){$1$}\htext(58 -3){$0$}
\htext(39.3 0){$\cdots$}
\htext(44 0){$\cdots$}
\move(-4 4)\move(94 -5)
\end{texdraw}
\end{center}

\noindent
We will write $[j]$ in place of $[j]^{(\ell)}$ whenever there is 
no danger of confusion. 

The following theorem gives one of the most fundamental isomorphisms
in the theory of crystals for the quantum affine algebra $\uq$.

\begin{thm}[cf.~\cite{KMN92}] \label{perfect}
For any dominant integral weight $\la=s\La_1 +(\ell-s) \La_0$
of level $\ell>0$, there exists a crystal isomorphism 
\begin{equation} 
\psi=\psi_{\la} : B(\la) \buildrel {\sim} \over \longrightarrow  
B(\s \la) \ot B^{(\ell)}
\end{equation}
such that 
\begin{equation}
u_{\la} \mapsto u_{\s \la} \ot [\ell-s],\label{EQU2}
\end{equation}
where $u_{\la}$ is the highest weight element of $B(\la)$. 
\end{thm}

Let $ | p_{\la} \rangle 
=(p_{\la}(\tts))_{\tts=1}^{\infty}$ be the sequence of elements 
in $B^{(\ell)}$ whose terms are given by
\begin{equation}
p_{\la}(\tts)= \begin{cases}  \  [\ell-s] \ \ &  \text {if $\tts$ is odd}, \\
\ [s] \ \ & \text {if $\tts$ is even}.
\end{cases}
\end{equation}
For each positive integer $N>0$,
there exists a crystal isomorphism 
\begin{equation}
\psi^{(N)} = \psi_{\la}^{(N)} :B(\la) \buildrel {\sim} \over 
\longrightarrow B(\s^N \la) 
\ot (B^{(\ell)})^{\ot N} 
\end{equation}
such that
\begin{equation}
u_{\la} \mapsto u_{\s^N \la} \ot p_{\la}(N)
\ot \cdots \ot p_{\la}(2) \ot p_{\la}(1).
\end{equation}

A sequence
$| p \rangle =(p(\tts))_{\tts=1}^{\infty}$ with $p(\tts) \in B^{(\ell)}$ is
called a \emph{$\la$-path in $B^{(\ell)}$} if $p(\tts)=p_{\la}(\tts)$
for all sufficiently large $\tts$. Let $P(\la)$ be
the set of all $\la$-paths in $B^{(\ell)}$.
Each $\la$-path
$| p \rangle=(p(\tts))_{\tts=1}^{\infty}$  is understood as the half-infinite
tensor product 
$|p \rangle = \cdots \ot p(\tts+1) \ot p(\tts) \ot \cdots \ot p(2) 
\ot p(1)$ and hence the set $P(\la)$ is given
a crystal structure for the quantum affine algebra $\uq$ by the 
tensor product rule for the crystals. 

Moreover, one can prove:

\begin{thm}[cf.~\cite{KMN92}] \label{KMN}
For each $b\in B(\la)$, there exists a positive integer $N>0$
such that 
\begin{equation}
\psi^{(N)}(b) \in u_{\s^N \la} \ot (B^{(\ell)})^{\ot N}.
\end{equation}
Hence we have the crystal isomorphism 
$B(\la) \buildrel{\sim} \over \longrightarrow  P(\la).$
\end{thm}

In~\cite{Nak96b}, in his study of 6 vertex models of
inhomogeneous type,  Nakayashiki considered the crystal
isomorphism $\psi$ in a more general setting. 

\begin{thm}[cf.~\cite{Nak96b,kaka97}]\label{Nak}
Let $\mu=b \La_1
+ (n-b) \La_0$ be a dominant integral weight of level $n>0$ 
and let $B^{(k)}$ be the perfect crystal of level $k>0$ for the
quantum affine algebra $\uq$. Then there exists a crystal 
isomorphism 
\begin{equation} 
\Psi: B^{(k)} \ot B(\mu) \buildrel {\sim} \over \longrightarrow 
B(\s \mu) \ot B^{(n+k)}
\end{equation}
such that
\begin{equation}
[j]^{(k)} \ot u_{\mu} \mapsto u_{\s \mu} \ot 
[j+n-b]^{(n+k)}.\label{MAP1}
\end{equation}
\end{thm}

Suppose $m> n$ and let $\la=a \La_{1} + (m-n-a) \La_{0}$
be a dominant integral weight of level $m-n$. Define a crystal 
isomorphism 
\begin{equation} \label {phi}
\begin{aligned}
\ \Phi =\Phi_{\la, \mu} & = \Psi \circ (\psi_{\la} \ot \psi_{\mu})
= (\id\ot\Psi\ot\id) \circ 
(\psi_{\la} \ot\id\ot\id) \circ (\id\ot\psi_{\mu}) : \\
\ & B(\la) \ot B(\mu) \xrightarrow {\id\ot \psi}
B(\la) \ot B(\s \mu) \ot B^{(n)}  \\
& \hskip 1cm \xrightarrow{\psi \ot \id \ot\id}
B(\s \la) \ot B^{(m-n)} \ot B(\s \mu) 
\ot B^{(n)} \\
& \hskip 1cm  \xrightarrow {\id \ot \Psi \ot \id} 
B(\s \la) \ot B(\mu) \ot B^{(m)}\ot B^{(n)}
\end{aligned}
\end{equation}
to be the composite of the crystal isomorphisms defined 
in Theorem~\ref{perfect} and Theorem~\ref{Nak}. 

By repeating the crystal isomorphism $\Phi$, we obtain the crystal 
isomorphism
$$\Phi^{(2)}: B(\la) \ot B(\mu) \buildrel {\sim} \over \longrightarrow
B(\la) \ot B(\mu) \ot (B^{(m)} \ot B^{(n)})^{\ot 2}$$
such that
$$u_{\la} \ot u_{\mu} \mapsto u_{\la} \ot u_{\mu} 
\ot [a+b]^{(m)} \ot [n-b]^{(n)} 
\ot [m-n-a+b]^{(m)} \ot [n-b]^{(n)}.$$

In general, let 
$|p_{\la, \mu} \rangle=(p_{\la, \mu}(\tts))_{\tts=1}^{\infty}$ be the 
sequence of elements in $B^{(m)}$ and $B^{(n)}$ whose terms are given by
\begin{equation} 
p_{\la, \mu}(\tts) =
\begin{cases}
{}[n-b] & \text{if $\tts$ is odd},\\
{}[a+b] & \text{if $\tts \equiv 0 \pmod 4$},\\
{}[m-n-a+b] & \text{if $\tts \equiv 2 \pmod 4$}.
\end{cases}
\end{equation}
For each positive integer $N>0$, we have a crystal isomorphism 
\begin{equation}
\Phi^{(N)}=\Phi_{\la, \mu}^{(N)}:
B(\la) \ot B(\mu) \buildrel {\sim} \over
\longrightarrow B(\s^N \la) \ot B(\mu) \ot 
(B^{(m)} \ot B^{(n)})^{\ot N}
\end{equation}
such that 
\begin{equation} 
u_{\la} \ot u_{\mu} \mapsto u_{\s^N \la} \ot u_{\mu} 
\ot p_{\la, \mu}(2N) \ot p_{\la, \mu}(2N-1) 
\ot \cdots \ot p_{\la, \mu}(2) \ot p_{\la, \mu}(1).
\end{equation}

A sequence $|p \rangle =(p(\tts))_{\tts=1}^{\infty}$ of elements
in $B^{(m)}$ and $B^{(n)}$ is called a \emph{$(\la, \mu)$-path}
if $p(\tts)=p_{\la, \mu}(\tts)$ for all sufficiently large $\tts$.
Let $P(\la, \mu)$ be the set of all $(\la, \mu)$-paths.
The crystal structure of $P(\la, \mu)$ is the
same as that of the path space $P_{a,b}$. 

Now, we would like to show that there exists a crystal isomorphism 
$B(\la) \ot B(\mu) \buildrel {\sim} \over \longrightarrow 
P(\la, \mu)$. As in the case with the crystal 
isomorphism $B(\la) \buildrel {\sim} \over \longrightarrow
P(\la)$, it suffices to prove that for each 
$v \ot w \in B(\la) \ot B(\mu)$, there exists a positive
integer $N>0$ such that
\begin{equation} \label{claim}
\Phi^{(N)}(v \ot w) \in u_{\s^N \la} \ot u_{\mu} 
\ot (B^{(m)} \ot B^{(n)})^{\ot N}.
\end{equation}

For this purpose, we need an explicit description of the crystal
isomorphism $R:B^{(m)} \ot B^{(n)} \isomo B^{(n)} \ot 
B^{(m)}$, called the \emph{combinatorial $R$-matrix}, normalised by the 
condition $R([0]^{(m)} \ot [0]^{(n)})=[0]^{(n)} \ot [0]^{(m)}$.
We rephrase \eqref{zerot} as follows.

\begin{lem} \label{R-matrix}
The normalised combinatorial $R$-matrix 
$R:B^{(m)} \ot B^{(n)} \isomo B^{(n)} \ot B^{(m)}$ is given by
\begin{equation}\label{EQU3}
R([i]^{(m)} \ot [j]^{(n)}) =
\begin{cases}
\ [i]^{(n)} \ot [j]^{(m)} \ \ & \text {if} \ \ 
i+j \leq m, n, \\
\ [2i+j-m]^{(n)} \ot [m-i]^{(m)} \ \ & \text {if} 
\ \ m \leq i+j \leq n, \\
\ [n-j]^{(n)} \ot [i+2j-n]^{(m)} \ \ & \text {if} 
\ \ n \leq i+j \leq m, \\
\ [n-m+i]^{(n)} \ot [j+m-n]^{(m)} \ \ & \text {if} 
\ \ i+j \geq m, n.
\end{cases}
\end{equation}
\end{lem}

The following lemma plays a crucial role in proving our isomorphism 
theorem.

\begin{lem}
Let $\mu=b \La_1 + (n-b) \La_0$. 
If $m> n$, the following diagram of crystal isomorphisms is commutative.

\begin{center}
\begin{texdraw}
\fontsize{10}{13}\selectfont
\setunitscale 1
\drawdim mm
\arrowheadsize l:2.4 w:1.1 \arrowheadtype t:F
\textref h:C v:C
\htext(0 0){$B^{(m-n)}\ot B(\mu)$}
\htext(0 -22){$B(\s\mu)\ot B^{(m)}$}
\htext(0 -44){$B(\mu)\ot B^{(n)}\ot B^{(m)}$}
\htext(70 0){$B^{(m-n)}\ot B(\s\mu) \ot B^{(n)}$}
\htext(70 -44){$B(\mu)\ot B^{(m)}\ot B^{(n)}$}
\move(14 0)\avec(50 0) \move(0 -3)\avec(0 -19) \move(0 -26)\avec(0 -41)
\move(70 -3)\avec(70 -41) \move(53 -44)\avec(17 -44)
\htext(-3 -11){$\Psi$}
\htext(-7 -33){$\psi\ot\id$}
\htext(32 2.5){$\id\ot\psi$}
\htext(35 -41.5){$R$}
\htext(77 -22){$\Psi\ot\id$}
\end{texdraw}
\end{center}
\end{lem}
\begin{proof}
Since the crystal $B^{(m-n)} \ot B(\mu)$ is connected
(see~\cite{Nak96b, kaka97}), it suffices to check the commutativity for a 
single element, say $[j]^{(m-n)} \ot u_{\mu} \in B^{(m-n)} \ot B(\mu)$
with $0 \leq j \leq m-n$.
Using \eqref{MAP1}, \eqref{EQU2}, and \eqref{EQU3}, we can show that
$[j]^{(m-n)} \ot u_{\mu}$
is mapped to $u_{\mu} \ot [b]^{(n)} \ot [j+n-b]^{(m)}$ in both ways.
\end{proof}

By applying the above lemma repeatedly, we obtain:

\begin{coro}\label{comm}
For each positive integer $N>0$, the following diagram of
crystal isomorphisms is commutative.

\begin{center}
\begin{texdraw}
\fontsize{10}{13}\selectfont
\setunitscale 1
\drawdim mm
\arrowheadsize l:2.4 w:1.1 \arrowheadtype t:F
\textref h:C v:C
\move(0 0)
\bsegment
\htext(0 0){$B(\la)\ot(B^{(m-n)})^{\ot N}\ot B(\mu)$}
\esegment
\move(75 0)
\bsegment
\htext(-2 2){$B(\s\la)\ot(B^{(m-n)})^{\ot (N+1)}$}
\htext(5 -3){$\ot B(\s\mu)\ot B^{(n)}$}
\esegment
\move(75 -50)
\bsegment
\htext(-4 2){$B(\s\la)\ot B^{(m-n)}\ot B(\s^{N+1}\mu)$}
\htext(4 -3){$\ot(B^{(m)})^{\ot N}\ot B^{(n)}$}
\esegment
\move(0 -50)
\bsegment
\htext(-4 2){$B(\s\la)\ot B^{(m-n)}\ot B(\s^{N+1}\mu)$}
\htext(4 -3){$\ot B^{(n)}\ot(B^{(m)})^{\ot N}$}
\esegment
\move(0 -25)
\bsegment
\htext(0 0){$B(\la)\ot B(\s^N\mu)\ot(B^{(m)})^{\ot N}$}
\esegment
\htext(34 -48){$R\circ\cdots\circ R$}
\htext(35 2.3){$\psi_\la\ot\psi_\mu$}
\htext(-9.5 -12){$\Psi\circ\cdots\circ\Psi$}
\htext(-9.5 -38){$\psi_\la\ot\psi_{\s^N\mu}$}
\htext(80 -25){$\Psi\circ\cdots\circ\Psi$}
\move(23 0)\avec(50 0)
\move(47 -50)\avec(20 -50)
\move(0 -4)\avec(0 -23)
\move(0 -29)\avec(0 -45)
\move(70 -7)\avec(70 -45)
\move(-20 4)\move(91 -56)
\end{texdraw}
\end{center}

\end{coro}

Let $v \ot w \in B(\la) \ot B(\mu)$. By Theorem~\ref{KMN},
there exists a positive integer $N>0$ such that 
$\psi_{\la}^{(N)} (v) \in u_{\s^N \la} \ot
(B^{(m-n)})^{\ot N}$ and 
$\psi_{\mu}^{(N)} (v) \in u_{\s^N \mu} \ot
(B^{(n)})^{\ot N}$.
Hence we obtain the crystal isomorphism 
\begin{equation}
\begin{aligned}
B(\la) \ot B(\mu) & \xrightarrow {\psi_{\la}^{(N)}
\ot \psi_{\mu}^{(N)}} 
B(\s^N \la) \ot (B^{(m-n)})^{\ot N}
\ot B(\s^N \mu) \ot (B^{(n)})^{\ot N} \\
& \xrightarrow {\Psi \circ \cdots \circ \Psi} 
B(\s^N \la) \ot B(\mu)
\ot (B^{(m)})^{\ot N} \ot (B^{(n)})^{\ot N} \\
& \xrightarrow {R \circ \cdots \circ R} 
B(\s^N \la) \ot B(\mu)
\ot (B^{(m)} \ot B^{(n)})^{\ot N}
\end{aligned}
\end{equation}
such that 
$v \ot w$ is mapped to an element in
$u_{\s^N \la} \ot u_{\mu} \ot 
(B^{(m)} \ot B^{(n)})^{\ot N}$.

Therefore, in order to prove our claim \eqref{claim}, it suffices to prove
that the following diagram of crystal isomorphisms is commutative.

\begin{center}
\begin{texdraw}
\fontsize{10}{13}\selectfont
\setunitscale 1
\drawdim mm
\arrowheadsize l:2.4 w:1.1 \arrowheadtype t:F
\textref h:C v:C
\htext(0 0){$B(\la)\ot B(\mu)$}
\move(70 0)
\bsegment
\htext(-0.5 2){$B(\s^N\la)\ot(B^{(m-n)})^{\ot N}$}
\htext(0.5 -3){$\ot B(\s^N\mu)\ot(B^{(n)})^{\ot N}$}
\esegment
\move(0 -25)
\bsegment
\htext(-1.6 2){$B(\s^N\la)\ot B(\mu)$}
\htext(1.6 -3){$\ot (B^{(m)}\ot B^{(n)})^{\ot N}$}
\esegment
\move(70 -25)
\bsegment
\htext(-3.5 2){$B(\s^N\la)\ot B(\mu)$}
\htext(3.5 -3){$\ot (B^{(m)})^{\ot N}\ot(B^{(n)})^{\ot N}$}
\esegment
\htext(35 -22.5){$R\circ\cdots\circ R$}
\htext(33 3.2){$\psi^{(N)}_\la\ot\psi^{(N)}_\mu$}
\htext(-5 -12){$\Phi^{(N)}$}
\htext(80 -12){$\Psi\circ\cdots\circ\Psi$}
\move(12 0)\avec(50 0)
\move(53 -25)\avec(17 -25)
\move(0 -3)\avec(0 -20)
\move(70 -6)\avec(70 -20)
\move(-14 4)\move(85 -32)
\end{texdraw}
\end{center}

We will prove our assertion by induction on $N$. 
If $N=1$, there is nothing to prove. 
Assume that our assertion is true for $N-1$. 
Then by the induction hypothesis, the following diagram of 
crystal isomorphisms is commutative.

\begin{center}
\begin{texdraw}
\fontsize{10}{13}\selectfont
\setunitscale 1
\drawdim mm
\arrowheadsize l:2.4 w:1.1 \arrowheadtype t:F
\textref h:C v:C
\htext(0 0){$B(\la)\ot B(\mu)$}
\rmove(0 -25)
\bsegment
\htext(-2 2){$B(\s^{N-1}\la)\ot B(\mu)$}
\htext(2 -3){$\ot (B^{(m)}\ot B^{(n)})^{\ot(N-1)}$}
\esegment
\rmove(0 -25)
\bsegment
\htext(-4.5 2){$B(\s^N\la)\ot B^{(m-n)}\ot B(\s\mu)\ot B^{(n)}$}
\htext(4.5 -3){$\ot(B^{(m)}\ot B^{(n)})^{\ot(N-1)}$}
\esegment
\rmove(0 -25)
\bsegment
\htext(-2 2){$B(\s^N\la)\ot B(\mu)$}
\htext(1 -3){$\ot(B^{(m)}\ot B^{(n)})^{\ot N}$}
\esegment
\move(80 0)
\bsegment
\htext(0 2){$B(\s^{N-1}\la)\ot(B^{(m-n)})^{\ot(N-1)}$}
\htext(0.5 -3){$\ot B(\s^{N-1}\mu)\ot(B^{(n)})^{\ot(N-1)}$}
\esegment
\rmove(0 -25)
\bsegment
\htext(-5.5 2){$B(\s^{N-1}\la)\ot B(\mu)$}
\htext(5.5 -3){$\ot(B^{(m)})^{\ot(N-1)}\ot(B^{(n)})^{\ot(N-1)}$}
\esegment
\rmove(0 -25)
\bsegment
\htext(-1.5 2){$B(\s^N\la)\ot B^{(m-n)}\ot B(\s\mu)\ot B^{(n)}$}
\htext(1.5 -3){$\ot(B^{(m)})^{\ot(N-1)}\ot(B^{(n)})^{\ot(N-1)}$}
\esegment
\rmove(0 -25)
\bsegment
\htext(-0.5 2){$B(\s^N\la)\ot B(\mu)\ot B^{(m)}\ot B^{(n)}$}
\htext(0.5 -3){$\ot(B^{(m)})^{\ot(N-1)}\ot(B^{(n)})^{\ot(N-1)}$}
\esegment
\htext(35 3){$\psi_\la^{(N-1)}\ot\psi_\mu^{(N-1)}$}
\htext(39 -23){$R\circ\cdots\circ R$}
\htext(38 -48){$R\circ\cdots\circ R$}
\htext(35 -73){$R\circ\cdots\circ R$}
\htext(-7 -12){$\Phi^{(N-1)}$}
\htext(-11 -37){$\psi_{\s^{N-1}\la}\ot\psi_\mu$}
\htext(-3 -62){$\Psi$}
\htext(90 -12){$\Psi\circ\cdots\circ\Psi$}
\htext(92 -37){$\psi_{\s^{N-1}\la}\ot\psi_\mu$}
\htext(83.5 -62){$\Psi$}
\move(11 0)\avec(56 0)
\move(59 -25)\avec(20 -25)
\move(50 -50)\avec(25 -50)
\move(53 -75)\avec(18 -75)
\move(0 -3)\avec(0 -20)
\move(0 -31)\avec(0 -45)
\move(0 -56)\avec(0 -70)
\move(80 -6)\avec(80 -20)
\move(80 -31)\avec(80 -45)
\move(80 -56)\avec(80 -70)
\move(-15 6)\move(101 -81)
\end{texdraw}
\end{center}

\noindent
Note that the commutativity of the first square follows from the induction
hypothesis and the commutativity of the other squares is trivial.

Next, observe that Corollary \ref{comm} yields the following
commutative diagram of crystal isomorphisms. 

\begin{center}
\begin{texdraw}
\fontsize{10}{13}\selectfont
\setunitscale 1
\drawdim mm
\arrowheadsize l:2.4 w:1.1 \arrowheadtype t:F
\textref h:C v:C
\move(0 0)
\bsegment
\htext(-0.5 2){$B(\s^{N-1}\la)\ot(B^{(m-n)})^{\ot(N-1)}$}
\htext(0.5 -3){$\ot B(\s^{N-1}\mu)\ot(B^{(n)})^{\ot(N-1)}$}
\esegment
\move(0 -25)
\bsegment
\htext(-5.5 2){$B(\s^{N-1}\la)\ot B(\mu)$}
\htext(5.5 -3){$\ot(B^{(m)})^{\ot(N-1)}\ot(B^{(n)})^{\ot(N-1)}$}
\esegment
\move(0 -50)
\bsegment
\htext(-1 2){$B(\s^N\la)\ot B^{(m-n)}\ot B(\s\mu)\ot B^{(n)}$}
\htext(1 -3){$\ot (B^{(m)})^{\ot(N-1)}\ot(B^{(n)})^{\ot(N-1)}$}
\esegment
\move(80 -50)
\bsegment
\htext(-0.5 2){$B(\s^N\la)\ot B^{(m-n)}\ot B(\s\mu)$}
\htext(0.5 -3){$\ot (B^{(m)})^{\ot(N-1)}\ot(B^{(n)})^{\ot N}$}
\esegment
\move(80 0)
\bsegment
\htext(-0.5 2){$B(\s^N\la)\ot(B^{(m-n)})^{\ot N}$}
\htext(0.5 -3){$\ot B(\s^N\mu)\ot(B^{(n)})^{\ot N}$}
\esegment
\htext(-10 -12){$\Psi\circ\cdots\circ\Psi$}
\htext(-11 -37){$\psi_{\s^{N-1}\la}\ot\psi_\mu$}
\htext(43 -48){$R\circ\cdots\circ R$}
\htext(89 -25){$\Psi\circ\cdots\circ\Psi$}
\htext(42 2.5){$\psi_{\s^{N-1}\la}\ot\psi_{\s^{N-1}\mu}$}
\move(24 0)\avec(61 0)
\move(57 -50)\avec(27 -50)
\move(0 -6)\avec(0 -20)
\move(0 -31)\avec(0 -45)
\move(80 -6)\avec(80 -45)
\move(-18 5)\move(99 -56)
\end{texdraw}
\end{center}

Moreover, the commutativity of the following diagram is trivial. 

\begin{center}
\begin{texdraw}
\fontsize{10}{13}\selectfont
\setunitscale 1
\drawdim mm
\arrowheadsize l:2.4 w:1.1 \arrowheadtype t:F
\textref h:C v:C
\move(0 0)
\bsegment
\htext(-1 2){$B(\s^N\la)\ot B^{(m-n)}\ot B(\s\mu)\ot B^{(n)}$}
\htext(1 -3){$\ot(B^{(m)})^{\ot(N-1)}\ot(B^{(n)})^{\ot(N-1)}$}
\esegment
\move(0 -20)
\bsegment
\htext(-0.5 2){$B(\s^N\la)\ot B(\mu)\ot B^{(m)}\ot B^{(n)}$}
\htext(0.5 -3){$\ot(B^{(m)})^{\ot(N-1)}\ot(B^{(n)})^{\ot(N-1)}$}
\esegment
\move(80 0)
\bsegment
\htext(-0.5 2){$B(\s^N\la)\ot B^{(m-n)}\ot B(\s\mu)$}
\htext(0.5 -3){$\ot(B^{(m)})^{\ot(N-1)}\ot(B^{(n)})^{\ot N}$}
\esegment
\move(80 -20)
\bsegment
\htext(-3.5 2){$B(\s^N\la)\ot B(\mu)$}
\htext(3.5 -3){$\ot(B^{(m)})^{\ot N}\ot(B^{(n)})^{\ot N}$}
\esegment
\move(56 0)\avec(28 0)
\move(63 -20.5)\avec(26 -20.5)
\move(0 -6)\avec(0 -15)
\move(80 -6)\avec(80 -15)
\htext(42 2.2){$R\circ\cdots\circ R$}
\htext(44 -18.4){$R\circ\cdots\circ R$}
\htext(-3 -10){$\Psi$}
\htext(83 -10){$\Psi$}
\move(-22 4)\move(103 -26)
\end{texdraw}
\end{center}
By combining all the commutative diagrams obtained above, 
we obtain the desired commutative diagram.

To summarise the above discussion, we obtain:

\begin{thm}
Suppose $m> n$ and let $\la=a \La_1 + (m-n-a) \La_0$
and $\mu=b \La_1 + (n-b) \La_0$ be dominant integral 
weights of level $m-n$ and $n$, respectively. Then there exists a
crystal isomorphism 
\begin{equation}
B(\la) \ot B(\mu) \buildrel {\sim} \over 
\longrightarrow P_{a,b}
\end{equation}
such that 
\begin{equation}
u_{\la} \ot u_{\mu} \longmapsto 
\cdots \ot [a+b] \ot [n-b] \ot [m-n-a+b] \ot [n-b].
\end{equation}
\end{thm}


\setcounter{equation}{0}
\section{Level-1 Intertwiners}\label{sec4}

In this section, we define some $\uq$ intertwiners and study them
in the level-$1$ case.
Commutation relations for these operators will be given and some of their
matrix elements will be calculated.
These will be used in the next section to reduce questions
about general level intertwiners to those of level-$1$ intertwiners.

The reader should recall some notation from Section~\ref{sec21}.
Let $\s$ denote the map exchanging the two fundamental weights
$\La_0\leftrightarrow \La_1$.
The set of dominant integral weights of level $k$ will be denoted by
$P^+_k$.
For $\la\in P^+_k$, $V(\la)$ will be the irreducible highest weight module
of highest weight $\la$.
Also set $\la_\pm = \la \pm (\La_1 - \La_0)$ for $\la\in P^+_k$.
We only consider the cases when $\la_\pm \in P^+_k$.
The following notation for various $\uq$ intertwiners will be used.
\bea
\Ps{\an}{\an+k} & : & \Vz{\an}\ot\vl \ra \vsl\ot\Vz{\an+k},\label{IT1}\\
\Phi^\pm & : & \vl \ra V(\la_\pm)\ot\Vz{1},\label{IT2}\\
\Psi^{*\pm} & : & \Vz{1}\ot\vl \ra V(\la_\pm).\label{IT3}
\ena
As we said in 2.1, we consider the evaluation modules in these formulas as
modules of finite rank over the ring of coefficients containing $\zeta$ and
$\zeta^{-1}$. This implies, in particular, that the intertwiners commute
with the multiplication of $\zeta$ and $\zeta^{-1}$. We remark also that in
each of \eqref{IT1}--\eqref{IT3} we need completion in the right hand side
(see~\cite{DJO93} for a detailed discussion).

As to the existence of the intertwiners \eqref{IT1}, we have the following
proposition.
The case $\an=0$ is given in~\cite{DJO93}, and the case $k=1$ in~\cite{Nak}.
See also~\cite{Nak96b} and~\cite{kaka97} for the crystal version.

\begin{prop}\label{pr:ex}
Let $\la, \nu \in P^+_k$.
There exists a $\uqp$ intertwiner from $\Vz{\an}\ot\vl$ to
$V(\nu)\ot\Vz{\an+k}$ if and only if $\nu = \s(\la)$.
When $\nu = \s(\la)$, the intertwiner is unique up to scalar multiple.
\end{prop}
\begin{proof}
We first prove
\begin{align*}
\Hom_{U^\prime_q}&(V^{(\an)}_\z\ot V(\la), V(\nu)\ot V^{(\an+k)}_\z)\\
&\cong W_\la^\nu
= \{ v\in V^{(\an+k)}_\z\ot V^{(\an)}_{\z q^{-(k+2)/2}};
e_i^{\nu(h_i)+1} v = 0, \operatorname{wt}(v) = \la-\nu \}.
\end{align*}
This is very similar to the proof in~\cite{Nak} which is an extended version
of the case $\an=0$ given in~\cite{DJO93}.
The steps are
\begin{align*}
\Hom_{U_q^\prime}&(V^{(\an)}_\z\ot V(\la), V(\nu)\ot V^{(\an+k)}_\z)\\
&\cong \Hom_{U^\prime_q}
   (V(\la), V^{(\an)*a^{-1}}_\z\ot V(\nu)\ot V^{(\an+k)}_\z)\\
&\cong \Hom_{U^\prime_q(b^+)}
   (\Q(q)u_\la, V^{(\an)*a^{-1}}_\z\ot V(\nu)\ot V^{(\an+k)}_\z)\\
&\cong \Hom_{U^\prime_q(b^+)}
   (V(\nu)^{*a}\ot V^{(\an)}_\z\ot\Q(q)u_\la, V^{(\an+k)}_\z)\\
&\cong \Hom_{U^\prime_q(b^+)}
   (V(\nu)^{*a}\ot \Q(q)u_\la\ot V^{(\an)}_{\z q^{-k/2}}, V^{(\an+k)}_\z)\\
&\cong \Hom_{U^\prime_q(b^+)}
   (V(\nu)^{*a}\ot \Q(q)u_\la, V^{(\an+k)}_\z\ot V^{(\an)*a}_{\z q^{-k/2}})\\
&\cong \Hom_{U^\prime_q(b^+)}
   (V(\nu)^{*a}\ot \Q(q)u_\la, V^{(\an+k)}_\z\ot V^{(\an)}_{\z q^{-(k+2)/2}})\\
&\cong W^\nu_\la.
\end{align*}
Here we used the $U^\prime_q(b^+)$ isomorphism
$V^{(\an)}_\z\ot\Q(q)u_\la\cong \Q(q)u_\la\ot V^{(\an)}_{\z q^{-k/2}}$:
\begin{equation*}
u^{(\an)}_j\ot u_\la\mapsto q^{-j(\la(h_1)-k/2)} u_\la\ot u^{(\an)}_j.
\end{equation*}
To complete the proof, we show $\dim(W_\la^\nu) = \de_{\nu,\s(\la)}$.
Suppose $\nu = j\La_1 + (k-j)\La_0$.
We solve for the vector $v\in V^{(\an+k)}_\z\ot V^{(\an)}_{\z q^{-(k+2)/2}}$
which satisfies $e_1^{j+1}v = 0$, $e_0^{k-j+1}v = 0$.
It is uniquely given up to constant multiple by
\begin{equation*}
y^{(k)}_j = \sum_{i=0}^{\an}(-1)^{i}
            \qbinom{k+i-j}{i}^{\half}\qbinom{\ell-i+j}{j}^{\half}
            u^{(\an+k)}_{\an-i+j}\ot u^{(\an)}_{i}.
\end{equation*}
These vectors span a space isomorphic to $V^{(k)}_{\z q^{\an/2}}$.
The weight of $y^{(k)}_j$ is $(k-2j)(\La_1-\La_0)$ and this is equal to
$\la-\nu$ if and only if $\la = \s(\nu)$.
\end{proof}

Similar existence and uniqueness theorems for the other two intertwiners are
also known.
We define $\ket{\la}$ to be the highest weight vector of $V(\la)$ and
take $\bra{\la}$ to be its dual.
With these, we normalise the intertwiners as follows:
\bea
\lefteqn{\bra{\La_1}\Ps{\an}{\an+1}\ket{\La_0}(u^{(\an)}_{\an}) =
u^{(\an+1)}_{\an+1},}\label{eq:459}\hspace{23mm}\\
\lefteqn{\bra{\La_0}\Ps{\an}{\an+1}\ket{\La_1}(u^{(\an)}_{0}) =
u^{(\an+1)}_{0},}\label{eq:460}\hspace{23mm}\\
\bra{\la_+}\Phi^+\ket{\la} = u^{(1)}_1, & &
\bra{\la_+}\Psi^{*+}\ket{\la}(u^{(1)}_0) = 1, \\
\bra{\la_-}\Phi^-\ket{\la} = u^{(1)}_0, & &
\bra{\la_-}\Psi^{*-}\ket{\la}(u^{(1)}_1) = 1.
\ena
Normalisation for the arbitrary level operators $\Ps{\an}{\am}$ will
be given later (page~\pageref{normal}, above Proposition~\ref{pr:reduc}).

The matrix elements of these intertwiners are Laurent polynomials. Therefore,
we can write the vertex operators in Laurent series expansions:
\begin{alignat}{2}
\Ps{\an}{\an+k}_{i,j}(\zeta)
  &= \sum_{\ttn} \zeta^{-\ttn}\Ps{\an}{\an+k}_{i,j,\ttn},&
\quad\Psz{\an}{\an+k}(u^{(\an)}_i\ot v)
  &= \sum_j \Ps{\an}{\an+k}_{i,j}(\zeta)v \ot u^{(\an+k)}_j,\\
\Phi^\pm_i(\zeta)
  &= \sum_{\ttn} \zeta^{-\ttn}\Phi^\pm_{i,\ttn},\label{4.12}&
\quad\Phi^\pm(\zeta)(v)
  &= \sum_{i=0,1} \Phi^\pm_i(\zeta)v\ot u^{(1)}_i,\\
\Psi^{*\pm}_i(\zeta)
  &= \sum_{\ttn} \zeta^{-\ttn}\Psi^{*\pm}_{i,\ttn},\label{4.13}&
\quad\Psi^{*\pm}(\zeta)(u^{(1)}_i\ot v)
  &= \Psi^{*\pm}_i(\zeta)v.
\end{alignat}
 From the normalisations, we see that for $\Phi^\ve_i$,
the sum is taken over $\ttn$ satisfying $\ve\cdot(-1)^{i+1} = (-1)^{\ttn}$
and that for $\Psi^{*\ve}_i$, it is taken over $\ttn$ satisfying
$\ve\cdot(-1)^{i} = (-1)^{\ttn}$.

Let us state some properties of the level-$1$ intertwiners.
Here, we suppress the appearance of the
$\pm$ superscripts on $\Phi^{\pm}(\z)$ and $\Psi^{\pm}(\xi)$.
Except where we state otherwise, these relations are valid for $\an\geq
0$, with the identification $\Phi^{(0,1)}(\z)=\Phi(\z)$.
\begin{prop}
{ \allowdisplaybreaks
\bea
\xi^{-D}\;\Ps{\an}{\an+1}_{s,t}(\zeta)\;\xi^{D}
    & = & \Ps{\an}{\an+1}_{s,t}(\zeta/\xi),\\
\de_{i+j, \an}
    & = & g^{(\an,\an+1)}
    \textstyle{\sum_{s+t=\an+1}} \Phi^{(\an,\an+1)}_{i,s}(-q^{-1}\zeta)
    \Phi^{(\an,\an+1)}_{j,t}(\zeta),\\
\Ps{\an}{\an+1}(\zeta)\Phi^{}(\xi)
    & = & R^{(1,\an+1)}(\xi/\zeta)\Phi(\xi)\Ps{\an}{\an+1}(\zeta),
\label{eq:phipm}\\
\Ps{\an}{\an+1}(\zeta)\Psi^{*}(\xi)
    & = & \Psi^{*}(\xi)\Ps{\an}{\an+1}(\zeta) R^{(\an,1)}(\zeta/\xi)
\quad \hb{for   }\an>0,\label{eq:psipm}\\
\Phi(\zeta)\Psi^{*}(\xi)
    & = & \Psi^{*}(\xi)\Phi(\zeta) \tau(\zeta/\xi),
\ena
\vspace*{-10mm}
\be
&&\!\!\!\!\!\!\!\!\!\!\!\!\!\!\!\!\!\!\!\!
\Ps{\an}{\an+1}(\zeta_1)\Ps{\am}{\am+1}(\zeta_2)
R^{(\am,\an)}(\zeta_2/\zeta_1)\nn\\
&&= \ \; R^{(\am+1,\an+1)}(\zeta_2/\zeta_1)
  \Ps{\am}{\am+1}(\zeta_2)\Ps{\an}{\an+1}(\zeta_1)\quad
\hb{for }\an,\an'>0,\label{eq:commrel}
\ena
} 
where
\bea
g^{(\an,\am)} = \frac{\qpf{q^{2\an+2}}}{\qpf{q^{2\am+2}}},\quad
\tau(\z)=\z^{-1}\frac{  \qpf{q\z^2}\qpf{q^3\z^{-2}} }%
    {\qpf{q\z^{-2}}\qpf{q^3\z^2}}
.\label{fix1}
\ena
\end{prop}
\begin{proof}
All but~\eqref{eq:commrel} appear in~\cite{MW97}.
The last one may be proved by applying~\eqref{eq:phipm} and~\eqref{eq:psipm}
on the fusion construction of $\Psz{\an}{\an+1}$ appearing in~\cite{Nak,MW97}.
\end{proof}

We wish to calculate various level-1, highest weight to highest weight,
matrix elements.

\begin{lem}\label{lm:42}\hfill
{\allowdisplaybreaks
\bea
&&\begin{aligned}
\bra{\La_0}\Ps{\an}{\an+1}(\zeta_1)\Ps{\am}{\am+1}(\zeta_2)\ket{\La_0}
  (u^{(\an)}_k \ot u^{(\am)}_j)\hspace{30mm}&\\
= h^{(\am+\an+5)}(\zeta_2/\zeta_1)
  \Big(\frac{1}{[\an+1]_q [\am+1]_q}\Big)^{\frac{1}{2}}\hspace{18mm}&\\
\times\Big\{
  q^{-\frac{1}{2}(\am+k-j)}[\an-k+1]_q^{\frac{1}{2}}[j+1]_q^{\frac{1}{2}}
  \;&u^{(\an+1)}_k\ot u^{(\am+1)}_{j+1}\\
- q^{\frac{1}{2}(\an-k+j)+1}[\am-j+1]_q^{\frac{1}{2}}[k+1]_q^{\frac{1}{2}}
  (\zeta_2/\zeta_1)\;&u^{(\an+1)}_{k+1} \ot u^{(\am+1)}_j \Big\},
\end{aligned}\\
&&\begin{aligned}
\bra{\La_1}\Ps{\an}{\an+1}(\zeta_1)\Ps{\am}{\am+1}(\zeta_2)\ket{\La_1}
  (u^{(\an)}_k \ot u^{(\am)}_j)\hspace{30mm}&\\
= h^{(\am+\an+5)}(\zeta_2/\zeta_1)
  \Big(\frac{1}{[\an+1]_q [\am+1]_q}\Big)^{\frac{1}{2}}\hspace{18mm}&\\
\times\Big\{
- q^{\frac{1}{2}(\am+k-j)+1}[\an-k+1]_q^{\frac{1}{2}}[j+1]_q^{\frac{1}{2}}
  (\zeta_2/\zeta_1)\;&u^{(\an+1)}_k\ot u^{(\am+1)}_{j+1}\\
+ q^{-\frac{1}{2}(\an-k+j)}[\am-j+1]_q^{\frac{1}{2}}[k+1]_q^{\frac{1}{2}}
  \;&u^{(\an+1)}_{k+1} \ot u^{(\am+1)}_j \Big\},
\end{aligned}
\ena
} 
where
$h^{(\an)}(\zeta) = \frac{\qpf{q^{\an+1}\zeta^2}}{\qpf{q^{\an-1}\zeta^2}}$.
\end{lem}
\begin{proof}
Let us just show the first one.
The second one follows from the first by applying the symmetry $\sigma$.
We set
$F = \bra{\La_0}\Ps{\an}{\an+1}(\zeta_1)\Ps{\am}{\am+1}(\zeta_2)\ket{\La_0}$
to simplify notations.
This is a map from $V^{(\an)}_{\zeta_1}\ot V^{(\am)}_{\zeta_2}$
to $V^{(\an+1)}_{\zeta_1}\ot V^{(\am+1)}_{\zeta_2}$.
By using $\bra{\La_0}e_1 = 0$, $f_0^2 \ket{\La_0} = 0$, and the fact that
$\Ps{\an}{\an+1}(\zeta_1)\Ps{\am}{\am+1}(\zeta_2)$ is an intertwiner,
we can show that
\begin{align}
F(e_1 u) - e_1 F(u) &= 0,\label{eq:470}\\
f_0^2 F(u) -[2]_q f_0 F(f_0 u) + F(f_0^2 u) &= 0\label{eq:471}
\end{align}
for any $u\in V^{(\an)}_{\zeta_1}\ot V^{(\am)}_{\zeta_2}$.
 From weight considerations, we know $F(u^{(\an)}_0\ot u^{(\am)}_0)$
is a linear combination of $u^{(\an+1)}_0\ot u^{(\am+1)}_1$ and
$u^{(\an+1)}_1\ot u^{(\am+1)}_0$.
The use of \eqref{eq:470} with $u = u^{(\an)}_0\ot u^{(\am)}_0$ allows us
to write
\begin{align*}
F(u^{(\an)}_0\ot u^{(\am)}_0) = 
f_{\an,\am}(\zeta_2/\zeta_1) \Big\{\phantom{-}
&q^{-\half\am}[\an+1]_q^{\half} \;u^{(\an+1)}_0\ot u^{(\am+1)}_1\\
- &q^{\half\an+1}[\am+1]_q^{\half} (\zeta_2/\zeta_1)
\;u^{(\an+1)}_1 \ot u^{(\am+1)}_0 \Big\}.
\end{align*}
Starting from this, and by
using \eqref{eq:470} and \eqref{eq:471}, we may inductively determine
$F$ up to the constant multiple $f_{\an,\am}(\zeta_2/\zeta_1)$.
To obtain the constant, we use \eqref{eq:commrel}, which shows
\begin{equation*}
\frac{f_{\an,\am}(\zeta_2/\zeta_1)}{f_{\am,\an}(\zeta_1/\zeta_2)} =
\frac{\qpf{q^{\an+\am+6}(\zeta_2/\zeta_1)^2}}
     {\qpf{q^{\an+\am+4}(\zeta_2/\zeta_1)^2}}
\frac{\qpf{q^{\an+\am+4}(\zeta_1/\zeta_2)^2}}
     {\qpf{q^{\an+\am+6}(\zeta_1/\zeta_2)^2}},
\end{equation*}
together with the normalisation \eqref{eq:459} and \eqref{eq:460}.
\end{proof}

\begin{lem}\hfill
{\allowdisplaybreaks
\bea
&&\begin{aligned}
\bra{\La_0}\Psz{\an}{\an+1}\Phi^+(\xi)&\ket{\La_0}(u^{(\an)}_j)\\
= h^{(\an+5)}(\xi/\zeta)&\Big(\frac{1}{[\an+1]_q}\Big)^{\frac{1}{2}}
   \Big\{ \phantom{ii}
    q^{-\frac{1}{2}j}[\an-j+1]_q^{\frac{1}{2}}
    \;u^{(\an+1)}_j \ot u^{(1)}_1\label{eq:408}\\
&\phantom{iiiiiiii}
  - q^{\frac{1}{2}(\an-j)+1}[j+1]_q^{\frac{1}{2}}(\xi/\zeta)
    \;u^{(\an+1)}_{j+1}\ot u^{(1)}_0
   \Big\},
\end{aligned}\\
&&\begin{aligned}
\bra{\La_1}\Psz{\an}{\an+1}\Phi^-(\xi)&\ket{\La_1}(u^{(\an)}_j)\\
= h^{(\an+5)}(\xi/\zeta)&\Big(\frac{1}{[\an+1]_q}\Big)^{\frac{1}{2}}
   \Big\{ \text{}
  - q^{\frac{1}{2}j+1}[\an-j+1]^{\frac{1}{2}}_q(\xi/\zeta)
    \;u^{(\an+1)}_j\ot u^{(1)}_1\label{eq:409}\\
&\phantom{iiiiiiiiiiiiiiiiiiiiiiiii}
  + q^{-\frac{1}{2}(\an-j)}[j+1]^{\frac{1}{2}}_q
    \;u^{(\an+1)}_{j+1}\ot u^{(1)}_0
   \Big\},
\end{aligned}\\
&&\begin{aligned}
\bra{\La_0}\Phi^-(\xi)\Psz{\an}{\an+1}&\ket{\La_0}(u^{(\an)}_j)\\
= h^{(\an+5)}(\zeta/\xi)&\Big(\frac{1}{[\an+1]_q}\Big)^{\frac{1}{2}}
   \Big\{ \phantom{ii}
    q^{-\frac{1}{2}(\an-j)}[j+1]_q^{\frac{1}{2}}
    \;u^{(1)}_0 \ot u^{(\an+1)}_{j+1}\\
&\phantom{iiiiii}
  - q^{\frac{1}{2}j+1}[\an-j+1]_q^{\frac{1}{2}}(\zeta/\xi)
    \;u^{(1)}_1 \ot u^{(\an+1)}_j
   \Big\},
\end{aligned}\\
&&\begin{aligned}
\bra{\La_1}\Phi^+(\xi)\Psz{\an}{\an+1}&\ket{\La_1}(u^{(\an)}_j)\\
= h^{(\an+5)}(\zeta/\xi)&\Big(\frac{1}{[\an+1]_q}\Big)^{\frac{1}{2}}
   \Big\{ \text{}
  - q^{\frac{1}{2}(\an-j)+1}[j+1]_q^{\frac{1}{2}}(\zeta/\xi)
    \;u^{(\an+1)}_0 \ot u^{(1)}_{j+1}\\
&\phantom{iiiiiiiiiiiiiiiiiiiiiiiiii}
  + q^{-\frac{1}{2}j}[\an-j+1]_q^{\frac{1}{2}}
    \;u^{(1)}_1 \ot u^{(\an+1)}_j
   \Big\}.
\end{aligned}
\ena
} 
\end{lem}
\begin{proof}
Just set $\an$ or $\am$ to zero in the preceding Lemma.
\end{proof}

Arguments similar to the proof of Lemma~\ref{lm:42} show:
\begin{lem}\hfill
{\allowdisplaybreaks
\bea
\bra{\La_0} \Psz{\an}{\an+1}\Psi^{*+}(\xi) \ket{\La_0}
      (u^{(\an)}_j \ot u^{(1)}_0)
   &\!\!=\!\!& h^{(\an+2)}(\xi/\zeta)
      \big(\textstyle{\frac{[\an+1-j]_q}{[\an+1]_q}}\big)^{\frac{1}{2}}
      q^{-\frac{1}{2}j}
      u^{(\an+1)}_j,\label{eq:412}\\
\bra{\La_0} \Psz{\an}{\an+1}\Psi^{*+}(\xi) \ket{\La_0}
      (u^{(\an)}_j \ot u^{(1)}_1)
   &\!\!=\!\!& h^{(\an+2)}(\xi/\zeta)
      \big(\textstyle{\frac{[j+1]_q}{[\an+1]_q}}\big)^{\frac{1}{2}}
      q^{\frac{1}{2}(\an-j)} (\xi/\zeta)
      u^{(\an+1)}_{j+1},\label{eq:413}\\
\bra{\La_1} \Psz{\an}{\an+1}\Psi^{*-}(\xi) \ket{\La_1}
      (u^{(\an)}_j \ot u^{(1)}_0)
   &\!\!=\!\!& h^{(\an+2)}(\xi/\zeta)
      \big(\textstyle{\frac{[\an+1-j]_q}{[\an+1]_q}}\big)^{\frac{1}{2}}
      q^{\frac{1}{2}j} (\xi/\zeta)
      u^{(\an+1)}_j,\label{eq:414}\\
\bra{\La_1} \Psz{\an}{\an+1}\Psi^{*-}(\xi) \ket{\La_1}
      (u^{(\an)}_j \ot u^{(1)}_1)
   &\!\!=\!\!& h^{(\an+2)}(\xi/\zeta)
      \big(\textstyle{\frac{[j+1]_q}{[\an+1]_q}}\big)^{\frac{1}{2}}
      q^{-\frac{1}{2}(\an-j)}
      u^{(\an+1)}_{j+1},\label{eq:415}\\
\bra{\La_0} \Psi^{*-}(\xi)\Psz{\an}{\an+1} \ket{\La_0}
      (u^{(1)}_0 \ot u^{(\an)}_j)
   &\!\!=\!\!& h^{(\an+2)}(\zeta/\xi)
      \big(\textstyle{\frac{[\an+1-j]_q}{[\an+1]_q}}\big)^{\frac{1}{2}}
      q^{\frac{1}{2}j} (\zeta/\xi)
      u^{(\an+1)}_j,\\
\bra{\La_0} \Psi^{*-}(\xi)\Psz{\an}{\an+1} \ket{\La_0}
      (u^{(1)}_1 \ot u^{(\an)}_j)
   &\!\!=\!\!& h^{(\an+2)}(\zeta/\xi)
      \big(\textstyle{\frac{[j+1]_q}{[\an+1]_q}}\big)^{\frac{1}{2}}
      q^{-\frac{1}{2}(\an-j)}
      u^{(\an+1)}_{j+1},\\
\bra{\La_1} \Psi^{*+}(\xi)\Psz{\an}{\an+1} \ket{\La_1}
      (u^{(1)}_0 \ot u^{(\an)}_j)
   &\!\!=\!\!& h^{(\an+2)}(\zeta/\xi)
      \big(\textstyle{\frac{[\an+1-j]_q}{[\an+1]_q}}\big)^{\frac{1}{2}}
      q^{-\frac{1}{2}j}
      u^{(\an+1)}_j,\\
\bra{\La_1} \Psi^{*+}(\xi)\Psz{\an}{\an+1} \ket{\La_1}
      (u^{(1)}_1 \ot u^{(\an)}_j)
   &\!\!=\!\!& h^{(\an+2)}(\zeta/\xi)
      \big(\textstyle{\frac{[j+1]_q}{[\an+1]_q}}\big)^{\frac{1}{2}}
      q^{\frac{1}{2}(\an-j)} (\zeta/\xi)
      u^{(\an+1)}_{j+1}.
\ena
} 
\end{lem}


\setcounter{equation}{0}
\section{Commutativity with the DVA Action}\label{sec5}

In this section, we show that $\Psz{\am}{\am+1}\Psz{\an}{\am}$ commutes with
the DVA action on $V(\la)\ot V(\La_i)$.
This will allow us to reduce questions about general level intertwiners to
those of level-$1$ intertwiners. Results on general level intertwiners will
be used in Section~\ref{sec63} to diagonalise the transfer matrix.

Let $\la\in P^+_k$ and define $\psi^{(\la,i,\pm)}(\xi)$ to be the composition
of operators given by:
\bea
V(\la)\ot V(\La_i)
\xrightarrow{\Phi^\pm(\xi)\ot\id}
V(\la_\pm)\ot V^{(1)}_\xi \ot V(\La_i)
\xrightarrow{\id\ot\Psi^*(\xi)}
V(\la_\pm)\ot V(\La_{1-i}).
\ena
This is equivalent to defining
\begin{equation}
\psi^{(\la,i,\pm)}(\xi)
 = \sum_{\ttn}\psi^{(\la,i,\pm)}_{\ttn} \xi^{-\ttn}
 = \sum_{j=0,1}\Phi^\pm_j(\xi)\ot\Psi^*_j(\xi).
\end{equation}
Each component
\begin{equation}
\psi^{(\la,i,\pm)}_{\ttn} :
 V(\la)\ot V(\La_i)\rightarrow V(\la_\pm)\ot V(\La_{1-i})
\end{equation}
is a $\uqp$ homomorphism.
In~\cite{JS97}, Jimbo and Shiraishi considered the irreducible decomposition
\bea
V(\la)\ot V(\La_i) = \bigoplus_\nu V(\nu)\ot \Omega^{\la,\La_i}_{\ \nu},
\label{eq:tdcmp}
\ena
and constructed an action of the deformed Virasoro algebra on
$\Omega^{\la,\La_i}_\nu$
by using the operator $\psi^{(\la,i,\pm)}(\xi)$.

Now, define $\phi^{(\an,\la,i)}(\zeta)$ to be the
composition of operators given by:
\bea
\lefteqn{ V^{(\an)}_\zeta\ot V(\la)\ot V(\La_i)
          \xrightarrow{\Psz{\an}{\an+k}\ot\id}
          V(\s(\la))\ot V^{(\an+k)}_\zeta \ot V(\La_i)}
          \hspace{40mm}\nn\\
& & \xrightarrow{\id\ot\Psz{\an+k}{\an+k+1}}
V(\s(\la))\ot V(\La_{1-i})\ot V^{(\an+k+1)}_\zeta.\nn
\ena
We shall show
\bea
\phi^{(\an,\la_\pm,1-i)}(\zeta)\circ\big(\id\ot\psi^{(\la,i,\pm)}(\xi)\big) =
\big(\id\ot\psi^{(\s(\la),1-i,\mp)}(\xi)\big)\circ\phi^{(\an,\la,i)}(\zeta).
\label{eq:DVAcomm}
\ena
This will imply the commutativity of $\phi^{(\an,\la,i)}$ with the
DVA action mentioned above.

Let us state a small lemma before considering the level-$1$ case.
\begin{lem}\label{lm:51}
Fix any $\uqp$ modules $V$ and $W$.
Let $\Theta : V(\la)\ot V \rightarrow W\ot V(\mu)$ be a $\uqp$ intertwiner.
Then any matrix element of $\Theta$ \textup{(}as an operator in
$\operatorname{End}(V,W)$\textup{)} may be written in the form
\bea
\sum \;\big(\makebox[11mm]{$\cdot$ $\cdot$ $\cdot$}\big)
       \circ\bra{\mu}\Theta\ket{\la}
       \circ\big(\makebox[11mm]{$\cdot$ $\cdot$ $\cdot$}\big)\; ,
\ena
where the parentheses signify appropriate $\uqp$ actions on $V$ and $W$,
respectively.
\end{lem}
\begin{proof}
The proof follows from the simple fact that $\Theta$ is an intertwiner.
\end{proof}

Now we show \eqref{eq:DVAcomm} in the level-$1$ case.

\begin{prop}\label{pr:51}
The equality,
\bea
\phi^{(\an,\La_{1-j},1-i)}(\zeta)\circ
\big(\id\ot\psi^{(\La_j,i,\pm)}(\xi)\big)
=
\big(\id\ot\psi^{(\La_{1-j},1-i,\mp)}(\xi)\big)\circ
\phi^{(\an,\La_j,i)}(\zeta)
\ena
holds for their matrix elements as Laurent series of $\zeta$ and $\xi$.
They contain no poles.
\end{prop}
\begin{proof}
With the help of equation~\eqref{eq:phipm} and~\eqref{eq:psipm}
we can show:
{\allowdisplaybreaks
\bea
\lefteqn{
    \big(\Psz{\an+1}{\an+2}\Psz{\an}{\an+1}\big)\circ
    \big(\Psi^*(\xi)\Phi(\xi)\big)
}\hspace{25mm}\nn\\
&=& \Psz{\an+1}{\an+2}\Psi^*(\xi)\Psz{\an}{\an+1}\Phi(\xi)\nn\\
&=& \Psz{\an+1}{\an+2}\Psi^*(\xi)
    R^{(1,\an+1)}(\xi/\zeta)\Phi(\xi)\Ps{\an}{\an+1}(\zeta)\nn\\
&=& \Psi^*(\xi)\Psz{\an+1}{\an+2}
    \Phi(\xi)\Ps{\an}{\an+1}(\zeta)\nn\\
&=& \big(\Psi^*(\xi)\Phi(\xi)\big)\circ
    \big(\Psz{\an+1}{\an+2}\Psz{\an}{\an+1}\big).\nn
\ena
} 
So the two sides agree as meromorphic functions.
Let us look at the structure of poles.
We have
\bea
\lefteqn{
\bra{\La_j}\ot\bra{\La_i}
\big(\Psz{\an+1}{\an+2}\Psz{\an}{\an+1}\big)\circ
\big(\Psi^*(\xi)\Phi(\xi)\big)
\ket{\La_j}\ot\ket{\La_i}}
\hspace{30mm}\nn\\
&=& \bra{\La_i}
   \Psz{\an+1}{\an+2}\Psi^*(\xi)
   \ket{\La_i}
\circ
   \bra{\La_j}
   \Psz{\an}{\an+1}\Phi(\xi)
   \ket{\La_j}.\nn
\ena
Using the equations~\eqref{eq:408}, \eqref{eq:409} and also equations
\eqref{eq:412} to \eqref{eq:415} with $\an$ replaced by $\an+1$, we see that
a pole can occur in the above only if $1-q^{\an+2}(\xi/\zeta)^2 = 0$.
If $\zeta$ and $\xi$ satisfy this relation, 
there exists a submodule isomorphic to some $V^{(\an)}_\mu$
lying inside $V^{(\an+1)}_\zeta \ot V^{(1)}_\xi$.
When $\xi = \zeta q^{-\half(\an+2)}$,
a submodule isomorphic to $V^{(\an)}_{\zeta q^{\half}}$ lying inside
$V^{(\an+1)}_\zeta\ot V^{(1)}_\xi$ is spanned by
\bea
\Big\{
  [\an+1-k]^{\half}_q \; u^{(\an+1)}_{k} \ot u^{(1)}_1
- [k+1]^{\half}_q \; u^{(\an+1)}_{k+1} \ot u^{(1)}_0
\Big\}_{k=0}^{\an}.
\ena
Again, from the same set of equations, we see that the image of
$\bra{\La_j} \Psz{\an}{\an+1}\Phi(\xi) \ket{\La_j}$
lies inside this submodule.
We also see that $\bra{\La_i} \Psz{\an+1}{\an+2}\Psi^*(\xi) \ket{\La_i}$
sends this submodule to zero.
Therefore, the above matrix element contains no pole.

In view of Lemma~\ref{lm:51}, the general matrix element can be written
in the form
\bea
\sum\;
   \big(\makebox[8mm]{$\cdot$ $\cdot$ $\cdot$}\big)
\circ
   \bra{\La_i}
   \Psz{\an+1}{\an+2}\Psi^*(\xi)
   \ket{\La_i}
\circ
   \big(\makebox[8mm]{$\cdot$ $\cdot$ $\cdot$}\big)
\circ
   \bra{\La_j}
   \Psz{\an}{\an+1}\Phi(\xi)
   \ket{\La_j}
\circ
   \big(\makebox[8mm]{$\cdot$ $\cdot$ $\cdot$}\big)\; ,
\ena
where the parentheses are to be filled with $\uqp$ actions.
As the action of $\uqp$ cannot produce additional poles,
the only possible poles will occur at $1-q^{\an+2}(\xi/\zeta)^2 = 0$.
Again, the image of the first map will lie inside some
submodule $V^{(\an)}_\mu$.
The $\uqp$-action will still preserve this submodule.
Then the second map will send this submodule to zero.
This shows that the general matrix element also contains no pole.
\end{proof}

We now show that the commutativity with the DVA action allows us to construct
$\Ps{\an}{\an+k}$ from lower level operators.
Assume from now on that equation~\eqref{eq:DVAcomm}
is true as Laurent series for $\la\in P^+_{k-1}$.
Now, $\Psz{\an+k-1}{\an+k}\Psz{\an}{\an+k-1}$ is a map from
\bea
V^{(\an)}_\zeta \ot V(\la) \ot V(\La_i)
= \bigoplus_\nu V^{(\an)}_\zeta \ot V(\nu) \ot \Omega^{\la,\La_i}_{\ \nu},
\ena
where the sum runs over all level-$k$ weights, to the space
\bea
V(\s(\la))\ot V(\La_{1-i}) \ot V^{(\an+k)}_\zeta
= \bigoplus_\nu V(\s(\nu))\ot V^{(\an+k)}_\zeta \ot \Omega^{\la,\La_i}_{\ \nu}.
\ena
Here, we have used the ${\bf Z}_2$-symmetry to write
$\Omega^{\s(\la),\La_{1-i}}_{\ \s(\nu)} = \Omega^{\la,\La_i}_{\ \nu}$.
Recalling Proposition~\ref{pr:ex}, we may write
\bea
\label{eq:524}
\Psz{\an+k-1}{\an+k}\Psz{\an}{\an+k-1}
= \bigoplus_\nu \Ps{\an}{\an+k}_\nu(\zeta)\ot\Xi^{\la,\La_i}_{\ \nu}.
\ena
The subscript $\nu$ in $\Ps{\an}{\an+k}_\nu(\zeta)$ has been added
to show which space it acts on.
Now, each $\Omega^{\la,\La_i}_{\ \nu}$ is irreducible.
Hence the commutativity with the DVA action shows that each
$\Xi^{\la,\La_i}_{\ \nu}$ is a constant map.
We normalise the higher level intertwiner $\Ps{\an}{\an+k}_\nu(\zeta)$
so that this constant is equal to $1$ for the highest\label{normal}
component, i.e., for $\nu=\la + \La_i$.
This normalisation is independent of the way we break up
the level-$k$ weight into level-$1$ weights, as can be seen by the use of
\bea
\bra{\La_1}\Psz{\an}{\an+1}\ket{\La_0}(u^{(\an)}_j) & = &
\Big(q^{j-\ell}\frac{[j+1]}{[\ell+1]}\Big)^{\half} u^{(\an+1)}_{j+1},\\
\bra{\La_0}\Psz{\an}{\an+1}\ket{\La_1}(u^{(\an)}_j) & = &
\Big(q^{-j}\frac{[\ell-j+1]}{[\ell+1]}\Big)^{\half} u^{(\an+1)}_j.
\ena
The next proposition is more of a definition.

\begin{prop}\label{pr:reduc}
The map $\Psz{\an+k-1}{\an+k}\Psz{\an}{\an+k-1}$ restricted to
$V^{(\an)}_\zeta\ot V(\la+\La_i) \ot \Omega^{\la,\La_i}_{\la+\La_i}$
is equal to $\Ps{\an}{\an+k}_{\la+\La_i}(\zeta)\ot \id$.
\end{prop}

Later on, we will say something about the coefficients of other components.
But for now, let us continue with proving the commutativity with the
DVA action.
Since we now know how to construct $\Psz{\an}{\an+k}$ from lower level
operators, we can find its commutation relations.
Except where we state otherwise, the relations we give in the following
proposition are valid for $\an\geq 0$.
\begin{prop}\label{pr:5.4}\hfill
{\allowdisplaybreaks
\bea
\xi^{-D}\;\Ps{\an}{\am}_{s,t}(\zeta)\;\xi^{D} &=&
\Ps{\an}{\am}_{s,t}(\zeta/\xi),\label{5.129}\\
\delta_{i+j, \an}
    &=& g^{(\an,\am)} \textstyle{\sum_{s+t=\am}}
        \Phi^{(\an,\am)}_{i,s}(-q^{-1}\zeta)
        \Phi^{(\an,\am)}_{j,t}(\zeta),\label{5.130}\\
\Psz{\an}{\am}\Phi^{\pm}(\xi)
    &=& R^{(1,\am)}(\xi/\zeta)\Phi^{\mp}(\xi)\Psz{\an}{\am},\label{eq:527}\\
\Psz{\an}{\am}\Psi^{*\pm}(\xi)
    &=& \Psi^{*\mp}(\xi)\Psz{\an}{\am} R^{(\an,1)}(\zeta/\xi)\quad
\hb{for } \an>0,\label{eq:psm}\\
\Psz{0}{\an}\Psi^{*\pm}(\xi)
    &=& \Psi^{*\mp}(\xi)\Psz{0}{\an} \tau(\zeta/\xi).
\label{eq:psmzero}
\ena
} 
\end{prop}
\begin{proof}
Let us prove equation~\eqref{eq:527} as an example. We will consider
just one set of the signs involved.
First consider the map
\bea
\id\ot\Phi^{+}(\xi) &:&
V(\la)\ot V(\La_0)
\rightarrow
V(\la)\ot V(\La_{1})\ot V^{(1)}_\xi.
\ena
As before, we may write this map as
\bea
\id\ot\Phi^{+}(\xi) &=&
\big( \bigoplus_\nu \Phi^{+}_\nu(\xi) \ot \Xi^+_\nu \big) +
\big( \bigoplus_\nu \Phi^{-}_\nu(\xi) \ot \Xi^-_\nu \big),
\label{eq:528}
\ena
where the maps on the right hand side are from
$V(\nu)\ot\Omega^{\la,\La_0}_{\ \nu}$ to
$V(\nu_\pm)\ot V^{(1)}_\xi \ot \Omega^{\la,\La_1}_{\ \nu_\pm}$.
As in the proof of Proposition~\ref{pr:reduc} we take the highest weight
matrix element and apply both sides on $u^{(\an)}_0$ to conclude
$\Xi^+_{\la+\La_0}$ is nonzero.
By ${\bf Z}_2$-symmetry, the map
\bea
\id\ot\Phi^{-}(\xi) &:&
V(\s(\la))\ot V(\La_{1})
\rightarrow
V(\s(\la))\ot V(\La_0)\ot V^{(1)}_\xi
\ena
breaks up as
\bea
\id\ot\Phi^-(\xi) &=&
\big( \bigoplus_\nu \Phi^-_{\s(\nu)}(\xi)\ot \Xi^+_\nu \big) +
\big( \bigoplus_\nu \Phi^+_{\s(\nu)}(\xi)\ot \Xi^-_\nu \big).
\label{eq:529}
\ena
With this much in hand, we proceed by induction.
By equation \eqref{eq:phipm} we have,
\bea
\lefteqn{\Psz{\an+k-1}{\an+k}\,
\Psz{\an}{\an+k-1}(\id\ot\Phi^+(\xi))}\hspace{28mm}\\
&=& \Psz{\an+k-1}{\an+k}\big(\id\ot\Phi^+(\xi)\big)\,\Psz{\an}{\an+k-1}\\
&=& R^{(1,\an+k)}(\xi/\zeta)\big(\id\ot\Phi^-(\xi)\big)
    \,\Psz{\an+k-1}{\an+k}\,\Psz{\an}{\an+k-1}.
\ena
If we substitute equations~\eqref{eq:524}, \eqref{eq:528}, and \eqref{eq:529}
into both sides and pick up the term initiating at
$V(\la+\La_0)\ot\Omega^{\la,\La_0}_{\la+\La_0}$ and terminating at
$V(\s(\la+\La_1))\ot\Omega^{\la, \La_1}_{\la+\La_1}$,
and apply Proposition~\ref{pr:reduc} to the outcome, we get:
\bea
\big(\Psz{\an}{\an+k}\Phi^+_{\la+\La_0}(\xi)\big)\ot\Xi^+_{\la+\La_0} =
\big(R^{(1,\an+k)}(\xi/\zeta)\Phi^-_{\s(\la+\La_0)}(\xi)
\Psz{\an}{\an+k}\big)\ot\Xi^+_{\la+\La_0}.
\ena
We already know $\Xi^+_{\la+\La_0}$ is nonzero, so dividing them out,
we have the result.
\end{proof}

In much the same way, we can also calculate the higher level matrix elements.
Here we only write down what is needed in proving the commutativity with the
DVA action.
\begin{lem}\label{lm}
Let $\la$ be of level $k-1$. Then
\bea
\lefteqn{
\bra{\s(\la)+\La_0} \Psz{\an}{\an+k}\Phi^+(\xi) \ket{\la+\La_0}}
\hspace{25mm}\nn\\
  &=& \bra{\La_0} \Psz{\an+k-1}{\an+k}\Phi^+ \ket{\La_0} \circ
     \bra{\s(\la)} \Psz{\an}{\an+k-1} \ket{\la},\nn\\
\lefteqn{
\bra{\s(\la)+\La_1} \Psz{\an}{\an+k}\Phi^-(\xi) \ket{\la+\La_1}}
\hspace{25mm}\nn\\
  &=& \bra{\La_1} \Psz{\an+k-1}{\an+k}\Phi^- \ket{\La_1} \circ
     \bra{\s(\la)} \Psz{\an}{\an+k-1} \ket{\la}.\nn
\ena
\end{lem}

We are now ready for the induction step in proving commutativity with the
DVA action.
\begin{thm}
The equality,
\bea
\phi^{(\an,\la_\pm,1-i)}(\zeta)\circ\big(\id\ot\psi^{(\la,i,\pm)}(\xi)\big)
&=&
\big(\id\ot\psi^{(\s(\la),1-i,\mp)}(\xi)\big)\circ\phi^{(\an,\la,i)}(\zeta),
\ena
holds for their matrix elements as Laurent series of $\zeta$ and $\xi$.
They contain no poles.
\end{thm}
\begin{proof}
We are assuming that the statement is true for levels less than $k$.
Hence Proposition~\ref{pr:reduc}, Proposition~\ref{pr:5.4} and Lemma~\ref{lm}
hold true for level $k$.
Applying \eqref{eq:527} and then \eqref{eq:psipm}, we show the equality
of both sides as meromorphic functions.
\bea
\text{LHS}
&=& \Psz{\an+k}{\an+k+1}\Psi^*(\xi)\Psz{\an}{\an+k}\Phi(\xi)\nn\\
&=& \Psz{\an+k}{\an+k+1}\Psi^*(\xi) R^{(1,\an+k)}(\xi/\zeta)
    \Phi(\xi)\Psz{\an}{\an+k}\nn\\
&=& \Psi^*(\xi)\Psz{\an+k}{\an+k+1}\Phi(\xi)\Psz{\an}{\an+k}\nn\\
&=& \text{RHS}.\nn
\ena

For the rest of the proof, we will explain the case when
the left hand side contains the $+$ sign.
The other case can be taken care of similarly.
For this case, we may write $\la = \la'+\La_0$ with $\la'$ of level $k-1$.
We use Lemma~\ref{lm} to show:
\bea
\lefteqn{
\bra{\s(\la_+)}\ot\bra{\La_i}
\big(\Psz{\an+k}{\an+k+1}\Psz{\an}{\an+k}\big)\circ
\big(\Psi^*(\xi)\Phi(\xi)\big)
\ket{\la}\ot\ket{\La_i}}
\hspace{10mm}\nn\\
&=& \bra{\La_i}\Psz{\an+k}{\an+k+1}\Psi^*(\xi)\ket{\La_i}
\circ\bra{\s(\la')+\La_0}\Psz{\an}{\an+k}\Phi(\xi)\ket{\la'+\La_0}\nn\\
&=& \bra{\La_i}\Psz{\an+k}{\an+k+1}\Psi^*(\xi)\ket{\La_i}\nn\\
& &\phantom{iiiiiiiiiiiiiii}
\circ \bra{\La_0} \Psz{\an+k}{\an+k+1}\Phi^+(\xi) \ket{\La_0}
\circ \bra{\s(\la')} \Psz{\an}{\an+k} \ket{\la'}.\nn
\ena
Now, we may argue as in the proof of Proposition~\ref{pr:51} to show
that the above identity contains no pole.
\end{proof}

We remark on the coefficients of
the components of $\Psz{\an+k-1}{\an+k}\Psz{\an}{\an+k}$
before closing this section.
\begin{prop}
\bea
\label{eq:531}
\Psz{\an+k-1}{\an+k}\Psz{\an}{\an+k-1}
&=& \bigoplus_\nu c_\nu \cdot \Ps{\an}{\an+k}_\nu(\zeta)\ot\id,
\ena
with each $(c_\nu)^2 = 1$.
\end{prop}
\begin{proof}
We have only to show $(c_\nu)^2 = 1$.
Using equation \eqref{5.130}, we have
\bea
\lefteqn{
g^{(\an,\an+k)}
\bigoplus_\nu \big\{ \sum_{s+t=\an+k} (c_\nu)^2 \cdot
\Ps{\an}{\an+k}_{\nu,\ i,s}(-q^{-1}\zeta) \Ps{\an}{\an+k}_{\nu,\ j,t}(\zeta)
\ot \id \big\} }
\hspace{60mm}\\
&=& \bigoplus_\nu \big\{ (c_\nu)^2 \cdot \id_{V(\nu)} \ot \id \big\}
\;\delta_{i+j, \an}.
\ena
If we calculate the same thing with the left hand side expression
of equation~\eqref{eq:531}, we get
\bea
\bigoplus_\nu \big\{ \id_{V(\nu)} \ot \id \big\}\;\delta_{i+j,\an}.
\ena
This shows that each $(c_\nu)^2 = 1$.
\end{proof}

\setcounter{equation}{0}
\section{Diagonalisation of the Transfer Matrix}\label{sec6}
In this section, we identify the space of states, and the half
and full transfer matrices of the alternating spin vertex model
in terms of the representation
theory of $U_q(\widehat{sl}_2)$. 
We diagonalise the full transfer matrix in terms of the spin-$0$
and spin-$\half$ states mentioned in the introduction, and compute 
two-particle S-matrix elements.
\subsection{The Space of States}

In Section~\ref{sec3}, we have shown that there is a crystal 
isomorphism $P_{a,b}\simeq 
B(\la^{(m-n)}_a)\ot B(\la^{(n)}_b)$. This leads us to conjecture that 
we can extend this isomorphism away from $q=0$, and identify
the space of eigenstates of the corner transfer matrix
$A_{NW}(\z)$ with $\cH\equiv \bigoplus_{a,b} \cH_{a,b}$, where 
$\cH_{a,b}=V(\la^{(m-n)}_a)\ot V(\la^{(n)}_b)$, and $0\leq a \leq m-n$,
$0\leq b \leq n$.
The operator $A_{NW}(\z)$ will act as $\cH_{a,b}\ra\cH_{a,b}$.
Then $\cF\equiv\cH\ot \cH^{*}$ will be the
space on which our full transfer matrix acts. 
Here, $\cH^{*}$ is the dual space, defined using the $\uqp$
anti-automorphism $b$ given in~\cite{Miki}.
The motivation for this definition, and the reason for the use
of this particular anti-automorphism are discussed in the similar
context of the pure spin-$\half$ vertex model in~\cite{JM}.

We can identify an element $\bra{f}\in \cF^*$ with an element 
$\ket{f} \in \cF$ via the pairing $\bra{f}g\rangle=\Tr_{\cH}(f \circ g)$.
Here, we have used the canonical isomorphism $\cF\simeq \End(\cH)$ 
to identify $f,g\in \cF$ as elements of $\End(\cH)$ in the trace formula.

\subsection{Half and full transfer matrices}

A half transfer matrix represents the insertion of a half-infinite column
of lattice vertices. There are two types of half transfer matrices for
the alternating spin model - those associated with the insertions of
columns with
spin-$\frac{n}{2}$ and spin-$\frac{m}{2}$ vertical lines.
These are shown in Figures 3
(a) and (b) respectively. 

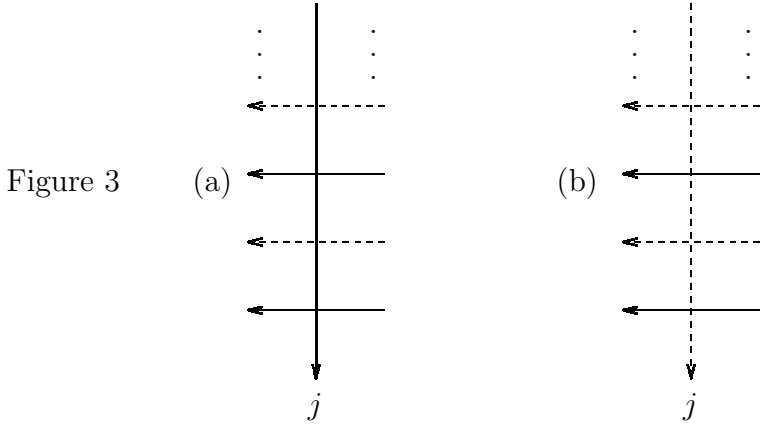
\setlength{\unitlength}{0.0006in}
\begingroup\makeatletter\ifx\SetFigFont\undefined%
\gdef\SetFigFont#1#2#3#4#5{%
  \reset@font\fontsize{#1}{#2pt}%
  \fontfamily{#3}\fontseries{#4}\fontshape{#5}%
  \selectfont}%
\fi\endgroup%
{\renewcommand{\dashlinestretch}{30}
\begin{picture}(5937,4500)(-800,-400)
\drawline(1425,861)(2625,861)
\drawline(1545.000,891.000)(1425.000,861.000)(1545.000,831.000)
\drawline(1425,2061)(2625,2061)
\drawline(1545.000,2091.000)(1425.000,2061.000)(1545.000,2031.000)
\dashline{60.000}(1425,1461)(2625,1461)
\drawline(1545.000,1491.000)(1425.000,1461.000)(1545.000,1431.000)
\dashline{60.000}(1425,2661)(2625,2661)
\drawline(1545.000,2691.000)(1425.000,2661.000)(1545.000,2631.000)
\drawline(4725,861)(5925,861)
\drawline(4845.000,891.000)(4725.000,861.000)(4845.000,831.000)
\drawline(4725,2061)(5925,2061)
\drawline(4845.000,2091.000)(4725.000,2061.000)(4845.000,2031.000)
\dashline{60.000}(4725,1461)(5925,1461)
\drawline(4845.000,1491.000)(4725.000,1461.000)(4845.000,1431.000)
\dashline{60.000}(4725,2661)(5925,2661)
\drawline(4845.000,2691.000)(4725.000,2661.000)(4845.000,2631.000)
\drawline(2025,3561)(2025,261)
\drawline(1995.000,381.000)(2025.000,261.000)(2055.000,381.000)
\dashline{60.000}(5325,3561)(5325,261)
\drawline(5295.000,381.000)(5325.000,261.000)(5355.000,381.000)
\put(1860,-50){\makebox(0,0)[lb]{\smash{{{\SetFigFont{12}{14.4}{\rmdefault}%
{\mddefault}{\updefault} $j$}}}}}
\put(5160,-50){\makebox(0,0)[lb]{\smash{{{\SetFigFont{12}{14.4}{\rmdefault}%
{\mddefault}{\updefault} $j$}}}}}
\put(4700,3300){\makebox(0,0)[lb]{\smash{{{\SetFigFont{12}{14.4}{\rmdefault}%
{\mddefault}{\updefault} .}}}}}
\put(4700,3100){\makebox(0,0)[lb]{\smash{{{\SetFigFont{12}{14.4}{\rmdefault}%
{\mddefault}{\updefault} .}}}}}
\put(4700,2900){\makebox(0,0)[lb]{\smash{{{\SetFigFont{12}{14.4}{\rmdefault}%
{\mddefault}{\updefault} .}}}}}
\put(5700,3300){\makebox(0,0)[lb]{\smash{{{\SetFigFont{12}{14.4}{\rmdefault}%
{\mddefault}{\updefault} .}}}}}
\put(5700,3100){\makebox(0,0)[lb]{\smash{{{\SetFigFont{12}{14.4}{\rmdefault}%
{\mddefault}{\updefault} .}}}}}
\put(5700,2900){\makebox(0,0)[lb]{\smash{{{\SetFigFont{12}{14.4}{\rmdefault}%
{\mddefault}{\updefault} .}}}}}
\put(1400,3300){\makebox(0,0)[lb]{\smash{{{\SetFigFont{12}{14.4}{\rmdefault}%
{\mddefault}{\updefault} .}}}}}
\put(1400,3100){\makebox(0,0)[lb]{\smash{{{\SetFigFont{12}{14.4}{\rmdefault}%
{\mddefault}{\updefault} .}}}}}
\put(1400,2900){\makebox(0,0)[lb]{\smash{{{\SetFigFont{12}{14.4}{\rmdefault}%
{\mddefault}{\updefault} .}}}}}
\put(2400,3300){\makebox(0,0)[lb]{\smash{{{\SetFigFont{12}{14.4}{\rmdefault}%
{\mddefault}{\updefault} .}}}}}
\put(2400,3100){\makebox(0,0)[lb]{\smash{{{\SetFigFont{12}{14.4}{\rmdefault}%
{\mddefault}{\updefault} .}}}}}
\put(2400,2900){\makebox(0,0)[lb]{\smash{{{\SetFigFont{12}{14.4}{\rmdefault}%
{\mddefault}{\updefault} .}}}}}
\put(-700,1911){\makebox(0,0)[lb]{\smash{{{\SetFigFont{12}{14.4}{\rmdefault}%
{\mddefault}{\updefault}Figure 3
\quad\quad    (a) }}}}}
\put(4050,1911){\makebox(0,0)[lb]{\smash{{{\SetFigFont{12}{14.4}{\rmdefault}%
{\mddefault}{\updefault}  (b)}}}}}
\end{picture}
}

As lattice insertions these will be the maps 
$\cH_{a,b}\ra \cH_{a,n-b}$ and $\cH_{a,b}\ra \cH_{m-n-a,n-b}$ respectively
(one can see this by an inspection of the ground state configuration
 shown in Figure 2).

As discussed in the introduction, we identify these lattice 
insertions with
components $\phi^A_j(\z)$ and $\phi^B_j(\z)$ of
the  following intertwiners.
\[\begin{array}{lllll}
\phi^A(\z)&:&V(\la_a^{(m-n)})\ot V(\la_b^{(n)}) &
\xrightarrow{\id\ot \Phi^{(0,n)}(\zeta)}&
V(\la_a^{(m-n)})\ot V(\s(\la_b^{(n)}))
         \ot V^{(n)}_{\z},\nn\\[3mm]
\phi^B(\z)&:&V(\la_a^{(m-n)})\ot V(\la_b^{(n)}) 
&
\xrightarrow{\Phi^{(0,m-n)}(\z)\ot\id}
& V(\s(\la_a^{(m-n)})) \ot V^{(m-n)}_{\z} \ot V(\la_b^{(n)})\nn \\
&& &
\xrightarrow{\id\ot \Phi^{(m-n,m)}(\z)}
& V(\s(\la_a^{(m-n)}))\ot V(\s(\la_b^{(n)}))
         \ot V^{(m)}_{\z}.
\end{array}\]
Here, $\Phi^{(k,l)}(\z)$ is the perfect intertwiner defined in
Section~\ref{sec4}.
If $v\ot v'\in V(\la_a^{(m-n)})\ot
V(\la_b^{(n)})$, then the components $\phi^A_j(\z)$ and
$\phi^B_j(\z)$ are defined by
\be\begin{array}{lllll}
\phi^A(v\ot v')&=& \sli_{j=0}^{n}\phi^A_j(v\ot v')\ot u^{(n)}_j,\quad
\phi^A_j(v\ot v')&=&v\ot \Phi^{(0,n)}_jv',\label{def1}\\[3mm]
\phi^B(v\ot v')&=& \sli_{j=0}^{m} \phi^B_j(v\ot v')\ot u^{(m)}_j,\quad
\phi^B_j(v\ot v')&=&\sli_{j'=0}^{m-n} \Phi^{(0,m-n)}_{j'}v\ot
\Phi^{(m-n,m)}_{j',j}v'.
\end{array}\ee
Here, for clarity, we have suppressed the $\z$ dependence of all 
our intertwiners.

Now consider the corresponding full transfer matrices, i.e., those 
associated with the insertion of full, double-infinite, columns 
of lattice vertices. Again there will be two such transfer matrices,
$T^A(\z)$ and $T^B(\z)$,
associated with spin-$\frac{n}{2}$ and spin-$\frac{m}{2}$ auxiliary spaces 
respectively. We identify these in terms of intertwiners that act on the tensor
product space $\cH_{a,b}\ot \cH_{a',b'}^*$ as follows
\begin{equation}
\cH_{a,b}\ot \cH_{a',b'}^* 
\xrightarrow{\phi^A(\z)\ot\id} \cH_{a,n-b}\ot V^{(n)}_{\z} \ot \cH_{a',b'}^*
\xrightarrow{\id\ot \phi^A(\z)^t} \cH_{a,n-b}\ot \cH_{a',n-b'}^*
\end{equation}
\begin{equation}
\begin{aligned}
\cH_{a,b}\ot \cH_{a',b'}^* 
\xrightarrow{\phi^B(\z)\ot\id} \cH_{m-n-a,n-b}&\ot V^{(m)}_{\z}
\ot \cH_{a',b'}^*\\
&\xrightarrow{\id\ot \phi^B(\z)^t} \cH_{m-n-a,n-b}\ot \cH_{m-n-a',n-b'}^*.
\end{aligned}
\end{equation}
Here $t$ denotes the transpose.
Specifically, we define
\be\begin{array}{lllll} 
     T^A(\z) &=& \sli_{j=0}^{n} T^A_j(\z),\quad  T^A_j(\z)&=& 
     g^{(0,n)} \, \phi^A_{j}(\z) \ot \phi^A_{n-j}(\z)^t,\label{Tdef1}\\
     T^B(\z) &=& \sli_{j=0}^{m}T^B_j(\z),\quad T^B_j(\z)&=&
     g^{(0,m)} \, \phi^B_{j}(\z) \ot \phi^B_{m-j}(\z)^t,\\
\end{array}\ee
where the constants $g^{(0,n)}$ and $g^{(0,m)}$ are given by
\eqref{fix1}.

\subsection{Diagonalisation of the full transfer matrices}\label{sec63}
A vacuum is, by definition, a largest eigenvalue eigenvector of the 
composition $\cT(\z)=T^B(\z) \circ T^A(\z)$. In our alternating spin model, 
there are $(m-n+1)(n+1)$ degenerate vacua 
$\vac_{a,b}\in \cH_{a,b}\ot \cH^*_{a,b}$. 
The expressions for these
vacua appear simple if we express them as elements of $\End(\cH)$.
We conjecture that the vacua $\vac_{a,b}$, and $_{a,b}\dvac$  are given
by
\begin{equation}
_{a,b}\langle{\rm vac}|=|{\rm vac}\rangle_{a,b}
=(\chi^{(m-n)}_{a} \chi^{(n)}_b)^{-\frac{1}{2}} (-q)^D
 \pi_{a,b},
\end{equation}
Here, $\chi^{(\ell)}_r$ is the character 
\be\chi^{(\ell)}_r=\Tr_{V(\la_r^{(\ell)})}(q^{2D}),\ee
whose appearance gives the normalisation $_{a,b}\langle \rm{vac} \vac_{a,b}=1$,
and $\pi_{a,b}\in{\rm End}({\cal H})$ is the projector to
${\cal H}_{a,b}$.

Let us consider the action of $\cT(\z)$ on $\vac_{a,b}$. First, note
that the action of
a map $\cO_1\ot \cO_2: \cH\ot \cH^* \ra \cH\ot \cH^*$ on an 
element $f\in \End(\cH)$  is given by $\cO_1 \circ f \circ \cO_2^t$.
Then, using \eqref{Tdef1}, \eqref{def1} and properties \eqref{5.129},
and \eqref{5.130} we have 
\bea
T^A(\z) \vac_{a,b} &=& \vac_{a,n-b},\label{TAvac}\\
T^B(\z) \vac_{a,b} &=& \vac_{m-n-a,n-b},\label{TBvac}\ena
and hence $\cT(\z)\vac_{a,b} = \vac_{m-n-a,b}$. 
To be precise about our use of the terminology
`eigenvector' or `eigenvalue',
the vacuum vector $|{\rm vac}\rangle_{a,b}$ is
not an eigenvector of ${\cal T}(\zeta)$
but of ${\cal T}(\zeta)^2$ or ${\cal T}(1)^{-1}{\cal T}(\zeta)$. However,
in the following we abuse this terminology, and call
$|{\rm vac}\rangle_{a,b}$
a vacuum eigenvector.

Let us show how 
to derive \eqref{TAvac}.
 From \eqref{Tdef1}, and \eqref{def1} we have,
\bea
T^A(\z) \vac_{a,b}
        &=& g^{(0,n)} \sli_{j=0}^n(\id\ot \Phi^{(0,n)}_{j}(\z))
       ((-q)^D \ot (-q)^D) \pi_{a,b} (\id\ot \Phi^{(0,n)}_{n-j}(\z)),\\
        &=& g^{(0,n)} \sli_{j=0}^n((-q)^D \ot
        \Phi^{(0,n)}_{j}(\z)(-q)^D\Phi^{(0,n)}_{n-j}(\z)) \pi_{a,n-b},\\
        &=& g^{(0,n)} \sli_{j=0}^n((-q)^D \ot
        (-q)^D)(\id\ot \Phi^{(0,n)}_{j}(-q^{-1}\z)\Phi^{(0,n)}_{n-j}(\z))
        \pi_{a,n-b},\label{step3}\\
 &=&   ((-q)^D\ot (-q)^D)\pi_{a,n-b} = \vac_{a,n-b}.\label{step4}
\ena
In executing steps \eqref{step3} and \eqref{step4}, 
we have used properties \eqref{5.129} and \eqref{5.130} respectively.
Equation \eqref{TBvac} can be shown similarly.

We will construct excited states by making use of the following
operators
\begin{equation}\label{IIdefs}
\begin{aligned}
\psi^{(\half)s}_\ve(\xi)
&= \sum_{\ttn} \psi^{(\half)s}_{\ve,\ttn} \xi^{-\ttn}
 = \Psi^{*s}_{\ve}(\xi)\ot\id,\\
\psi^{(0) s,s'}(\xi)
&= \sum_{\ttn} \psi^{(0) s,s'}_{\ttn} \xi^{-\ttn}
 = \sum_{\ve=0,1} \Phi^s_\ve(\xi) \ot \Psi^{*s'}_\ve(\xi),
\end{aligned}
\end{equation}
where $\Psi^{*s}_{\ve}$ and $\Phi^s_\ve$ are defined by
\eqref{4.12} and \eqref{4.13}.
Again $s,\bs=\pm$ (or equivalently $\pm 1$).

These components act as follows.
\begin{align*}
\psi^{(\half)s}_{\ve,\ttn} &: V(\la_a^{(m-n)})\ot V(\la_b^{(n)})
\rightarrow V(\la_{a+s}^{(m-n)}) \ot V(\la_b^{(n)}),\\
\psi^{(0) s,s'}_{\ttn} &: V(\la_a^{(m-n)}) \ot V(\la_b^{(n)})
\rightarrow V(\la_{a+s}^{(m-n)}) \ot V(\la_{b+s'}^{(n)}).
\end{align*}

Using the commutation relations \eqref{5.129}--\eqref{eq:psmzero},
and \eqref{unit}, we arrive at
\bea
 T^A(\z) \psi^{(\half)\,s}_{\ep}(\xi)\vac_{a,b}&=&
 \psi^{(\half)\,s}_{\ep}(\xi)\vac_{a,n-b},\nn\\[2mm]
T^B(\z) \psi^{(\half)\,s}_{\ep}(\xi)\vac_{a,b}&=&\tau(\z/\xi)
 \psi^{(\half)\,-s}_{\ep}(\xi)\vac_{m-n-a,n-b},\nn\\[2mm]
 T^A(\z) \psi^{(0)\,s,\bs}(\xi)\vac_{a,b}&=&\tau(\z/\xi)
 \psi^{(0)\,s,-\bs}(\xi)\vac_{a,n-b},\nn\\[2mm]
T^B(\z) \psi^{(0)\,s,\bs}(\xi)\vac_{a,b}&=&
 \psi^{(0)\,-s,-\bs}(\xi)\vac_{m-n-a,n-b},\nn
\ena
and hence that
\bea \cT(\z) \psi^{(\half)\,s}_{\ep}(\xi)\vac_{a,b}&=&\tau(\z/\xi)
     \psi^{(\half)\,-s}_{\ep}(\xi)\vac_{m-n-a,b},\nn\\
\cT(\z) \psi^{(0)\,s,\bs}(\xi)\vac_{a,b}&=&\tau(\z/\xi)
     \psi^{(0)\,-s,\bs}(\xi)\vac_{m-n-a,b}.\nn\ena

The vectors $\psi^{(\half)\,s}_{\ep}(\xi)\vac_{a,b}$ and
$\psi^{(0)\,s,\bs}(\xi)\vac_{a,b}$ are
the spin-$\half$ and spin-$0$ eigenstates mentioned in the introduction.
Note that both states are degenerate with respect to $\cT(\z)$, but that
$\psi^{(\half)\,s}_{\ep}(\xi)\vac_{a,b}$ has an eigenvalue of 1 for 
$T^A(\xi)$, and $\psi^{(0)\,s,\bs}_{\ep}(\xi)\vac_{a,b}$ an eigenvalue
of 1 for $T^B(\xi)$. This is consistent with the 
Bethe Ansatz calculations for the
alternating spin-$\half$/spin-$1$ model presented in~\cite{VMN94}.

Further eigenstates may be constructed by acting
with any composition of
$\psi^{(\half)\,s_i}(\xi_i)$ and $\psi^{(0)\,s'_i,\bs'_i}(\xi'_i)$ 
on $\vac_{a,b}$. The eigenvalues of $\cT(\z)$ are the product of all
the $\tau(\z/\xi_i)$ and $\tau(\z/\xi'_i)$ factors.

\subsection{The S-matrix}

The S-matrix for our model is specified by the exchange relations
of $\psi^{(\half)\,s_i}(\xi_i)$ and $\psi^{(0),s'_i,\bs'_i}(\xi'_i)$
with themselves and with each other.
These intertwiners are defined in terms of the
intertwiners $\Phi^s(\z)$ and $\Psi^{*,s}(\z)$ of irreducible modules
in \eqref{IIdefs}. If we act with both sides on the level $\ell$ highest-weight
module $V(\la_r^{(\ell)})$, then  
the commutation relations of the $\Phi^s(\z)$ and $\Psi^{*,s}(\z)$ are
\bea
&&\!\!\!\!\!\!\!\!\!\!\!\!\!\!\!\!\!\!\!\!\!\!\!\!\!\!\!\!\!\!\!
\sli_{\epp_1,\epp_2} \Rbar^{(1,1)}(\xi)^{\epp_1,\epp_2}_{\ep_1,\ep_2}
\Phi^{s_1}_{\epp_1}(\xi_1) \Phi^{s_2}_{\epp_2}(\xi_2) = 
\sli_{s'_1,s'_2}
\Phi^{s'_2}_{\ep_2}(\xi_2)
\Phi^{s'_1}_{\ep_1}(\xi_1)
\wo{r}{r+s_2}{r+s'_1}{r+s_1+s_2}{\xi}{\ell},\label{comI-I}\\[2mm]
&&\!\!\!\!\!\!\!\!\!\!\!\!\!\!\!\!\!\!\!\!\!\!\!\!\!\!\!\!\!\!\!
\Psi^{*\,s_1}_{\ep_1}(\xi_1) \Psi^{*\,s_2}_{\ep_2}(\xi_2)=
\sli_{s'_1,s'_2,\epp_1,\epp_2} 
\Psi^{*\,s'_2}_{\epp_2}(\xi_2) \Psi^{*\,s'_1}_{\epp_1}(\xi_1)
\Rbar^{(1,1)}(\xi)^{\epp_1,\epp_2}_{\ep_1,\ep_2}
\wt{r}{r+s_2}{r+s'_1}{r+s_1+s_2}{\xi}{\ell},\label{comII-II}\\[2mm]
&&\!\!\!\!\!\!\!\!\!\!\!\!\!\!\!\!\!\!\!\!\!\!\!\!\!\!\!\!\!\!\!
\Phi^{s_1}_{\ep_1}(\xi_1) \Psi^{*\,s_2}_{\ep_2}(\xi_2)=
\sli_{s'_1,s'_2}
\Psi^{*\,s'_2}_{\ep_2}(\xi_2) \Phi^{s'_1}_{\ep_1}(\xi_1) 
\wth{r}{r+s_2}{r+s'_1}{r+s_1+s_2}{\xi}{\ell}.\label{comI-II}
\ena
Here, the sum over $s'_1$ and $s'_2$ is restricted to the values for which
$s'_1+s'_2=s_1+s_2$.  $W^I_{\ell}$, $W^{II}_{\ell}$ and $W^*_{\ell}$ are
given in terms
of the RSOS Boltzmann weight $\overline{W}^1_{\ell}$ (given for example in
equation (B.2) in~\cite{JMO93}) as follows:
\bea
\wo{r}{s}{u}{t}{\xi}{\ell}&=&
\frac{X(p^2\xi^{-2})}{X(p^2\xi^{2})}
\;\wec{\la^{(\ell)}_r}{\la^{(\ell)}_s}{\la^{(\ell)}_u}{\la^{(\ell)}_t}{\xi^2} 
\xi^{
\de_{t,s+1} -\,
\de_{r,u-1}},\nn \\[6mm]
\wt{r}{s}{u}{t}{\xi}{\ell}&=&\frac{X(\xi^{-2})}{X(\xi^{2})}
\;\wec{\la^{(\ell)}_r}{\la^{(\ell)}_s}{\la^{(\ell)}_u}{\la^{(\ell)}_t}{\xi^2} 
\xi^{
\de_{t,s+1} -\,
\de_{r,u-1}},\nn \\[6mm]
\wth{r}{s}{u}{t}{\xi}{\ell}&=&
\frac{X(p\xi^{-2})}{X(p\xi^{2})}
\;\wec{\la^{(\ell)}_r}{\la^{(\ell)}_s}{\la^{(\ell)}_u}{\la^{(\ell)}_t}{p^{-1}
\xi^2} 
(-\xi q^{-(1+r)})^{\de_{t,s+1} -\,
\de_{r,u-1}}
\;q^{\de_{r,t}}, \nn
\ena
where 
$\xi=\xi_1/\xi_2$,
$X(z)= \frac{(z;p^2,q^4)_{\infty} (q^4 z;p^2,q^4)_{\infty}}%
{(q^2 z;p^2,q^4)_{\infty}^2}$ and $p=q^{\ell+2}$ 
(note that the $p$ of~\cite{idzal93} is equal to our $p^2$).
Relations \eqref{comI-I} and \eqref{comII-II} come from~\cite{idzal93},
where they were obtained (for a homogeneous evaluation
representation) by solving the $q$-KZ equation.
We obtained \eqref{comI-II} by making use of the
technique mentioned in Proposition A.4 (ii)
of~\cite{idzal93} (and due originally to Okado).

Using these commutation relations, the definitions \eqref{IIdefs}, and
the unitarity property \eqref{unit}, it is then simple to show that on
$V(\la_a^{(m-n)})\ot V(\la_b^{(n)})$ (and hence on $\vac_{a,b}$) we have
\bea
&&\!\!\!\!\!\!\!\!\!\!\!\!\!\!\!
\psi^{(\half)\,s_1}_{\ep_1}(\xi_1) \psi^{(\half)\,s_2}_{\ep_2}(\xi_2)\!\!
=\!\!\!\!\!\!\!
\sli_{\epp_1,\epp_2,s'_1,s'_2}\!\!\!\!\!\!\!
\psi^{(\half)\,s'_2}_{\epp_2}(\xi_2) \psi^{(\half)\,s'_1}_{\epp_1}(\xi_1)
\bar{R}^{(1,1)}(\xi)^{\epp_1,\epp_2}_{\ep_1,\ep_2}\; 
\wt{a}{a+s_2}{a+s'_1}{a+s_1+s_2}{\xi}{m-n},\nn\\
\label{commy1} \\
&&\!\!\!\!\!\!\!\!\!\!\!\!\!\!\!
\psi^{(0)\,s_1,\bs_1}(\xi_1) \psi^{(0)\,s_2,\bs_2}(\xi_2) =
\sli_{s'_1,s'_2,\bs'_1,\bs'_2} 
 \psi^{(0)\,s'_2,\bs'_2}(\xi_2) \psi^{(0)\,s'_1,\bs'_1}(\xi_1) \;
\nn\\&&\times \wo{a}{a+s_2}{a+s'_1}{a+s_1+s_2}{\xi}{m-n}
\wt{b}{b+\tilde{s}_2}{b+\tilde{s}'_1}%
{b+\tilde{s}_1+\tilde{s}_2}{\xi}{n}
,\label{commy2} \\[6mm]
&&\!\!\!\!\!\!\!\!\!\!\!\!\!\!\!
\psi^{(0)\,s_1,\bs_1}(\xi_1) \psi^{(\half)\,s_2}_{\ep_2}(\xi_2) =
\sli_{s'_1,s'_2} \psi^{(\half)\,s'_2}_{\ep_2}(\xi_2) \psi^{(0)\,s'_1,\bs_1}
(\xi_1)
\;
\wth{a}{a+s_2}{a+s'_1}{a+s_1+s_2}{\xi}{m-n}.
\label{commy3}\ena
Again, the sums are restricted so that  $s'_1+s'_2=s_1+s_2$, and 
$\bs'_1+\bs'_2=\bs_1+\bs_2$.

When $n=1$, $m=2$, our model consists of alternating spin-$\half$ and
spin-$1$ lines. In this case we have
\bea
\psi^{(\half)}_{\ep_1}(\xi_1) \psi^{(\half)}_{\ep_2}(\xi_2)
&=& - \sli_{\epp_1,\epp_2} 
\psi^{(\half)}_{\epp_2}(\xi_2) \psi^{(\half)}_{\epp_1}(\xi_1)
R^{(1,1)}(\xi)^{\epp_1,\epp_2}_{\ep_1,\ep_2},\nn\\
\psi^{(0)}(\xi_1) \psi^{(0)}(\xi_2)
&=& - \psi^{(0)}(\xi_2) \psi^{(0)}(\xi_1),\nn\\[1mm]
\psi^{(0)}(\xi_1) \psi^{(\half)}_{\ep_2}(\xi_2) &=&
\tau(\xi)
 \psi^{(\half)}_{\ep_2}(\xi_2) \psi^{(0)}(\xi_1).\nn\ena
Here, the intertwiners act on the tensor product of level-1
irreducible highest weight modules. So, the $s$ and $\bs$
superscripts depend solely on the choice $i$ and $j$, and
we suppress them. These relations are consistent with Bethe
Ansatz calculations of the S-matrix for this example~\cite{VMN94}.


\setcounter{equation}{0}
\section{The Domain Wall Description of the Path Space and the Particle
         Picture}\label{sec7}

\subsection{Domain walls}\label{sec7.1}

Let us now use $\ket{p}$ to denote a double infinite path
$\ket{p}=\cdots \, p(2)\,  p(1)\,  p(0)\, p(-1)\,p(-2)\, \cdots$, for which
\begin{alignat}{2}
p(\tts)&\in \{0, 1, \cdots, n \} \quad& &\text{if $\tts$ is odd,}\\
p(\tts)&\in \{0, 1, \cdots, m \} \quad& &\text{if $\tts$ is even.}
\end{alignat}
Define 
\be
\calP =\oplus_{a,b;a',b'} P_{a,b;a',b'},
\ee
where $P_{a,b;a',b'}$ is the set
\be
P_{a,b;a',b'}=\{\ket{p};\,p(\tts)=\bar{p}(\tts;a,b), \tts\gg0;\,
p(\tts)=\bar{p}(\tts;a',b'), \tts\ll0\}.
\ee
The ground state path $\bar{p}(\tts;a,b)$ was defined by~\eqref{gs2}
(note, however, that $\tts$ may now be negative).

In this section, we construct a domain wall description of the space
$\calP$ and give rules for the induced crystal action on this set of
domain walls.

First, we label a \emph{domain} of a path $\ket{p}\in \calP$ by a 
pair of integers  $(a,b)$, which can take the values
$0\leq a\leq m-n$ and  $0\leq b\leq n$. 
Suppose we start with a path $\ket{p}\in \calP$ and try to associate a 
particular domain $(a(\tts),b(\tts))$ with each $\tts$, such that
$p(\tts)=\bar{p}(\tts;a(\tts),b(\tts))$. There are clearly 
different choices of how to do this.
For example, suppose we choose $\tts\equiv0\pmod4$.
Then because $\bar{p}(\tts;a,b)=a+b$, we could associate any of the domains
$(p(\tts)-b,b)$, such that $0\leq p(\tts)-b\leq m-n$ and  $0\leq b\leq n$,
with $\tts$.

In order to fix uniquely which domain $(a(\tts),b(\tts))$
to associate with a  particular
$\tts$ such that $p(\tts)=\bar{p}(\tts;a(\tts),b(\tts))$, we use the
following rules.
\begin{enumerate}
\item Choose $\tts$ odd. 
\item If $n\leq p(\tts+1)+p(\tts)\leq m$, let
\begin{align}
b(\tts+1)&=b(\tts)=n-p(\tts),\label{7.5}\\
a(\tts+1)&=a(\tts)=
\begin{cases}
m-p(\tts)-p(\tts+1)& \text{if $\tts\equiv1\pmod4$;}\\
p(\tts)+p(\tts+1)-n& \text{if $\tts\equiv3\pmod4$.}
\end{cases}\label{7.6}
\end{align}
\item If  $p(\tts+1)+p(\tts)>m$, let
\begin{align}
b(\tts+1)&=p(\tts+1)-m+n, \ws b(\tts)=n-p(\tts),\\
a(\tts+1)&=a(\tts)=
\begin{cases}
0   &\text{if $\tts\equiv1\pmod4$;}\\
m-n &\text{if $\tts\equiv3\pmod4$.}
\end{cases}
\end{align}
\item If $p(\tts+1)+p(\tts)<n$, let
\begin{align}
b(\tts+1)&=p(\tts+1), \ws b(\tts)=n-p(\tts),\\
a(\tts+1)&=a(\tts)=
\begin{cases}
m-n &\text{if $\tts\equiv1\pmod4$;}\\
0   &\text{if $\tts\equiv3\pmod4$.}
\end{cases}
\end{align}
\end{enumerate}
By following these rules for all odd $\tts$, we can associate
a unique domain $(a(\tts),b(\tts))$ for all $\tts\in \Z$.
Then $p(\tts)$ will be given by 
\be p(\tts)=\bar{p}(\tts;a(\tts),b(\tts)).\label{ident4}\ee
We can write the resulting $(a(\tts),b(\tts))_{\tts\in \Z}$ as a sequence
of domains $(a_{N+1},b_{N+1})\cdots (a_1,b_1)$ and a 
sequence of integers $\tts_N > \tts_{N-1} > \cdots > \tts_{1}$.
The identification is that
\be
(a(\tts),b(\tts))=(a_i,b_i)\hb{ for } \tts_i\geq \tts > \tts_{i-1}
\hb{ (with } \tts_{N+1}=\infty, \tts_0=-\infty).\label{ident5}
\ee

\begin{dfn}
Let $\calD$ be the set of elements, each of which is specified by a
domain sequence $(a_{N+1},b_{N+1})\cdots (a_1,b_1)$ and integers
$\tts_N>\tts_{N-1}>\cdots > \tts_1$, 
where $N\in \Z_{\geq 0}$, $0\leq a_i\leq m-n$, $0\leq b_i\leq n$, 
$(a_{i+1},b_{i+1})\neq (a_i,b_i)$,  and
\be &\tts_i\in\left\{
\br{lll}
& 2\Z\cup (1+4\Z)  \ws&\hb{ if    }\ws a_{i+1}=a_i=0, \; b_{i+1}> b_i,\\
& 2\Z \cup (1+4\Z) \ws&\hb{ if    }\ws a_{i+1}=a_i=m-n, \; b_{i+1}< b_i,\\
& 2\Z \cup (3+4\Z) \ws&\hb{ if    }\ws a_{i+1}=a_i=0, \; b_{i+1}< b_i,\\
& 2\Z \cup (3+4\Z) \ws&\hb{ if    }\ws a_{i+1}=a_i=m-n, \; b_{i+1}> b_i,\\ 
& 2\Z &\hb{ otherwise}.
\er\right.
\ee
\end{dfn}

Then rules (1)-(4) and equation~\eqref{ident5} specify an injection
$M_1:\calP\ra \calD$, and 
\eqref{ident4} specifies a map $M_2:\calD\ra\calP$ which is the
left inverse of $M_1$, i.e., $M_2\circ M_1 =\id_{\calP}$.
\begin{prop}
$M_1:\calP\ra \calD$ is a bijection.
\end{prop}
\begin{proof}
We will prove that the left inverse $M_2:\calD\ra\calP$ is an 
injection, from which the proposition follows.
Consider two elements,  $(a(\tts),b(\tts))_{\tts\in \Z}$ and 
$(a'(\tts),b'(\tts))_{\tts\in \Z}$  of $\calD$
(we can specify them in this way by making use of~\eqref{ident5}).
Let $\tts_0\equiv1\pmod4$. Then from the definition of
$\calD$, one of the following must be true
\be 
{\rm I.} & a(\tts_0+1)=a(\tts_0),\ws & b(\tts_0+1)=b(\tts_0),\nn\\
{\rm II.} & a(\tts_0+1)=a(\tts_0)=0,\ws & b(\tts_0+1)>b(\tts_0),\nn\\
{\rm III.} & a(\tts_0+1)=a(\tts_0)=m-n,\ws & b(\tts_0+1)<b(\tts_0).\nn\ee
One of three similar conditions must hold for 
$a'(\tts_0+1)$, $a'(\tts_0)$, $b'(\tts_0+1)$,
and $b'(\tts_0)$.
Under the map $M_2$ we have
$(a(\tts),b(\tts))_{\tts\in \Z}\to
\big(\bar{p}(\tts;a(\tts),b(\tts))\big)_{\tts\in \Z}$. 
The requirements that
\be
\bar{p}(\tts_0;a(\tts_0),b(\tts_0))
  &=&\bar{p}(\tts_0;a'(\tts_0),b'(\tts_0)),\\
\bar{p}(\tts_0+1;a'(\tts_0+1),b'(\tts_0+1))
  &=&\bar{p}(\tts_0+1;a'(\tts_0+1),b'(\tts_0+1))
\ee
are equivalent to
\be
b(\tts_0)&=& b'(\tts_0),\label{con1}\\
a(\tts_0+1)-b(\tts_0+1)&=&b'(\tts_0+1)-a'(\tts_0+1)\label{con2}
\ee
respectively.
Combining~\eqref{con1}, \eqref{con2}, one of I, II, III for
$a(\tts_0+1)$, $a(\tts_0)$, $b(\tts_0+1)$, $b(\tts_0)$ and one of the
similar conditions I, II, III for $a'(\tts_0+1)$, $a'(\tts_0)$,
$b'(\tts_0+1)$, $b'(\tts_0)$, we get nine
possible sets of equations in eight unknowns.
It is only possible to get a solution to three of these sets of
equations, namely those we get when
$a(\tts_0+1)$, $a(\tts_0)$, $b(\tts_0+1)$, $b(\tts_0)$ and
$a'(\tts_0+1)$, $a'(\tts_0)$, $b'(\tts_0+1)$, $b'(\tts_0)$ both satisfy I, 
or both satisfy II, or both satisfy III. The single solution for all three
sets is 
\be
a(\tts_0+1)=a'(\tts_0+1),\ws
a(\tts_0)=a'(\tts_0),\ws
b(\tts_0+1)=b'(\tts_0+1),\ws
b(\tts_0)=b'(\tts_0).\label{sol}
\ee
A similar argument leads to the same solution~\eqref{sol} in the case
when $\tts_0\equiv3\pmod4$.
This completes the proof.
\end{proof}

The next step is to understand the induced crystal action on $\calD$.
If we refer to the position at which two domains meet as a
\emph{domain wall}, then the general picture is that the crystal action
moves domain walls around.
In order to describe this action we first identify certain types of domain
wall as \emph{elementary}.
The following is a complete list of elementary walls.


\begin{tabular}{|l|l|l|}\hline\hline
$(a_{i+1},b_{i+1})(a_i,b_i)$ &  $\tts_i$& symbol \\ \hline\hline
$(a-1,b)(a,b)$ & 0 mod 4& $\Com$ \\
$(a+1,b)(a,b)$ & 0 mod 4& $\Czp$ \\
$(0,b+1)(0,b)$ & 0 mod 4 & $\Bzp$\\
$(m-n,b-1)(m-n,b)$ & 0 mod 4& $\Tom$\\
$(a\pm 1,b\mp 1)(a,b)$ & 0 mod 4& $\bpmmp$ \\
\hline
$(0,b+1)(0,b)$ & 1 mod 4 & $\Bop$\\
$(m-n,b-1)(m-n,b)$ & 1 mod 4& $\Tzm$\\
\hline
$(a-1,b)(a,b)$ & 2 mod 4& $\Czm$ \\
$(a+1,b)(a,b)$ & 2 mod 4& $\Cop$ \\
$(0,b-1)(0,b)$ & 2 mod 4 & $\Bom$\\
$(m-n,b+1)(m-n,b)$ & 2 mod 4 & $\Tzp$\\
$(a\pm 1,b\pm 1)(a,b)$ & 2 mod 4& $\bpmpm$ \\
\hline
$(0,b-1)(0,b)$ & 3 mod 4 & $\Bzm$\\
$(m-n,b+1)(m-n,b)$ & 3 mod 4 & $\Top$\\
\hline
\end{tabular}
\vspace*{5mm}

\noindent We write $\,[\,$ to mean either of $\,\Bb\,$ or $\,\Tb\,$.
We shall refer to $\,|_\ssj^\ssa\,$, $\,[_\ssj^\ssb\,$ as spin-$\half$
elementary walls, and to $\,\bu^{\ssa,\ssb}\,$ as spin-$0$ elementary walls. 

We wish to decompose each domain wall of an element in $\calD$
into \emph{ordered} elementary domain walls.
We use the fact that when $\tts_i$ is even, $(a_{i+1},b_{i+1})(a_i,b_i)$
may be written in terms of a unique sequence of intermediate domains such that
the corresponding intermediate domain walls are elementary, and ordered as
\be
\eo \cdots \eo \co \cdots \co \bu \cdots \bu    \quad \hb{or} \quad
 \ez \cdots \ez \cz \cdots \cz \bu \cdots \bu ,
\ee
where the $\bu$ are taken to be of one kind only.
When $\tts_i$ is odd, $(a_{i+1},b_{i+1})(a_i,b_i)$ may also be written uniquely
in terms of an ordered sequence of elementary domain walls of the form
\be
\co \cdots \co  \quad \hb{or} \quad \cz \cdots \cz .
\ee
The whole sequence of ordered elementary walls will then said to have been
\emph{normally ordered}.
It is simple to prove the uniqueness of these ordered decompositions,
but perhaps more illuminating to consider some examples.
\begin{quote}
1) $m=6$, $n=2$, $\tts_i\equiv0\pmod4$:\\
$(0,0)(3,1) = (0,0)(1,0)(2,0)(3,0)(4,0)(3,1) \sim \Com\Com\Com\Com\bpm$,\\
$(1,2)(1,0) = (1,2)(0,2)(0,1)(1,0) \sim \Czp\Bzp\bmp$,\\
$(4,0)(4,2) = (4,0)(4,1)(4,2) \sim \Tom\Tom$.

2) $m=6$, $n=2$, $\tts_i\equiv2\pmod4$:\\
$(0,0)(3,1) = (0,0)(1,0)(2,0)(3,1) \sim \Czm\Czm\bmm$,\\
$(1,2)(1,0) = (1,2)(2,2)(3,2)(2,1)(1,0) \sim \Czm\Czm\bpp\bpp$,\\
$(4,0)(4,2) = (4,0)(3,0)(2,0)(3,1)(4,2) \sim \Cop\Cop\bmm\bmm$.

3) $m=6$, $n=2$, $\tts_i\equiv1\pmod4$:\\
$(0,2)(0,0) = (0,2)(0,1)(0,0) \sim \Bop\Bop$,\\
$(4,0)(4,2) = (4,0)(4,1)(4,2) \sim \Tzm\Tzm$.

4) $m=6$, $n=2$, $\tts_i\equiv3\pmod4$:\\
$(0,0)(0,2) = (0,0)(0,1)(0,2) \sim \Bzm\Bzm$,\\
$(4,2)(4,0) = (4,2)(4,1)(4,0) \sim \Top\Top$.
\end{quote}

Explicitly, the ordered walls turn out as follows.
{\allowdisplaybreaks
\begin{itemize}
\item[] \underline{$(a_2, b_2)(a_1, b_1)$ at $\tts\equiv0\pmod4$.}
  \begin{align}
  &(\Czp)^{a_2}(\Bzp)^{b_2-b_1-a_1}(\bmp)^{a_1}
   \quad\text{if $b_2-b_1 > a_1$;}\label{e13}\\
  &(\Czp)^{a_2+b_2-a_1-b_1}\left\{
   \begin{matrix}
   (\bmp)^{b_2-b_1} \;\text{($b_2\geq b_1$);}\\
   (\bpm)^{b_1-b_2} \;\text{($b_2\leq b_1$)}
   \end{matrix}
   \right.
   \quad\text{if $a_2+b_2\geq a_1+b_1\geq b_2$;}\label{e14}\\
  &(\Com)^{a_1+b_1-a_2-b_2}\left\{
   \begin{matrix}
   (\bmp)^{b_2-b_1} \;\text{($b_2\geq b_1$);}\\
   (\bpm)^{b_1-b_2} \;\text{($b_2\leq b_1$)}
   \end{matrix} 
   \right.
   \quad\text{if $m-n+b_2\geq a_1+b_1\geq a_2+b_2$;}\label{e15}\\
  &(\Com)^{m-n-a_2}(\Tom)^{a_1+b_1-b_2-m+n}(\bpm)^{m-n-a_1}
   \quad\text{if $b_1-b_2 > m-n-a_1$.}\label{e16}
  \end{align}
\item[] \underline{$(0, b_2)(0, b_1)$ at $\tts\equiv1\pmod4$.}
  \begin{equation}
  (\Bop)^{b_2-b_1}
  \phantom{ ()^{m-n-a_2}(\Tom)^{a_1+b_1-b_2-m+n}(\bpm)^{m-n-a_1}
            \quad\text{if $b_1-b_2 > m-n-a_1$.} }
  \label{e21}
  \end{equation}
\item[] \underline{$(m-n, b_2)(m-n, b_1)$ at $\tts\equiv1\pmod4$.}
  \begin{equation}
  (\Tzm)^{b_1-b_2}
  \phantom{ ()^{m-n-a_2}(\Tom)^{a_1+b_1-b_2-m+n}(\bpm)^{m-n-a_1}
            \quad\text{if $b_1-b_2 > m-n-a_1$.} }
  \label{e22}
  \end{equation}
\item[] \underline{$(a_2, b_2)(a_1, b_1)$ at $\tts\equiv2\pmod4$.}
  \begin{align}
  &(\Cop)^{a_2}(\Bom)^{b_1-b_2-a_1}(\bmm)^{a_1}
   \quad\text{if $b_1-b_2 > a_1$;}\label{e17}\\
  &(\Cop)^{b_1-a_1-b_2+a_2}\left\{
   \begin{matrix}
   (\bmm)^{b_1-b_2} \;\text{($b_2\leq b_1$);}\\
   (\bpp)^{b_2-b_1} \;\text{($b_2\geq b_1$)}
   \end{matrix}
   \right.
   \quad\text{if $b_2\geq b_1-a_1\geq b_2-a_2$;}\label{e18}\\
  &(\Czm)^{b_2-a_2-b_1+a_1}\left\{
   \begin{matrix}
   (\bmm)^{b_1-b_2} \;\text{($b_2\leq b_1$);}\\
   (\bpp)^{b_2-b_1} \;\text{($b_2\geq b_1$)}
   \end{matrix}
   \right.
   \quad\text{if $b_2-a_2\geq b_1-a_1\geq b_2-m+n$;}\label{e19}\\
  &(\Czm)^{m-n-a_2}(\Tzp)^{b_2-b_1+a_1-m+n}(\bpp)^{m-n-a_1}
   \quad\text{if $b_2-b_1 > m-n-a_1$.}\label{e20}
  \end{align}
\item[] \underline{$(0, b_2)(0, b_1)$ at $\tts\equiv3\pmod4$.}
  \begin{equation}
  (\Bzm)^{b_1-b_2}
  \phantom{ ()^{m-n-a_2}(\Tom)^{a_1+b_1-b_2-m+n}(\bpm)^{m-n-a_1}
            \quad\text{if $b_1-b_2 > m-n-a_1$.} }
  \label{e23}
  \end{equation}
\item[] \underline{$(m-n, b_2)(m-n, b_1)$ at $\tts\equiv3\pmod4$.}
  \begin{equation}
  (\Top)^{b_2-b_1}
  \phantom{ ()^{m-n-a_2}(\Tom)^{a_1+b_1-b_2-m+n}(\bpm)^{m-n-a_1}
            \quad\text{if $b_1-b_2 > m-n-a_1$.} }
  \label{e24}
  \end{equation}
\end{itemize}
} 
After normally ordering all the walls in an element $\ket{d}\in\calD$,
the rules for the crystal action are relatively simple
(we give the rule for the action of $\tf_i$, the action of $\te_i$
can be reconstructed in terms of the inverse of this).
Suppose we have a total of $N$ elementary spin-$\half$ walls
with subscript $\ssj_N,\ssj_{N-1},\cdots,\ssj_{1}$ (and any number of
spin-$0$ walls). 
Now consider the vector
$[\ssj_N]^{(1)}\ot [\ssj_{N-1}]^{(1)}\ot\cdots\ot [\ssj_{1}]^{(1)}
 \in (B^{(1)})^{\ot N}$.
The operator $\tf_i$ acts on 
$[\ssj_N]^{(1)}\ot [\ssj_{N-1}]^{(1)}\ot\cdots\ot [\ssj_{1}]^{(1)}$
by changing a single $\ssj_j \ra {1-\ssj_j}$ (or by
sending the vector to zero). 
Which $\ssj_j$ is changed depends on whether $i=0$ or $1$.
The action of $\tf_i$ on the element $\ket{d}$ 
is to change only the single  
elementary spin-$\half$ wall with the corresponding $\ssj_j$ index (or it
sends the path to zero).
The change that occurs for this elementary domain wall depends on its type
and position $\tts$ in the following way:
\begin{alignat}{7}{}
\tts+2\phantom{|||}& & \tts+1\phantom{|}&
  & \tts\phantom{||||}& & & & \tts+2\phantom{||||||}&
  & \tts+1\phantom{|||}& & \tts\phantom{|||}& \notag \\
(\cdots)\underbrace{\bu\cdots\bu}_c& & & & |_{\ssj}(\cdots)& & \ra&
  & (\cdots)|_{1-\ssj}\underbrace{\bu\cdots\bu}_{c}& & & & (\cdots)&
  \label{7.25} \\
(\cdots)\underbrace{\bu\cdots\bu}_c& & & & |_{\ssj}(\cdots)& & \ra&
  & (\cdots)[_{1-\ssj}\underbrace{\bu\cdots\bu}_{c-1}& & & & (\cdots)&
  \label{7.24} \\
& & (\cdots)\phantom{|}& & \quad [_{\ssj}(\cdots)& & \quad\ra&\quad
  & \phantom{\underbrace{\bu\cdots\bu}_{c}}&
  & \quad(\cdots)[_{1-\ssj}& & \quad(\cdots)&
  \label{7.26} \\
(\cdots)\underbrace{\bu\cdots\bu}_{c}& & [_{\ssj}(\cdots)& & & & \ra&
  & (\cdots)|_{1-\ssj}\underbrace{\bu\cdots\bu}_{c+1}&
  & (\cdots)\phantom{|||}& & &
  \label{7.28} \\
(\cdots)\underbrace{\bu\cdots\bu}_{c}& & \quad[_{\ssj}(\cdots)& & & & \ra&
  & (\cdots)[_{1-\ssj}\underbrace{\bu\cdots\bu}_{c}&
  & (\cdots)\phantom{|||}& & &
  \label{7.27}
\end{alignat}
Here, we have taken $\tts$ to be even.
Also, for~\eqref{7.24} and~\eqref{7.27}, we are assuming that the domain
appearing on the left of $\,\bu\cdots\bu\,$ is at the appropriate boundary,
i.e., $(0,*)$ if $\,[_{1-\ssj}\,$ that may appear at the position $\tts+2$
is a $\,\Bb\,$ and $(m-n,*)$ if it is a $\,\Tb\,$.

Before showing how these rules for the crystal action were obtained,
let us give some simple examples of how this general rule works.
The following two examples capture the two possible ways in which a 
spin-$\half$ wall can pass a spin-$0$ wall under the crystal action. 

First, suppose $m=6$, $n=2$ and that we have an element of $\calD$
described by 3 domains $(4,2)(3,1)(4,1)$ and positions $\tts_2=2$,
$\tts_1=0$.
The $(4,2)(3,1)$ wall at $\tts_2=2$ is a $\bpp$ elementary wall.
The $(3,1)(4,1)$ wall at $\tts_1=0$ is a $\Com$ elementary wall.
Using $M_2$ (as specified by~\eqref{ident4}),
we can write out the section of path in which these walls lie.
The path is
\begin{equation}
\raisebox{-0.37\height}{\begin{texdraw}
\fontsize{10}{13}\selectfont
\setunitscale 1
\drawdim mm
\textref h:C v:C
\move(0 0)
\bsegment
\setsegscale 1.2
\htext(-10 0){$\cdots$} \htext(0 0){$6$} \htext(10 0){$0$} \htext(20 0){$2$}
\htext(30 0){$0$} \htext(40 0){$6$} \htext(50 0){$0$} \htext(55 0){$\bu^{++}$}
\linewd 0 \rmove(-1.7 -0.4)\rlvec(0 0) 
\htext(60 0){$2$} \htext(70 0){$1$} \htext(75 0){$|_1^-$} \htext(80 0){$5$}
\htext(90 0){$1$} \htext(100 0){$\cdots$}
\move(102 2)\move(-12 -2.5)
\esegment
\end{texdraw}}
\end{equation}
Using the above rules for the crystal action on elementary domain walls we 
find that $\tf_1$ sends this path to $0$, and $\tf_0$ sends it to
\begin{equation}\label{dec}
\raisebox{-0.37\height}{\begin{texdraw}
\fontsize{10}{13}\selectfont
\setunitscale 1
\drawdim mm
\textref h:C v:C
\move(0 0)
\bsegment
\setsegscale 1.2
\htext(-10 0){$\cdots$} \htext(0 0){$6$} \htext(10 0){$0$} \htext(20 0){$2$}
\htext(30 0){$0$} \htext(40 0){$6$} \htext(50 0){$0$} \htext(55 0){$\Tzp$}
\htext(60 0){$1$} \htext(70 0){$1$} \htext(80 0){$5$} \htext(90 0){$1$}
\htext(100 0){$\cdots$}
\linewd 0 \rmove(-0.2 0.2)\rlvec(0 0) 
\move(102 2)\move(-12 -2.5)
\esegment
\end{texdraw}}
\end{equation}
Here, we have used~\eqref{7.24}.
This is a path associated with  a single domain wall $(4,2)(4,1)$ at
$\tts_1=2$. 
Working out the subsequent action of $\tf_1$, $\tf_0$, and so on, one finds,
\begin{equation}
\raisebox{-0.48\height}{\begin{texdraw}
\fontsize{10}{13}\selectfont
\setunitscale 1
\drawdim mm
\textref h:C v:C
\move(0 0)
\bsegment
\setsegscale 1.2
\htext(-7 0){$\cdots$} \htext(0 0){$6$} \htext(10 0){$0$} \htext(20 0){$2$}
\htext(30 0){$0$} \htext(40 0){$6$} \htext(50 0){$0$} \htext(55 0){$\bu^{++}$}
\htext(60 0){$2$} \htext(70 0){$1$} \htext(75 0){$|_1^-$} \htext(80 0){$5$}
\htext(90 0){$1$} \htext(97 0){$\cdots$}
\linewd 0 \move(53.0 -0.5)\rlvec(0 0) 
\esegment
\move(0 -7)
\bsegment
\setsegscale 1.2
\htext(-15 1.5){$\xrightarrow{\tf_0}$} \htext(-7 0){$\cdots$} \htext(0 0){$6$}
\htext(10 0){$0$} \htext(20 0){$2$} \htext(30 0){$0$} \htext(40 0){$6$}
\htext(50 0){$0$} \htext(55 0){$\Tzp$} \htext(60 0){$1$} \htext(70 0){$1$}
\htext(80 0){$5$} \htext(90 0){$1$} \htext(97 0){$\cdots$}
\esegment
\move(0 -14)
\bsegment
\setsegscale 1.2
\htext(-15 1.5){$\xrightarrow{\tf_1}$} \htext(-7 0){$\cdots$} \htext(0 0){$6$}
\htext(10 0){$0$} \htext(20 0){$2$} \htext(30 0){$0$} \htext(40 0){$6$}
\htext(45 0){$\Top$} \htext(50 0){$1$} \htext(60 0){$1$} \htext(70 0){$1$}
\htext(80 0){$5$} \htext(90 0){$1$} \htext(97 0){$\cdots$}
\esegment
\move(0 -21)
\bsegment
\setsegscale 1.2
\htext(-15 1.5){$\xrightarrow{\tf_0}$} \htext(-7 0){$\cdots$} \htext(0 0){$6$}
\htext(10 0){$0$} \htext(20 0){$2$} \htext(30 0){$0$} \htext(35 0){$\Czp\bmp$}
\htext(40 0){$5$} \htext(50 0){$1$} \htext(60 0){$1$} \htext(70 0){$1$}
\htext(80 0){$5$} \htext(90 0){$1$} \htext(97 0){$\cdots$}
\esegment
\move(0 -28)
\bsegment
\setsegscale 1.2
\htext(-15 1.5){$\xrightarrow{\tf_1}$} \htext(-7 0){$\cdots$} \htext(0 0){$6$}
\htext(10 0){$0$} \htext(15 0){$\Cop$} \htext(20 0){$3$} \htext(30 0){$0$}
\htext(35 0){$\bmp$} \htext(40 0){$5$} \htext(50 0){$1$} \htext(60 0){$1$}
\htext(70 0){$1$} \htext(80 0){$5$} \htext(90 0){$1$} \htext(97 0){$\cdots$}
\esegment
\move(-21 2)\move(116 -31)
\end{texdraw}}\label{action1}
\end{equation}
etc.
Here we have just shown the elementary wall decomposition at each stage.
The final sequence of domains is $(4,2)(3,2)(4,1)$.
Notice, the action of $\tf_0$ on the fourth line used~\eqref{7.28} and
not~\eqref{7.27}.
This is because the domain on the left of the (non-existent)
$\,\bu\cdots\bu\,$ is at
the boundary, but not the relevant one.

Now consider the rather similar example when $m=6,n=2$ and we have an element
$\ket{d}\in\calD$ described by 3 domains $(3,2)(2,1)(3,1)$ and positions
$\tts_2=2$, $\tts_1=0$.
Again, the $(3,2)(2,1)$ wall at $\tts_2=2$ is a $\bpp$ elementary
wall and the $(2,1)(3,1)$ wall at $\tts_1=0$ is a $\Com$ elementary wall.
Using the rules for the arrows, we get the following sequence.
\begin{equation}
\raisebox{-0.48\height}{\begin{texdraw}
\fontsize{10}{13}\selectfont
\setunitscale 1
\drawdim mm
\textref h:C v:C
\move(0 0)
\bsegment
\setsegscale 1.2
\htext(-7 0){$\cdots$} \htext(0 0){$5$} \htext(10 0){$0$} \htext(20 0){$3$}
\htext(30 0){$0$} \htext(40 0){$5$} \htext(50 0){$0$} \htext(55 0){$\bpp$}
\htext(60 0){$3$} \htext(70 0){$1$} \htext(75 0){$\Com$} \htext(80 0){$4$}
\htext(90 0){$1$} \htext(97 0){$\cdots$}
\linewd 0 \move(53.0 -0.5)\rlvec(0 0) 
\esegment
\move(0 -7)
\bsegment
\setsegscale 1.2
\htext(-15 1.5){$\xrightarrow{\tf_0}$} \htext(-7 0){$\cdots$} \htext(0 0){$5$}
\htext(10 0){$0$} \htext(20 0){$3$} \htext(30 0){$0$} \htext(40 0){$5$}
\htext(50 0){$0$} \htext(55 0){$\Czm\bpp$} \htext(60 0){$2$} \htext(70 0){$1$}
\htext(80 0){$4$} \htext(90 0){$1$} \htext(97 0){$\cdots$}
\esegment
\move(0 -14)
\bsegment
\setsegscale 1.2
\htext(-15 1.5){$\xrightarrow{\tf_1}$} \htext(-7 0){$\cdots$} \htext(0 0){$5$}
\htext(10 0){$0$} \htext(20 0){$3$} \htext(30 0){$0$} \htext(35 0){$\Com$}
\htext(40 0){$6$} \htext(50 0){$0$} \htext(55 0){$\bpp$} \htext(60 0){$2$}
\htext(70 0){$1$} \htext(80 0){$4$} \htext(90 0){$1$} \htext(97 0){$\cdots$}
\esegment
\move(0 -21)
\bsegment
\setsegscale 1.2
\htext(-15 1.5){$\xrightarrow{\tf_0}$} \htext(-7 0){$\cdots$} \htext(0 0){$5$}
\htext(10 0){$0$} \htext(15 0){$\Czm$} \htext(20 0){$3$} \htext(30 0){$0$}
\htext(40 0){$6$} \htext(50 0){$0$} \htext(55 0){$\bpp$} \htext(60 0){$2$}
\htext(70 0){$1$} \htext(80 0){$4$} \htext(90 0){$1$} \htext(97 0){$\cdots$}
\esegment
\move(-21 2)\move(116 -24)
\end{texdraw}}
\end{equation}
The final sequence of domains is $(3,2)(4,2)(3,1)$.
We see that the spin-$0$ wall remains fixed in this case, whereas
in~\eqref{action1} it was moved 2 spaces to the left by the passage of the 
spin-$\half$ wall. 

Let us now show how the rules~\eqref{7.25}--\eqref{7.27} were obtained.
We shall consider $\tf_1$ only.
Let $\ket{p}$ be any path.
Following the rule~\eqref{eq3.1}, we associate a sequence of $1$'s and $0$'s
to each $p(\tts)$.
Then, using the usual rule, we simplify it in such a way that
each domain wall carries $(1)^c$ or $(0)^c$.
This is determined locally at each wall and called the localisation.
It is convenient to think that there always exists a domain wall between
$\tts+1$:odd and $\tts$:even.
If it is not a real one, the localisation is trivial, i.e., $c=0$.
Let us explain this more carefully, starting with two examples.
\begin{itemize}
\item[] \underline{$\tts\equiv1\pmod4$ and $n\leq p(\tts+1)+p(\tts)\leq m$.}\\
  Recalling the process for fixing the domains, given at the beginning of this
  section, we see that, in this case, there is no wall between $p(\tts+1)$
  and $p(\tts)$.
  The domain is $(a,b)$, where $p(\tts+1) = m-n-a+b$ and $p(\tts) = n-b$
  (see~\eqref{7.5} and~\eqref{7.6}).
  We have
  \begin{equation}
  (1)^{m-n-a+b}(0)^{n+a-b}(1)^{n-b}(0)^b\;\sim\;(1)^{m-n-a+b}(0)^{a+b}.
  \end{equation}
  Distribute $(1)^{m-n-a+b}$ to the left wall, and $(0)^{a+b}$ to the right.
\item[] \underline{$\tts\equiv1\pmod4$ and $p(\tts+1)+p(\tts) < n$.}\\
  The domain changes at the centre.
  The domains are given by $(m-n,b_2)(m-n,b_1)$ with $p(\tts+1) = b_2$ and
  $p(\tts) = n-b_1$.
  We have $b_2 < b_1$.
  \begin{equation}
  (1)^{b_2}(0)^{m-b_2}(1)^{n-b_1}(0)^{b_1}\;\sim\;(1)^{b_2}(0)^{m-n+2b_1-b_2}.
  \end{equation}
  Distribute $(1)^{b_2}$ to the left, $(0)^{b_1-b_2}$ to the centre,
  and $(0)^{m-n+b_1}$ to the right.
\end{itemize}
We carry out a similar procedure for all other cases.
Now, consider a wall between $\tts+1$:odd and $\tts$:even.
Suppose $(0)^{c_1}$ is distributed from the left and $(1)^{c_2}$ from the
right.
If $c_1 > c_2$, the localisation is $(0)^{c_1-c_2}$.
If $c_1 = c_2$, there is no (real) wall.
If $c_1 < c_2$, the localisation is $(1)^{c_2-c_1}$.
For a wall between $\tts+1$:even and $\tts$:odd, the localisation is already
given in the form $(0)^c$ or $(1)^c$.
In fact, we have the following simple rule (for the $\tf_1$ case).
\begin{equation}
\begin{matrix}
\text{domain} & \text{position} & \text{localisation}\\
(a_2, b_2)(a_1, b_1) & \tts\equiv0 & (0)^{a_2+b_2}(1)^{a_1+b_1}\\
(0, b_2)(0, b_1) & \tts\equiv1 & (1)^{b_2-b_1}\\
(m-n, b_2)(m-n, b_1) & \tts\equiv1 & (0)^{b_1-b_2}\\
(a_2, b_2)(a_1, b_1) & \tts\equiv2 & (0)^{m-n-a_2+b_2}(1)^{m-n-a_1+b_1}\\
(0, b_2)(0, b_1) & \tts\equiv3 & (0)^{b_1-b_2}\\
(m-n, b_2)(m-n, b_1) & \tts\equiv3 & (1)^{b_2-b_1}
\end{matrix}
\end{equation}
We now consider the action of $\tf_1$.
Suppose it acts on the part of a path $x=p(\tts+1)$ and $y=p(\tts)$ with
$\tts\equiv3\pmod4$.
Suppose that by the action of $\tilde f_1$ we have the change:
\begin{equation}\label{change1}
x\ra x+1.
\end{equation}
In the $1$ and $0$ notation, this part is equivalent to starting from
\begin{equation}
(1)^x(0)^{m-x}(1)^y(0)^{n-y}
\end{equation}
and changing the leftmost $0$ in $(0)^{m-x}$ to $1$.
It implies $m-x>y$.
We have two cases.
\begin{itemize}
\item[] \underline{$n\leq x+y<m$.}\\
  Both $x$ and $y$ belong to the same domain, say, $(a_2,b_2)$.
  We have $x=a_2+b_2$ and $y=n-b_2$.
  Therefore, we have
  \bea\label{EQ1}
  a_2+n=x+y<m.
  \ena
  Let $(a_1,b_1)$ be the domain on the right of $y$, $(a_3,b_3)$ on the left
  of~$x$.

  It is necessary that the localisation at the wall between $(a_2,b_2)$ and
  $(a_1,b_1)$ is $(0)^{c_1}$ with $c_1>0$.
  Therefore, we have $a_1-b_1-a_2+b_2>0$ and this wall is of the
  form~\eqref{e19} or~\eqref{e20}.
  Because of~\eqref{EQ1}, we see that the number of $|_0$ is at least one.
  It is also necessary that the localisation at the wall between
  $(a_3,b_3)$ and $(a_2,b_2)$ is $(1)^{c_2}$ with $c_2\geq0$.
  Therefore, we have $a_3+b_3-a_2-b_2\geq0$ and this wall is of the
  form~\eqref{e15} or~\eqref{e16}.

  The change~\eqref{change1} is equivalent to the change $a_2\ra a_2+1$.
  Using the explicit wall descriptions \eqref{e19}, \eqref{e20}, \eqref{e15},
  and \eqref{e16}, we see that it corresponds to~\eqref{7.25} or~\eqref{7.24}.

\item[] \underline{$x+y<n$.}\\
  There is a wall between $x$ and $y$.
  We have the domains $(a_3,b_3)(0,b_2)(0,b_1)$, where $x$ and $y$ belong
  to $(0,b_2)$ and $(0,b_1)$, respectively.
  We have $x=b_2$ and $y=n-b_1$, and therefore $b_1-b_2=n-x-y>0$.

  The wall between $(0,b_2)$ and $(0,b_1)$ is of the form~\eqref{e23}.
  It is also necessary that the localisation at the wall
  between $(a_3,b_3)$ and $(0,b_2)$ is of the form $(1)^{b_2-a_3-b_3}$.
  Therefore, we have $b_2-a_3-b_3\geq0$, and, in particular, $b_2-b_3\geq0$.
  This wall is of the form~\eqref{e15} (the lower line) or~\eqref{e16}.

  The change~\eqref{change1} is equivalent to the change $b_2\ra b_2+1$.
  It corresponds to~\eqref{7.28} (for~\eqref{e15}) or~\eqref{7.27}
  (for~\eqref{e16}).
\end{itemize}
We may also consider the $\tf_1$ action as sending $y\ra y+1$.
This will bring about the remaining case $x+y\geq m$ and corresponds
to~\eqref{7.26}.
The case $\tts\equiv1\pmod4$ may be similarly analysed to confirm the
results \eqref{7.25}--\eqref{7.27}.

Before ending this section, let us consider one consequence of the rules for
the crystal action on~$\calD$.
Let $\ket{d}\in\calD$ have normally ordered elementary domain walls at
positions $\tts_K, \cdots, \tts_1$.
Define
\begin{equation}\label{7.34}
n(\tts_i,|_{\ssj}) = n(\tts_i,\bu) = -\tts_i/2,
\quad n(\tts_i,[_{\ssj}) = -\tts_i.
\end{equation}
It is simple to check from the rules for the crystal action that 
$\sli_{i=1}^K n(\tts_i,t_i)$ decreases by $1$ under the action of $\tf_i$,
and increases by 1 under the action of $\te_i$.
So, the action of the principal grading operator $\rho$ is given by
\begin{equation}\label{7.34bis}
\rho(\ket{d})=\sli_{i=1}^K n(\tts_i,t_i)\ket{d},
\end{equation}
where $t_i$ refers to the `type' $\,|_{\ssj}\,$, $\,[_{\ssj}\,$, or $\,\bu\,$
of the elementary wall.

\subsection{The particle picture}

In describing a path $\ket{p}\in\calP$ in terms of either local
spin variables $p(\tts)$ or a sequence of domains and domain walls,
we have been using the \emph{local picture}.
We shall now go on to explain the \emph{particle picture} of the
space $\calP\simeq \calD$.
Let $\cD_{a,b;a',b'}$ denote the range of
the restricted map $M_1|_{\cP_{a,b;a',b'}}: \cP_{a,b;a',b'} \ra \cD$.
As a crystal,
$\cD_{a,b;a',b'}$ will decompose into a  (usually infinite) number
of connected components. We wish to understand these connected
components as the crystals created by the creation operators 
$\psi^{(\half)\ssa}_{\ssj,\ttn}$ and $\psi^{({0})\ssa,\ssb}_{\ttn}$
of spin-$\half$ and spin-$0$ particles.
We call this the \emph{particle picture}.
The operators $\psi^{(\half)\ssa}_{\ssj,\ttn}$ and
$\psi^{({0})\ssa,\ssb}_{\ttn}$ will be given as the $q\ra 0$ limit 
of the corresponding operators defined in Section~\ref{sec6} (we conjecture
that this limit is well-defined). 
In the particle picture,
any sequence of the operators $\psi^{(\half)\ssa}_{\ttn}$ and
$\psi^{({0})\ssa,\ssb}_{\ttn}$ is allowed, but with the condition that the
corresponding sequence of the highest weights (represented by $(a_i,b_i)$)
satisfies $0\leq a_i\leq m-n$ and $0\leq b_i\leq n$.
This condition will always be assumed when we talk of a sequence of these
operators.
However, they are not linearly independent because of the commutation
relations \eqref{commy1}--\eqref{commy3} of the particle creation operators.
The particle pictures for the pure spin-$\frac{1}{2}$ model, pure
spin-$\frac{n}{2}$ models and RSOS
fusion models were constructed in references \cite{DFJMN}, 
\cite{naya95a} and \cite{naya95b} respectively.

Before looking at the space spanned by the particles, we prepare some details
about affine crystals.
Suppose we have a $\uq$ crystal $B$ which takes weights in
$P_{cl}= \Z \La_0\oplus\Z \La_1$.
Then the affinization of this crystal, denoted by $\Aff(B)$, takes weights
in $P=\Z\La_1\oplus \Z\La_0 \oplus \Z\de$ (See~\cite{KMN92} for a definition).
Here, we use $\Aff(B)$ defined in the principal gradation.
For example, $\Affo$ is given by either of the following diagrams.
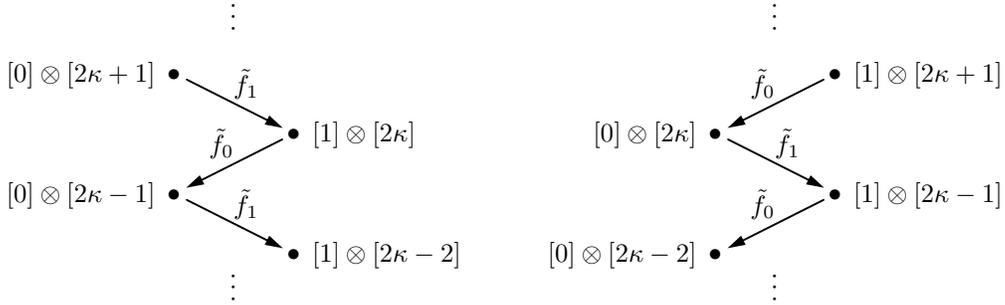
\begin{center}
\begin{texdraw}
\fontsize{10}{13}\selectfont
\setunitscale 0.8
\drawdim mm
\arrowheadsize l:3 w:1.7 \arrowheadtype t:F
\move(0 0)
\bsegment
\textref h:C v:C
\htext(-10 0){$\bu$}\rmove(2 -1)\ravec(16 -8)
\htext(10 -10){$\bu$}\rmove(-2 -1)\ravec(-16 -8)
\htext(-10 -20){$\bu$}\rmove(2 -1)\ravec(16 -8)
\htext(10 -30){$\bu$}\vtext(0 10){$\cdots$}\vtext(0 -35){$\cdots$}
\textref h:R v:C
\htext(-13 0){$[0]\ot[2\kappa+1$]}
\htext(-13 -20){$[0]\ot[2\kappa-1]$}
\textref h:L v:C
\htext(13 -10){$[1]\ot[2\kappa]$}
\htext(13 -30){$[1]\ot[2\kappa-2]$}
\textref h:R v:B
\htext(0 -14){$\tf_0$}
\textref h:L v:B
\htext(0 -4){$\tf_1$}
\htext(0 -24){$\tf_1$}
\esegment
\move(90 0)
\bsegment
\textref h:C v:C
\htext(10 0){$\bu$}\rmove(-2 -1)\ravec(-16 -8)
\htext(-10 -10){$\bu$}\rmove(2 -1)\ravec(16 -8)
\htext(10 -20){$\bu$}\rmove(-2 -1)\ravec(-16 -8)
\htext(-10 -30){$\bu$}\vtext(0 10){$\cdots$}\vtext(0 -35){$\cdots$}
\textref h:L v:C
\htext(13 0){$[1]\ot[2\kappa+1$]}
\htext(13 -20){$[1]\ot[2\kappa-1]$}
\textref h:R v:C
\htext(-13 -10){$[0]\ot[2\kappa]$}
\htext(-13 -30){$[0]\ot[2\kappa-2]$}
\textref h:L v:B
\htext(0 -14){$\tf_1$}
\textref h:R v:B
\htext(0 -4){$\tf_0$}
\htext(0 -24){$\tf_0$}
\esegment
\move(-40 15)\move(130 -40)
\end{texdraw}
\end{center}

Let us now consider the states spanned by just one particle.
 From the definition~\eqref{IIdefs} and the remarks following~\eqref{4.13},
we see that $\psi^{(\half)\ssa}_{\ssj,\ttn}$ is meaningful only if
$\ssa\cdot(-1)^\ssj = (-1)^{\ttn}$ and that $\psi^{(0)\ssa,\ssb}_{\ttn}$ is
meaningful only if $-\ssa\cdot\ssb = (-1)^{\ttn}$.
Considering the degree given by~\eqref{7.34bis} also, we identify
\begin{equation}\label{opart}
\begin{aligned}
|_\ssj^{\ssa} \ \text{at $\tts$}\quad
&\longleftrightarrow \quad\psi^{(\half)\ssa}_{\ssj,-\frac{\tts}{2}},\\
\bu^{\ssa,\ssb} \ \text{at $\tts$}\quad
&\longleftrightarrow \quad\psi^{(0)\ssa,\ssb}_{-\frac{\tts}{2}}.
\end{aligned}
\end{equation}
Recalling the rules for the crystal action on the elementary walls,
we see that each set of $\psi^{(\half)\ssa}_{\ssj,\ttn}$ with $\ssa$ fixed
and other indices satisfying $\ssa\cdot(-1)^\ssj = (-1)^{\ttn}$ brings
about a crystal isomorphic to $\Affo$.
Each set of $\psi^{(0)\ssa,\ssb}_{\ttn}$ with both $\ssa$ and $\ssb$ fixed
is a crystal isomorphic to $\Affz$.

To consider spaces spanned by more than one particle, we have to
study the linear dependence relations in the particle picture
more carefully.
We take the $q\ra 0$ limit of the relations \eqref{commy1}--\eqref{commy3}
and write out the results componentwise.
When acting on $(a,b)$, we have
{\allowdisplaybreaks
\begin{align}
\psi^{(\half)\ssa}_{\ssj,\ttn}\psi^{(\half)\ssa'}_{\ssj',\ttm} &=
- \psi^{(\half)\ssa}_{\ssj,\ttm+\nu}
\psi^{(\half)\ssa'}_{\ssj',\ttn-\nu}\quad\ \,
\text{with $\nu = \de_{\ssa,\ssa'} + \de_{\ssj,\ssj'}$},\label{7.37}\\
\psi^{(0)\ssa,\ssb}_{\ttn}\psi^{(0)\ssa',\ssb'}_{\ttm} &=
- \psi^{(0)\ssa,\ssb}_{\ttm+\nu}
\psi^{(0)\ssa',\ssb'}_{\ttn-\nu}\quad
\text{with $\nu = \de_{\ssa,\ssa'} + \de_{\ssb,\ssb'}$},\\
\psi^{(0)\ssa,\ssb}_{\ttn} \psi^{(\half)\ssa}_{\ssj,\ttm} &=
\psi^{(\half)\ssa}_{\ssj,\ttm} \psi^{(0)\ssa,\ssb}_{\ttn},\label{7.39}\\
\psi^{(0) +,\ssb}_{\ttn} \psi^{(\half)-}_{\ssj,\ttm} &=
\begin{cases}
\psi^{(\half)-}_{\ssj,\ttm} \psi^{(0) +,\ssb}_{\ttn}
&\text{if $a\neq m-n$},\\
\psi^{(\half)+}_{\ssj,\ttm+1} \psi^{(0) -,\ssb}_{\ttn-1}
&\text{if $a= m-n$},
\end{cases}\label{7.40}\\
\psi^{(0) -,\ssb}_{\ttn} \psi^{(\half)+}_{\ssj,\ttm} &=
\begin{cases}
\psi^{(\half)+}_{\ssj,\ttm} \psi^{(0) -,\ssb}_{\ttn}
&\text{if $a\neq 0$},\\
\psi^{(\half)-}_{\ssj,\ttm+1} \psi^{(0) +,\ssb}_{\ttn-1}
&\text{if $a=0$}.\label{7.41}
\end{cases}
\end{align}
} 
Relations \eqref{7.39}--\eqref{7.41} tell us that we may always order the
particles so that all the $\psi^{(0)}$ are to the right of all the
$\psi^{(\half)}$.
Relation~\eqref{7.37} shows that
$\psi^{(\half)\ssa}_{\ssj,\ttn}
 \psi^{(\half)\ssa'}_{\ssj',\ttn-\nu} = 0$.
A little more scrutiny at~\eqref{7.37} shows that we may always order
any nonzero $\psi^{(\half)\ssa}_{\ssj,\ttn}\psi^{(\half)\ssa'}_{\ssj',\ttm}$
so that $\ttn < \ttm$, or in the case $(\ssa,\ssj) = (\ssa',\ssj')$,
$\ttn \leq \ttm$.
A similar statement is true for
$\psi^{(0)\ssa,\ssb}_{\ttn}\psi^{(0)\ssa',\ssb'}_{\ttm}$.
We have shown:

\begin{prop}\label{pr:7.1}
Any sequence of $M$ spin-$\half$ particles $\psi^{(\half)\ssa}_{\ssj,\ttn}$
and $N$ spin-$0$ particles $\psi^{(0)\ssa,\ssb}_{\ttn}$ initiating and
terminating at two given domains may be written in the form
\begin{equation}\label{7.4i}
\psi^{(\half)\ssa'_1}_{\ssj_1,\ttm_1}\cdots\psi^{(\half)\ssa'_M}_{\ssj_M,\ttm_M}
\psi^{(0)\ssa_1,\ssb_1}_{\ttn_1}\cdots\psi^{(0)\ssa_N,\ssb_N}_{\ttn_N},
\end{equation}
modulo sign, if it is not equal to zero.
Here, we require the indices to satisfy
\begin{align}
\ttm_i<\ttm_{i+1} &\ \text{or}\ 
(\ttm_i,\ssa'_i,\ssj_i) = (\ttm_{i+1},\ssa'_{i+1},\ssj_{i+1}),\\
\ttn_i<\ttn_{i+1} &\ \text{or}\ 
(\ttn_i,\ssa_i,\ssb_i) = (\ttn_{i+1},\ssa_{i+1},\ssb_{i+1}).
\end{align}
\end{prop}

A sequence of particles of the form given by this proposition will be called
\emph{separately ordered}.
The name comes from the way the spin-$\half$ particles and spin-$0$
particles have been grouped separately.
This is to be contrasted with the \emph{normally ordered} sequence
to be defined in Section~\ref{sec7.3}.

Let us now consider the vector space which is spanned by the sequence of
particles.
We fix $M$, the number of spin-$\half$ operators
$\psi^{(\half)\ssa}_{\ssj,\ttn}$,
$N$, the number of spin-$0$ operators $\psi^{(0)\ssa,\ssb}_{\ttn}$,
and the initial and final domains.
We will denote the space by $\calA$.
We do not impose the commutation relations \eqref{7.37}--\eqref{7.41} in
$\calA$.
As before, we identify $\psi^{(\half)\ssa}_{\ssj,\ttn}$ and
$\psi^{(0)\ssa,\ssb}_{\ttn}$
with the elements of $\Affo$ and $\Affz$.
Hence the monomial basis of $\calA$ is a crystal isomorphic to a union of
mixed tensor products of $M$-many $\Affo$ and $N$-many $\Affz$.
We call it the \emph{crystal part of $\calA$}.

The subspace spanned by the relations will be denoted by $\calR$.
It is easy to prove, using the tensor product rule for crystals bases,
that the set of relations \eqref{7.37}--\eqref{7.41} is
preserved under the crystal action.
Hence, the monomial basis of $\calA/\calR$ is given a crystal structure.
We call it the \emph{crystal part of $\calA/\calR$}.
We are interested in this crystal structure.

Denote by $\calN$, the set of separately ordered sequence of particles.
It is easy to show that $\calN$ is also preserved under the crystal action.
Hence, $\calN$ is a subcrystal of the crystal part of $\calA$.
We aim to show that $\calN$ forms a basis of $\calA/\calR$ so that the
crystal $\calN$ is, in fact, the crystal part of $\calA/\calR$.

We first give a partial ordering to the set of particles.
Two particles are said to satisfy $\psi^A < \psi^B$ if and only if
$\psi^A\neq\psi^B$ and $\psi^A\psi^B$ is separately ordered.
Then the monomial basis elements of $\calA$ are given the lexicographical
order using the order on the particles.
We define an action of $S_{M+N}$, the symmetric group of order
$M+N$, on $\calA$.
Since all the relations \eqref{7.37}--\eqref{7.41}
are of the form $\psi^A\psi^B = \pm\psi^C\psi^D$,
we may define the action of the transposition $\s_i = (i,i+1)$ on a sequence
of particles by substituting $\psi^A\psi^B$ at the $i$-th and $(i+1)$-th
position with the appropriate $\pm\psi^C\psi^D$.
It is easy to show that this defines an action of $S_{M+N}$ on $\calA$.
We prove two lemmas concerning these definitions.

\begin{lem}\label{lm:small}
Suppose $M+N\geq2$.
Let $A = \psi^{A_1}\cdots\psi^{A_{M+N}}$ and
$B = \psi^{B_1}\cdots\psi^{B_{M+N}}$.
If $A$ is separately ordered, and $\s_1\s_2\cdots\s_{r-1} A = \pm B$
\textup{(}$r\leq M+N$\textup{)}, then $\psi^{A_1} < \psi^{B_1}$.
\end{lem}
\begin{proof}
It suffices to show this for the case $r=M+N$.
We use induction on $r$.
For $r=2$, this may be done by checking each case.
So suppose $r>2$.
Let $\s_{r-1} A = \pm\psi^{A_1}\cdots\psi^{A_{r-2}}\psi^{C}\psi^{B_r}$.
We know from the $r=2$ case that $\psi^{A_{r-1}} < \psi^C$.
Hence, $\psi^{A_1}\cdots\psi^{A_{r-2}}\psi^C$ is separately ordered.
We may now apply induction hypothesis to conclude
$\psi^{A_1} < \psi^{B_1}$.
\end{proof}

\begin{lem}\label{bigger}
Suppose $M+N\geq2$.
Let $A = \psi^{A_1}\cdots\psi^{A_{M+N}}$ and
$B = \psi^{B_1}\cdots\psi^{B_{M+N}}$.
If $A$ is separately ordered, $\pi\in S_{M+N}$ is different from the
identity element, and $\pi A = \pm B$, then $A < B$.
\end{lem}
\begin{proof}
We use induction on $M+N$.
This is easy to check when $M+N = 2$.
If $M+N > 2$ and $\pi(1) = 1$, then we may apply the induction hypothesis
to $A' = \psi^{A_2}\cdots\psi^{A_{M+N}}$ and
$B' = \psi^{B_2}\cdots\psi^{B_{M+N}}$.
So suppose $\pi(r) = 1$ with $r>1$.
Then, we may write $\pi = \pi' \s_1\s_2\cdots\s_{r-1}$ for some
$\pi'\in S_{M+N}$ with $\pi'(1) = 1$.
But, then Lemma~\ref{lm:small} shows, $\psi^{A_1} < \psi^{B_1}$ and
hence $A < B$.
\end{proof}

The next easy corollary to this lemma shows that the expression~\eqref{7.4i}
is unique for each product of particles different from zero.

\begin{coro}\label{coro}
Let $x$ be separately ordered and choose any $\pi\in S_{M+N}$.
Then, $\pi(x)$ is separately ordered if and only if $\pi=\id$.
\end{coro}

We can now finally prove:

\begin{prop}
The set of separately ordered elements, $\calN$, forms a basis for
$\calA/\calR$.
\end{prop}
\begin{proof}
By Proposition~\ref{pr:7.1}, it suffices to show the linear independence
of $\calN$.
Let $(\ |\ )$ denote the natural orthonormal bilinear form on $\calA$.
Define $\bar{x} = \sum_{\pi\in S_{M+N}} \pi(x)$ for any $x\in\calA$.
Noting
\begin{equation*}
\calR = \operatorname{Span}\{ \pi(x) - x ; \pi\in S_{M+N}, x\in\calA \},
\end{equation*}
we have $(\bar{x}|\calR) = 0$ for any $x\in\calN$.
Hence, $(\bar{x}|\,\cdot\,)$ defines a linear functional on $\calA/\calR$.
Using Corollary~\ref{coro},
we may easily check that $(\bar{x}|y)_{x,y\in\calN} = \de_{x,y}$.
This proves that the set $\calN$ is linearly independent.
\end{proof}

So the space described by the particles initiating and terminating at given
domains is the crystal $\calN$ of separately ordered sequence of particles.
We have obtained a clear view of the particle picture given in terms of
the affine crystals $\Affo$ and $\Affz$.

\subsection{Connection between the local and particle pictures}\label{sec7.3}

Let us first describe a map from the domain wall description to the
particle picture.
We have already identified the walls $\,|_\ssj^\ssa\,$ and
$\,\bu^{\ssa,\ssb}\,$ with the particles in~\eqref{opart}.
Writing out the domain wall description in the path form, at
$\tts\equiv 0 \pmod{4}$, we can check
\begin{equation*}
(0,b+1)\Bb^+_0(0,b) = (0,b+1)\bu^{-+}|_0^+(0,b).
\end{equation*}
We may similarly write other $\,[\,$ at even $\tts$ as a combination of
$\,\bu^{\ssa,\ssb}\,$ and $\,|_\ssj^\ssb\,$.
With this and the identification~\eqref{opart}, we map
\begin{equation}\label{twopart}
\begin{aligned}
\Bb_\ssj^\ssb \ \text{at even $\tts$}\quad&\longmapsto\quad
\psi^{(0) -,\ssb}_{-\frac{\tts}{2}}\psi^{(\half)+}_{\ssj,-\frac{\tts}{2}},\\
\Tb_\ssj^\ssb \ \text{at even $\tts$}\quad&\longmapsto\quad
\psi^{(0) +,\ssb}_{-\frac{\tts}{2}}\psi^{(\half)-}_{\ssj,-\frac{\tts}{2}}.
\end{aligned}
\end{equation}
To map the remaining four elementary walls, we return to example
\eqref{action1}.
\begin{equation}
\raisebox{-0.48\height}{\begin{texdraw}
\fontsize{10}{13}\selectfont
\setunitscale 1
\drawdim mm
\textref h:C v:C
\move(0 0)
\bsegment
\setsegscale 1.2
\htext(3 0){$\cdots$} \htext(10 0){$0$} \htext(20 0){$2$}
\htext(30 0){$0$} \htext(40 0){$6$} \htext(50 0){$0$} \htext(55 0){$\bu^{++}$}
\htext(60 0){$2$} \htext(70 0){$1$} \htext(75 0){$|_1^-$} \htext(80 0){$5$}
\htext(87 0){$\cdots$}
\htext(101 0){$\psi^{(0) +,+}_{-1}\psi^{(\half) -}_{1,0}$}
\linewd 0 \move(53.0 -0.5)\rlvec(0 0) 
\esegment
\move(0 -7)
\bsegment
\setsegscale 1.2
\htext(-5 2){$\xrightarrow{\tf_0}$} \htext(3 0){$\cdots$}
\htext(10 0){$0$} \htext(20 0){$2$} \htext(30 0){$0$} \htext(40 0){$6$}
\htext(50 0){$0$} \htext(55 0){$\Tzp$} \htext(60 0){$1$} \htext(70 0){$1$}
\htext(80 0){$5$} \htext(87 0){$\cdots$}
\htext(101 0){$\psi^{(0) +,+}_{-1}\psi^{(\half) -}_{0,-1}$}
\esegment
\move(0 -14)
\bsegment
\setsegscale 1.2
\htext(-5 2){$\xrightarrow{\tf_1}$} \htext(3 0){$\cdots$}
\htext(10 0){$0$} \htext(20 0){$2$} \htext(30 0){$0$} \htext(40 0){$6$}
\htext(45 0){$\Top$} \htext(50 0){$1$} \htext(60 0){$1$} \htext(70 0){$1$}
\htext(80 0){$5$} \htext(87 0){$\cdots$}
\htext(101 0){\fbox{\quad?\quad}\ }
\esegment
\move(0 -21)
\bsegment
\setsegscale 1.2
\htext(-5 2){$\xrightarrow{\tf_0}$} \htext(3 0){$\cdots$}
\htext(10 0){$0$} \htext(20 0){$2$} \htext(30 0){$0$} \htext(35 0){$\Czp\bmp$}
\htext(40 0){$5$} \htext(50 0){$1$} \htext(60 0){$1$} \htext(70 0){$1$}
\htext(80 0){$5$} \htext(87 0){$\cdots$}
\htext(101 0){$\psi^{(\half) +}_{0,-2}\psi^{(0) -,+}_{-2}$}
\esegment
\move(0 -28)
\bsegment
\setsegscale 1.2
\htext(-5 2){$\xrightarrow{\tf_1}$} \htext(3 0){$\cdots$}
\htext(10 0){$0$} \htext(15 0){$\Cop$} \htext(20 0){$3$} \htext(30 0){$0$}
\htext(35 0){$\bmp$} \htext(40 0){$5$} \htext(50 0){$1$} \htext(60 0){$1$}
\htext(70 0){$1$} \htext(80 0){$5$} \htext(87 0){$\cdots$}
\htext(101 0){$\psi^{(\half) +}_{1,-3}\psi^{(0) -,+}_{-2}$}
\esegment
\move(-9 3)\move(129 -32)
\end{texdraw}}
\end{equation}
This time, we have written the particles to the right using~\eqref{opart}
and~\eqref{twopart}.
What should go in the box?
Coming down from the top, we can guess it to be
$\psi^{(0) +,+}_{-1}\psi^{(\half) -}_{1,-2}$.
Going up from the bottom, it should be
$\psi^{(\half) +}_{1,-1}\psi^{(0) -,+}_{-2}$.
We are dealing with the $a=m-n$ case, and as~\eqref{7.40} with
$\ssb=+$ shows, they are actually equal.
Generalising this, we map
\begin{equation}\label{7.46}
\begin{aligned}
\Bb_\ssj^\ssb \ \text{at odd $\tts$}\quad&\longmapsto\quad
\psi^{(\half)+}_{\ssj,-\half(\tts-1)}\psi^{(0) -,\ssb}_{-\half(\tts+1)},\\
\Tb_\ssj^\ssb \ \text{at odd $\tts$}\quad&\longmapsto\quad
\psi^{(\half)+}_{\ssj,-\half(\tts-1)}\psi^{(0) -,\ssb}_{-\half(\tts+1)}.
\end{aligned}
\end{equation}
We have defined a map from the domain wall description to the particle
picture.

We now define the inverse map.
To do this, we construct a new basis of $\calA/\calR$.
We say a sequence of particles is \emph{normally ordered} if each successive
pair is one of the following:
\begin{enumerate}
\item $\psi^{(0)\ssa_1,\ssb_1}_{\ttn_1}\psi^{(0)\ssa_2,\ssb_2}_{\ttn_2}$
      where $\ttn_1<\ttn_2$ or
      $(\ttn_1,\ssa_1,\ssb_1)=(\ttn_2,\ssa_2,\ssb_2)$.
\item $\psi^{(1/2)\ssa_1}_{\ssj_1,\ttn_1}\psi^{(1/2)\ssa_2}_{\ssj_2,\ttn_2}$
      where $\ttn_1<\ttn_2$ or $(\ttn_1,\ssj_1,\ssa_1)=(\ttn_2,\ssj_2,\ssa_2)$.
\item $\psi^{(1/2)\ssa'}_{\ssj,\ttm}\psi^{(0)\ssa,\ssb}_{\ttn}$ where
      $\ttm\leq\ttn$.
\item $\psi^{(1/2)-\ssa}_{\ssj,\ttn}\psi^{(0)\ssa,\ssb}_{\ttn-1}$
      where these are placed at the boundary, i.e., for $\ssa=+$, it acts
      on the domain $(0,*)$, for $\ssa=-$, it acts on the domain $(m-n,*)$.
      \label{it4}
\item $\psi^{(0)\ssa,\ssb}_{\ttn}\psi^{(1/2)\ssa'}_{\ssj,\ttm}$ where
      $\ttn<\ttm$.
\item $\psi^{(0)-\ssa,\ssb}_{\ttn}\psi^{(1/2)\ssa}_{\ssj,\ttn}$
      where these are placed at the boundary, i.e., for $\ssa=+$, it acts
      on the domain $(0,*)$, for $\ssa=-$, it acts on the domain $(m-n,*)$.
      \label{it6}
\end{enumerate}
The set of normally ordered sequence of particles will be denoted by~$\calNo$.
The relations \eqref{7.37}--\eqref{7.41} show that we may always bring any
sequence of particles to a normally ordered sequence.
The linear independence of the normally ordered sequence may be proved as in
the proof for Proposition~\ref{pr:7.1}.
So the normally ordered sequences form a basis for $\calA/\calR$.
The set of normally ordered sequences of particles, $\calNo$, is certainly
a crystal, the crystal action being ``first, act as an element of $\calA$,
then, normally order.''
The map from the particle picture to the local picture may now be taken
by first applying the inverse of~\eqref{twopart} and~\eqref{7.46}
to~(\ref{it6}) and~(\ref{it4}), respectively, and
then applying~\eqref{opart} to the remaining particles.
It is easy to check that the image is an ordered sequence of
elementary domain walls.
The defined map is certainly inverse to the map from the local picture
to the particle picture defined earlier.

\begin{thm}
The local picture and the particle picture are isomorphic as crystals.
\end{thm}
\begin{proof}
It suffices to show that the two maps defined in this section respect
the crystal structures.
So, let us study the rules for the crystal action on $\calNo$.
We shall consider $\tf_1$ only.
The action of $\tf_1$ will change some $\psi^{(\half)\ssa}_{0,\ttm}$
to $\psi^{(\half)\ssa}_{1,\ttm-1}$.
After this change the product may not be normally ordered.
In that case, we must normally order it by using the commutation relations.
The rules come out as follows:

If the product contains $\psi^{(0)-\ssa,+}_{\ttn}\psi^{(\half)\ssa}_{0,\ttn}$
at the boundary, i.e., it acts on the domain $(0,*)$ for the $\ssa=+$ case and
$(m-n,*)$ for the $\ssa=-$ case, the change is
\begin{equation}\label{A}
\psi^{(0)-\ssa,+}_{\ttn}\psi^{(\half)\ssa}_{0,\ttn}
\ra
\psi^{(\half)-\ssa}_{1,\ttn}\psi^{(0)\ssa,+}_{\ttn-1}.
\end{equation}
Otherwise and if the product contains
$(\psi^{(0)-\ssa,-}_{\ttn-1})^c \psi^{(\half)\ssa}_{0,\ttn}$
for some $c\geq1$ and the domain on the left of this part of the product is
at the boundary, the change is
\begin{equation}\label{B}
(\psi^{(0)-\ssa,-}_{\ttn-1})^c \psi^{(\half)\ssa}_{0,\ttn}
\ra
\psi^{(0)-\ssa,-}_{\ttn-1}\psi^{(\half)\ssa}_{1,\ttn-1}
(\psi^{(0)-\ssa,-}_{\ttn-1})^{c-1}.
\end{equation}
Otherwise, let $c\geq0$ be the maximal integer such that
$(\psi^{(0)\ssa',\ssa'\ssa}_{\ttn-1})^c \psi^{(\half)\ssa}_{0,\ttn}$
is contained in the product.
Then, the change is
\begin{equation}\label{C}
(\psi^{(0)\ssa',\ssa'\ssa}_{\ttn-1})^c \psi^{(\half)\ssa}_{0,\ttn}
\ra
\psi^{(\half)\ssa}_{1,\ttn-1}(\psi^{(0),\ssa',\ssa'\ssa}_{\ttn-1})^c.
\end{equation}

In the domain wall language, the case~\eqref{A} corresponds to~\eqref{7.26}.
The case~\eqref{B} corresponds to~\eqref{7.24} and~\eqref{7.27}.
The last case~\eqref{C} corresponds to~\eqref{7.25} and~\eqref{7.28}.
\end{proof}

We have thus related the path space $\calP$ with a crystal given explicitly
in terms of $\Affz$ and $\Affo$.
Namely, we have established the crystal isomorphisms between
$\calP$ and $\calD$, $\calD$ and $\calNo$, $\calNo$ and $\calN$.
And the crystal $\calN$ is given as a union of subcrystals of
$\Affo^{\ot M}\otimes\Affz^{\ot N}$.


\setcounter{equation}{0}
\section{Summary}\label{sec8}

Let us summarise very briefly the main results of our analysis
of infinite-volume alternating-spin vertex models. 
We identify the space on which the transfer
matrices of the alternating spin-$\frac{m}{2}$/
spin-$\frac{n}{2}$ model act as the direct sum of
\bea
&&\Hom_\C\Bigl(V(\la^{(m-n)}_a)\otimes V(\la^{(n)}_b),
V(\la^{(m-n)}_{a'})\otimes V(\la^{(n)}_{b'})\Bigr)\nonumber\\
&&\simeq
V(\la^{(m-n)}_{a'})\otimes V(\la^{(n)}_{b'})
\otimes
\Bigl(V(\la^{(m-n)}_a)\otimes V(\la^{(n)}_b)\Bigr)^*.
\ena
The transfer matrices themselves are constructed in terms of certain
$\uqp$ intertwiners defined on this space (see~\eqref{Tdef1}).
These transfer matrices are diagonalised by making
use of another set of intertwiners given by~\eqref{IIdefs}. 
The vacua are given by $(-q)^D$; the excited states are 
multi-particle states
consisting of a number of spin-$0$ particles
and a number of spin-$\half$ particles.
The two-particle S-matrices are given by~\eqref{commy1} to~\eqref{commy3}.

In~\cite{HKMW98a}, we show how to construct correlation functions
of these models. We derive there the relation between simple
correlation functions of the alternating model and those of the pure
spin-$\frac{n}{2}$ and pure spin-$\frac{m}{2}$ models.
In this, and in the diagonalisation
of the transfer matrix, we make use of the commutativity of one of
our intertwiners (see Section~\ref{sec5} of the current paper) with the
action of the deformed Virasoro algebra considered in~\cite{JS97}. 

In Sections~\ref{sec3} and~\ref{sec7} we consider the crystal limit 
(i.e., $q\ra 0$ limit) of our model in detail. 
In this limit, the corner transfer matrix acts diagonally on
the (half-infinite) path space $P_{a,b}$ associated with a particular
boundary condition $(a,b)$. We prove that there is a crystal 
isomorphism $P_{a,b}\simeq B(\lambda^{(m-n)}_a)\ot
B(\lambda^{(n)}_b)$.
We go on to consider the double infinite path space $\calP$.
We construct a crystal isomorphism between this space and the
space $\calD$ defined in terms of domain walls. $\calP$ and
$\calD$ are both considered as \emph{local picture} descriptions of
the space. We then construct two \emph{particle picture} descriptions,
$\calNo$ and $\calN$, by making use of the $q\ra 0$ limit of
the intertwiners which diagonalise our transfer matrix.
We finally establish an equivalence between $\calP$ in the local picture
and $\calN$ in the particle picture.
The latter, in turn, has 
a description in terms of tensor products of the crystals $\Affz$ and 
$\Affo$.

The observations in this paper and in~\cite{HKMW98a} might be applied
and extended in various directions. Two of them are:
\begin{enumerate}
\item It is possible to derive
difference equations for correlation functions and form factors of the
alternating spin model using techniques analogous to those described
in~\cite{JM}. It should also be possible to evaluate these quantities
by making use of the free field realisation of $\uq$.

\item The approach should generalise in a straightforward manner to
alternating spin models with three or more different alternating spins.
\end{enumerate}




\vspace*{5mm}
\noindent
\texttt{jhong@math.snu.ac.kr, sjkang@math.snu.ac.kr,\\
miwa@kurims.kyoto-u.ac.jp, r.a.weston@ma.hw.ac.uk
}

\end{document}